\documentclass[11pt]{amsart}

\numberwithin{equation}{section} \oddsidemargin=-.0cm
\usepackage{framed}
\usepackage{graphicx}
\usepackage{amsfonts, amssymb, amsmath, graphics, subfigure,mathtools}
\usepackage[usenames,dvipsnames,table]{xcolor}
\usepackage{caption}
\usepackage{booktabs} 

\usepackage[colorlinks=true, pdfstartview=FitV, linkcolor=blue,
            citecolor=blue, urlcolor=blue]{hyperref}  

\numberwithin{equation}{section}
\definecolor{rred}{rgb}{0.7,0,0.1}
\definecolor{greenrb}{rgb}{0.2,0.6,0.2}

\newcommand{\mk}{\color{black}}
\newcommand{\mkr}{\color{black}}

\newcommand{\att}{\color{black}}

\newcommand{\hh}{\color{black}} 
\newcommand{\hl}{\color{black}} 
\newcommand{\HL}{\color{black}}

\def\bt{\begin{thm}}
\def\et{\end{thm}}
\def\bl{\begin{lem}}
\def\el{\end{lem}}
\def\bd{\begin{defi}}
\def\ed{\end{defi}}
\def\bc{\begin{cor}}
\def\ec{\end{cor}}
\def\bp{\begin{proof}}
\def\ep{\end{proof}}
\def\br{\begin{rem}}
\def\er{\end{rem}}

\def\bprop{\begin{prop}}
\def\eprop{\end{prop}}

\def\s{\mathfrak{s}}
\def\xii{\xi}
\def\c{\mathfrak{c}}

\def\ur{u_{R}}
\newcommand{\urc}[1]{u_{R,#1}} 

\def\bi{\begin{itemize}}
\def\ei{\end{itemize}}
\def\ben{\begin{enumerate}}
\def\een{\end{enumerate}}

\def\mathbi#1{\textbf{\em #1}}

\def\Forall{\text{ } \forall \:}
\def\d{\mathrm{d}}

\def\be{\begin{equation}}
\def\ee{\end{equation}}
\def\bes{\begin{equation*}}
\def\ees{\end{equation*}}
\def\bea{\begin{equation} \begin{aligned}}
\def\eea{\end{aligned} \end{equation}}
\def\beas{\begin{equation*} \begin{aligned}}
\def\eeas{\end{aligned} \end{equation*}}

\def\u{u}

\newtheorem{thm}{Theorem}[section]
\newtheorem{lem}{Lemma}[section]
\newtheorem{defi}{Definition}[section]

\newtheorem{prop}[thm]{Proposition}

\newtheorem{rem}{Remark}[section]
\newtheorem{cor}{Corollary}[section]

\graphicspath{{./figures/},{../Figures/}}

\title[Finite-Horizon  Parameterizing Manifolds and Suboptimal Control of Nonlinear PDEs]{Finite-Horizon Parameterizing Manifolds, and Applications to Suboptimal Control of Nonlinear Parabolic PDEs}

\author[Mickael Chekroun]{Micka\"el D. Chekroun}
\address[MC]{Department of Mathematics, University of Hawai`i at M$\overline{\mbox{a}}$noa, Honolulu, HI 96822, USA, and 
Department of Atmospheric \& Oceanic Sciences, University of California, Los Angeles, CA 90095-1565, USA} 
\email{mdchekroun@math.hawaii.edu}
\email{mchekroun@atmos.ucla.edu}
\author[Honghu Liu]{Honghu Liu}
\address[HL]{Department of Atmospheric \& Oceanic Sciences, University of California, Los Angeles, CA 90095-1565, USA}
\email{hliu@atmos.ucla.edu}


\begin{document}

\begin{abstract}

This article proposes a new approach for the design of low-dimensional suboptimal controllers to optimal control problems of nonlinear partial differential equations (PDEs) of parabolic type.  The approach fits into the long tradition of seeking for slaving relationships between  the small scales and the large ones (to be controlled) but differ by the introduction of a new type of manifolds to do so, namely the {\it finite-horizon parameterizing manifolds} (PMs). 
Given a finite horizon $[0,T]$ and a low-mode truncation of the PDE, a PM provides an approximate parameterization of the high modes by the controlled low ones so that  the unexplained high-mode energy is reduced \textemdash \, in a mean-square sense over $[0,T]$ \textemdash \,  when this parameterization is applied.

Analytic formulas of  such PMs are derived by application of the method of pullback approximation of the high-modes introduced in \cite{CLW13a}.
These formulas allow for an effective derivation of  reduced systems of ordinary differential equations (ODEs), aimed to model  the evolution of the low-mode truncation of the controlled state variable, where the high-mode part is approximated by the PM function applied to the low modes.  
The design of low-dimensional suboptimal controllers is then obtained by (indirect) techniques from 
finite-dimensional optimal control theory, applied to the PM-based reduced ODEs.

{\it A priori} error estimates  between the resulting PM-based low-dimensional suboptimal controller $u_R^\ast$ and the optimal controller $u^*$ are derived {\mkr under a second-order sufficient optimality condition}. 
These estimates demonstrate that the closeness of $u_R^\ast$   to $u^*$ is mainly conditioned on two factors: (i) the parameterization defect of a given PM, associated respectively with the  suboptimal controller $u_R^\ast$ and the optimal controller $u^*$;  and (ii) the energy kept in the high modes of the PDE solution either driven by  $u_R^\ast$ or $u^*$ itself. 

The practical performances of such PM-based suboptimal controllers are numerically assessed for optimal control problems associated with a Burgers-type equation; the locally as well as globally distributed cases being both considered.  The numerical results show that a PM-based reduced system allows for the design of suboptimal controllers with good performances provided that the associated parameterization defects and energy kept in the high modes are small enough, in agreement with the rigorous results.

\end{abstract}

\maketitle

\tableofcontents

\section{Introduction}

In this article, we propose a new approach for the synthesis  of {\mkr low-dimensional} suboptimal controllers for optimal control problems of nonlinear partial differential equations (PDEs) of parabolic type. Optimal control of PDEs has been extensively studied in the past few decades due largely to its broad applications in both engineering and various scientific disciplines, and fruitful results have been obtained; see {\it e.g.} the monographs \cite{BCD08,bensoussan2007representation,christofides2008control,Fattorini99,Fur00,HPUU09,Lions71,Tro10}. 

Due to the complexity of most applications, optimal control problems {\mk of parabolic PDEs} are often solved numerically.  Among the commonly used methods one finds methods that solve at once the associated optimality system using techniques  such as the Newton or quasi-Newton methods \cite{Betts10,HPUU09,ito2008lagrange},
or methods that use optimization algorithms {\mkr involving for instance} an approximation to the gradient of the cost functional; see {\it e.g.} \cite{Betts98survey,HPUU09,ito2008lagrange,Tro10}. In this case, the gradient can be approximated by using sensitivity methods or methods based on the adjoint equation;  see {\it e.g.}\cite{AT90,bewley2001dns,BTZ00,gunzburger1999sensitivities,gunzburger2000adjoint,ito1998optimal,Med04,medjo2008optimal}. Efficient (and accurate) solutions can be designed by such methods \cite{AT90,baker2000,bewley2001dns,Choi_al93,HK00,Med04,medjo2008optimal} which may lead however to high-dimensional problems that can turn out to be computationally expensive to solve, especially for fluid flows applications.  The task becomes even  more challenging  when  a dynamic programming approach is adopted, involving 
typically to solve (infinite-dimensional) Hamilton-Jacobi-Bellman (HJB) equations \cite{BCD08,Beeler_al00,cannarsa1996infinite,crandall1992user,da2000dynamic,da2000dynamic2,da2002second}.

{\mk As an alternative,} various reduction techniques have been {\mk proposed  in the literature to} seek {\mk instead} for {\mk low-dimensional suboptimal controllers}. {\mk The main issue  related to such techniques relies however on the {\mkr ability to} design suboptimal solutions close enough to the genuine optimal one \cite{Ded10,GK11,Hinze_al05,IK08,TV09}, while keeping  cheap enough  the numerical efforts to do so.}  A general class of model reduction techniques used extensively in this context is the so-called reduced-order modeling (ROM) approach, based on approximating the nonlinear dynamics by a Galerkin technique relying on  basis functions, possibly empirical \cite{Franke12,HK98,HK00,ravindran2000reduced}. Various ROM techniques differ in the choice of the basis functions. One popular method that falls into this category is the so-called proper orthogonal decomposition (POD); see among many others \cite{AK01,bergmann2008optimal,Hinze_al05,Holmes_al12,KV99,KV02,LT01,ravindran2002adaptive}, and \cite{GK11,ito2001reduced,IS01} for other methods in constructing the reduced basis.  {\mk We refer also to  \cite{kunisch2004hjb} for suboptimal controllers designed from the solutions of low-dimensional HJB equations associated with POD-based Galerkin reduced-order models. }


Such Galerkin/ROM-based techniques can lead to a synthesis of very efficient suboptimal controllers once, at a given truncation, the disregarded high-modes do not contribute significantly to the dynamics of the low modes. However, when this is not the case, the seeking of parameterizations of the disregarded modes in terms of the low ones becomes central {\mkr for} the design of surrogate low-dimensional models {\mkr of good performances}.  The idea of seeking for slaving relationships between the unstable or least stable modes with the more stable ones has a long tradition in control theory of large-dimensional or distributed-parameter systems. {\mk For instance}, by use of methods from singular perturbation theory, the authors in \cite{kokotovic1968singular,kokotovic1976singular,kokotovic1984applications,kokotovic99} investigated  
the construction of such slaving functions for slow-fast systems in terms of invariant (slow) manifolds.\footnote{See also \cite[Chap.~5]{lions1972some}  and \cite{Lions73} for the use of singular perturbation techniques for optimal control of PDEs.} Such manifolds are then used to decouple the slow and fast parts of the dynamics and to feed
back the slow component of the state only. This is especially important since the fast components of the state are in general difficult to measure/estimate and consequently to feedback.  

Complementary to singular perturbation methods, the authors of \cite{chen1992accelerated}  used tools of center manifold and normal
form theory to design a nonlinear controller and obtained a 
closed-loop center manifold for a truncated distributed-parameter system; in their case
proximity to a bifurcation constitutes a guarantee to the separation of relevant
time scales of the problem.  In \cite{christofides1996nonlinear,CPD97}, the authors have gone beyond the finite-dimensional singular
perturbation work of \cite{kokotovic99} or center-manifold-based work of \cite{chen1992accelerated} to exploit approximate inertial manifolds (AIMs) \cite{FMT88} in the infinite-dimensional case; the latter are global manifolds in phase space that can be thought of as generalizations of slow/center manifolds.
Using AIMs, the authors of  \cite{chen1992accelerated} designed then observer-based nonlinear feedback controllers (through the corresponding
closed-loop AIMs) and demonstrated their performance.

 The potential usefulness of inertial manifolds (IMs) \cite{CFNT89,FST88,temam1990inertial} or  AIMs in control theory of {\mk nonlinear parabolic PDEs} was actually quickly {\mk identified} after  IM theory started to be established \cite{Bru91,christofides1996nonlinear,CPD97,SK95}; see {\it e.g.} \cite{RT97,SK98} for a state-of-the-art of the literature at the end of the $90$s. However since these works, IMs or AIMs have been mainly employed to derive low-dimensional vector fields for the design of feedback controllers \cite{AC00,Rosa03}. To the exception of \cite{AC02,IK08}, the use of IMs or AIMs to design suboptimal solutions to optimal control problems have been much less considered.

The main purpose of this article is to introduce  a general {\mk framework \textemdash\, in the continuity but different from the AIM approach \textemdash\,} for the {\mk effective} derivation of {\it suboptimal low-dimensional}  solutions to {\mk optimal control problems associated with nonlinear PDE such as \eqref{PDE1intro} given below}.    To be more specific, given an ambient Hilbert space, $\mathcal{H}$,  the control problems of PDEs we will consider hereafter take the following abstract form:
\begin{equation}  \label{PDE1intro}
\frac{\d y}{\d t} = L y + F(y) +  \mathfrak{C} u(t),  \qquad t \in (0, T],
\end{equation}
where $L$ denotes a linear operator, $F$ some nonlinearity, and  
$\mathfrak{C}$  denotes a {\it bounded linear operator} on $\mathcal{H}$; the state variable $y$ and the controller $u$  living  both in  $L^2(0,T; \mathcal{H})$ for a given horizon $T>0$;  see  Section \ref{Sect_functfram} for more details.

The underlying idea consists  of seeking for manifolds $\mathfrak{M}$ aimed to provide \textemdash\, over a finite horizon $[0,T]$ \textemdash\, an {\it approximate parameterization} of the small scales of the solutions to the {\it uncontrolled PDE} associated with Eq.~\eqref{PDE1intro}, namely 
\begin{equation}   \label{SEE-null} 
\frac{\d y}{\d t} = L y + F(y), 
\end{equation}
in terms of their {\it large scales}, so that $\mathfrak{M}$ allows in turn to derive low-dimensional reduced models from which suboptimal controllers 
can be  efficiently designed by standard methods of finite-dimensional optimal control theory {\mkr such as found in {\it e.g.}}   \cite{bonnard2003singular,Bryson_al75,Kirk12,Kno81,SL12}.  In that respect,  the notion of {\it finite-horizon parameterizing manifold} (PM)  is introduced in Definition \ref{def:PM} below.  Finite-horizon PMs distinguish from the {\mkr more} classical AIMs in the sense that they provide approximate parametrization of the small scales by the large ones  in {\HL the $L^2$-sense} (over $[0,T]$) rather than a hard $\epsilon$-approximation to be valid for each time $t\in [0,T]$, cf. \cite{FMT88}.  In particular, a finite-horizon PM allows to reduce the (cumulative) unexplained high-mode energy (over $[0,T]$) from the low modes to be controlled, in a way different from other slaving relationships considered so far; the  high-mode energy being reduced in a {\it mean-square sense} in the case of finite-horizon PMs.   

Obviously, the difficulty relies still on the ability of such an approach to give access to suboptimal controllers of good performance. A priori the task in not easy and  a key feature to ensure that a ``good'' performance is achieved from such a suboptimal {\mk low-dimensional} controller, $u_R^*$, relies on the ability of the manifold $\mathfrak{M}$  derived from the uncontrolled problem to still achieve a sufficiently ``small'' {\it parameterization defect} (over the horizon $[0,T]$) of the small scales by the large ones  once a controller $u_R^*$ is used to drive the PDE \eqref{PDE1intro}; see \eqref{Eq_QualPM} in Definition \ref{def:PM}. This point is rigorously formulated as  Theorem \ref{THM_controller} {\mk in Section \ref{Sect_reduced model}} (see also Corollary \ref{Cor_2}), which provides {\mkr \textemdash\, under a second-order sufficient optimality condition \textemdash\,} error estimates on how ``close'' a low-dimensional suboptimal controller {\mk $u_R^\ast$}, designed  from a PM-based reduced system, is to the optimal controller $u^*$.  The error estimates \eqref{thm1:goal} and \eqref{cor2:goal} show in particular that the closeness of {\mk $u_R^\ast$}  to $u^*$ is mainly conditioned on two factors: (i) the parameterization defect of a given PM, associated respectively with the 
suboptimal controller $u_R^\ast$ and the optimal controller $u^*$;  and (ii) the energy kept in the high modes of the PDE solution {\mk either driven by  $u_R^\ast$ or $u^*$ itself.}

The article is organized as follows. The functional framework associated with optimal control problems related to \eqref{PDE1intro} is introduced in Section~\ref{Sect_functfram}. The definition of finite-horizon PMs and a practical procedure to get access to such PMs are introduced in Section~\ref{Sect_PM}. In particular analytic formulas of leading-order PMs are provided; the latter being subject to a cross {\it non-resonance condition} \eqref{NR} to be satisfied between the high and the low modes; see Section \ref{ss:h1}. Section~\ref{Sect_reduced model} is devoted, given an arbitrary  PM,  to the derivation of rigorous {\it a priori} error estimates between a low-dimensional PM-based  suboptimal controller and the optimal one; see Theorem \ref{THM_controller} and Corollary \ref{Cor_2}. 
The performance of the resulting PM-based reduction approach is numerically investigated on a Burgers-type equation in the context of globally and locally distributed control  laws; see Sections~\ref{Sect_Burgers}--\ref{Sect_Burgers_h2}, and Section~\ref{Sect_Burgers_local}. As a main byproduct, the numerical results strongly indicate that a PM-based reduced system allows for a design of suboptimal controllers with good performances provided that the  aforementioned parameterization defects and the energy contained in the high modes are small enough, in agreement with the theoretical predictions of Theorem \ref{THM_controller} and Corollary \ref{Cor_2}.  This is particularly  demonstrated  in Section \ref{Sect_Burgers_h2}, where analytic formulas derived in Theorem \ref{THM_h2} give access to higher-order PMs  with reduced parameterization defects compared to those of the leading-order PMs introduced in Section \ref{Sect_PM}.  In all the cases, the analytic formulas of the PMs used hereafter allows for an efficient design of suboptimal controllers by standard (and {\mkr simple}) application of  the  Pontryagin maximum principle \cite{bonnard2003singular,Bonnard_al06,Kirk12,PBGM64} to the PM-based reduced systems.

\section{Optimal Control of Nonlinear PDEs, and Functional Framework}\label{Sect_functfram}
The functional framework for the optimal control problem considered in this article takes place in Hilbert spaces.  Let us first introduce the {\mk class of} partial differential equations (PDEs) to be controlled. For a given Hilbert space $\mathcal{H}$,  we consider $\mathcal{H}_1$ to be a subspace compactly and densely embedded in $\mathcal{H}$ such that $A:\mathcal{H}_1\rightarrow \mathcal{H}$ is a sectorial operator \cite[Def.~1.3.1]{Hen81} satisfying
\bes
-A \mbox{ is stable in the sense that its spectrum satisfies }  \mathrm{Re} (\sigma(-A) )< 0.  
\ees
To include in our framework PDEs for which the nonlinear terms are responsible of a loss of regularity compared to the ambient space $\mathcal{H}$, we consider standard interpolated spaces $\mathcal{H}_\alpha$ between $\mathcal{H}_1$ and $\mathcal{H}$ (with $\alpha \in [0, 1))$\footnote{depending on the problem at hand; see {\it e.g.} \cite{Hen81}.} along with  perturbations of the linear operator $-A$ given by a one-parameter family, $\{B_{\lambda}\}_{\lambda \in \mathbb{R}}$, of bounded linear operators from  $\mathcal{H}_\alpha$ to  $\mathcal{H}$,  that depend continuously on  a real parameter $\lambda$.

By defining 
\bes
L_{\lambda}:=-A+B_{\lambda},
\ees
we are thus left  with a one-parameter family of sectorial operators $\{-L_\lambda\}_{\lambda \in \mathbb{R}}$, each of them mapping  $\mathcal{H}_1$ into $\mathcal{H}$. Finally, $F: \mathcal{H}_\alpha \rightarrow \mathcal{H}$ will denote a {\it continuous $k$-linear mapping} ($k\geq 2$) for some $\alpha \in [0,1)$.\footnote{In particular, nonlinearities including a loss of regularity compared to the ambient space $\mathcal{H}$, are allowed; see {\it e.g.} Section \ref{Sect_Burgers} below.} 

The nonlinear evolution equation to be controlled takes then the following abstract form:
\begin{equation}  \label{PDE1}
\frac{\d y}{\d t} = L_\lambda y + F(y) +  \mathfrak{C} u(t),  \qquad t \in (0, T],
\end{equation}
where $y \in  L^2(0,T; \mathcal{H})$ denotes the state variable, $u \in  L^2(0,T; \mathcal{H})$ denotes the controller; $T > 0$ being  a fixed horizon, and 
\be
\mathfrak{C}: \mathcal{H} \rightarrow \mathcal{H}
\ee
denoting a {\it bounded {\mk (and non-zero) linear control} operator}. In particular, we will be mainly concerned with distributed control problems (control inside the domain) and not with problems involving a control on the boundary which leads typically to an unbounded control operator; see {\it e.g.} \cite[Part V, Chap. 2 and 3]{bensoussan2007representation} and \cite{fattorini1968,Fattorini99,flandoli84}. The parameter $\lambda$ governs typically the presence of (linearly) unstable modes for \eqref{PDE1}. In the application considered in Sections \ref{Sect_Burgers}--\ref{Sect_Burgers_local}, it will be chosen so that the linear operator, $L_{\lambda}$, admits  {\mk large-scale unstable modes}.

We introduce next the cost functional $J: L^2(0,T; \mathcal{H}) \times L^2(0,T; \mathcal{H}) \rightarrow \mathbb{R}$ given by
\be  \label{J_intro}
J(y,u) := \int_0^T [\mathcal{G}(y(t)) + \mathcal{E}(u(t)) ] \d t,
\ee 
where $\mathcal{G}: \mathcal{H} \rightarrow \mathbb{R}^+$ and $\mathcal{E}: \mathcal{H} \rightarrow \mathbb{R}^+$ are assumed to be continuous, and to satisfy the following conditions:
\be\label{C1}\tag{C1}
\mathcal{G}  \mbox{ is  uniformly Lipschitz on bounded sets of } \mathcal{H},
\ee
and
\be\label{C2}\tag{C2}
\| u\|\leq \|v\| \Longrightarrow \mathcal{E}(u)\leq \mathcal{E}(v),
\ee
where $\|\cdot \|$ denotes the $\mathcal{H}$-norm.

{\mk Given such a cost functional,\footnote{We refer to Sections \ref{Sect_Burgers}--\ref{Sect_Burgers_local} for other type of cost functional including a terminal cost.} we will} consider in this article {\mk the following type of optimal control problem}: 
\be  \label{P1}  \tag {$\mathcal{P}$}
\begin{aligned}
 \min \, J(y, u)  \quad \text{ s.t. } \quad &  (y, u) \in L^2(0,T; \mathcal{H}) \times  L^2(0,T; \mathcal{H}) \text{ solves  Eq.} ~\eqref{PDE1} \\
&\text{ subject to }  \qquad   y(0)  = y_0 \in \mathcal{H}.
\end{aligned}
\ee

To simplify the presentation, we will make the following assumptions on $L_\lambda$ and $F$ throughout this article:

\vspace{1ex}
{\bf Standing Hypothesis.}  {\it $L_\lambda$ is self-adjoint, whose eigenvalues {\mk (arranged in descending order)} are denoted by $\{\beta_i(\lambda)\}_{i \in \mathbb{N}}$;  and the eigenvectors $\{e_i(\lambda)\}_{i \in \mathbb{N}}$ of $L_\lambda$ form a Hilbert basis of $\mathcal{H}$.  {\HL The eigenvectors are regular enough such that $e_i(\lambda) \in \mathcal{H}_\alpha$ for all $i\in \mathbb{N}$.} The nonlinearity $F: \mathcal{H}_\alpha \rightarrow \mathcal{H}$ is a continuous $k$-linear mapping for some $k \ge 2$, and for some $\alpha \in [0,1)$. {\HL In particular, $F(0) = 0$.}}

\vspace{1ex}

We also assume that for any initial datum $y_0\in \mathcal{H}$, any $T>0$,  and any given $u \in L^2(0, T; \mathcal{H})$, the Cauchy problem 
\begin{equation}   \label{SEE} 
\frac{\d y}{\d t} = L_\lambda y + F(y) +  \mathfrak{C} u(t),  \quad y(0)  = y_0 \in \mathcal{H},
\end{equation}
has a unique solution $y(\cdot,y_0;u) \in C([0,T]; \mathcal{H}) \cap L^2(0,T; \mathcal{H}_\alpha)$, which lives furthermore in the space $C^1((0,T]; \mathcal{H}) \cap C([0,T]; \mathcal{H}_\alpha) \cap L^2(0,T; \mathcal{H}_1)$ {\mk when} $y_0 \in \mathcal{H}_\alpha$; 
see {\it e.g.}  \cite[Chap.~3]{Hen81} and \cite[Chap.~7]{Lun95} for conditions under which such properties are guaranteed. {\mk Section \ref{Sect_TP_existence} below deals with such an example.}

\section{Finite-Horizon Parameterizing Manifolds: Definition, Pullback Characterization and Analytic Formulas} \label{Sect_PM}

This section is devoted to the definition of finite-horizon parameterizing manifolds (PMs) for a given PDE of type \eqref{SEE} and a general {\mk method} to give access to {\mk explicit formulas of} such finite-horizon PMs in practice through pullback limits associated with certain backward-forward systems built from the uncontrolled Eq.~\eqref{SEE-null}. 

The key idea {\mk takes its roots in the notion of} (asymptotic) parameterizing manifold {\mk introduced} in \cite{CLW13a}\footnote{{\mk mainly in a stochastic context; see however \cite[Section 8.5]{CLW13a} for  the deterministic setting.}}, which {\mkr reduces} here of  {\mkr approximating} --- over {\mk some prescribed} finite time interval $[0,  T]$ --- the modes with ``high'' wave numbers  as a pullback limit depending on the time-history of {\mkr (some approximation of)} the {\mkr dynamics of the} modes with ``low'' wave numbers. The cut between what is ``low'' {\mk and} what is ``high'' is organized {\mk in an abstract setting} as follows; {\mk we refer to Section \ref{Sect_Burgers_local} for a more concrete specification of such a cut in the case of locally distributed controls. 
The subspace $\mathcal{H}^{\c} \subset \mathcal{H}$ defined by, 
\bea \label{Hc}
\mathcal{H}^{\c} := \mathrm{span}\{e_1, \cdots, e_m\}, 
\eea 
spanned by the $m$-leading modes will be considered as our subspace associated with the low modes. 
Its topological complements, $\mathcal{H}^{\s}$ and $\mathcal{H}^{\s}_\alpha$, in  respectively $\mathcal{H}$ and $\mathcal{H}_\alpha$, will be considered as associated with the high modes,  leading to the following decomposition
\be \label{Hs}
\mathcal{H} = \mathcal{H}^{\c} \oplus \mathcal{H}^{\s}, \qquad \mathcal{H}_\alpha = \mathcal{H}^{\c} \oplus \mathcal{H}^{\s}_\alpha.
\ee}
We will use $P_{\c}$ and $P_{\s}$ to denote the canonical projectors associated with $\mathcal{H}^{\c}$ and $\mathcal{H}^{\s}$, respectively. Here, the usage of the eigenbasis in the decomposition of the phase space is employed for the sake of analytic formulations derived hereafter. 
In practice, the methodology presented below can be (numerically) adapted when the phase space $\mathcal{H}$ is decomposed by using other bases; see also Remark~\ref{Rmk_POD} (ii).

\subsection{ Finite-horizon parameterizing manifolds} 

Let $t^\ast>0$ be fixed, $\mathcal{V}$ be an open set in $\mathcal{H}_\alpha$, and $\mathcal{U}$ an open set in $L^2(0,t^\ast; \mathcal{H})$.  For a given PDE of type \eqref{SEE},
a {\it finite-horizon parameterizing manifold}  $\mathfrak{M}$ over the interval $[0, t^\ast]$ is defined as the graph of a function $h^{\mathrm{pm}}$ from  $\mathcal{H}^{\c}$ to $\mathcal{H}^{\s}_\alpha$, which is aimed to provide, for any $y(t, y_0; u)$ solution of  \eqref{SEE} with initial datum $y_0 \in \mathcal{V}$ and control $u \in \mathcal{U}$,  an approximate parameterization  of {\mk its ``high-frequency''}  part, $y_{\s}(t, y_0;u)=P_{\s} y(t, y_0;u)$,  in terms of  {\mk its ``low-frequency'' part}, $y_{\c}(t, y_0;u)=P_{\c} y(t, y_0; u)$,  so that  the mean-square error, $\int_0^{t^\ast} \bigl \|y_{\s}(t, y_0; u) -  h^{\mathrm{pm}}(y_{\c}(t, y_0; u)) \bigr \|_\alpha^2 \, \d t $, is {\mk strictly} smaller than the {\mkr high-mode energy}
of $y_{\s}$,  $\int_0^{t^\ast}  \|y_{\s}(t, y_0; u)\|_\alpha^2 \; \d t$. 
{\mk Here the frequencies are understood in a spatial sense, {\it i.e.} in terms of wave numbers\footnote{In particular, the reduction techniques developed in this article should not be confused with the reduction techniques based on the {\it slow manifold theory} which have been used to deal with the reduction of optimal control problems arising in {\it slow-fast systems}, where the separation of the dynamics holds in time rather than in space; see {\it e.g.} \cite{kokotovic99,lebiedz2013numerical,motte2000slow}. Furthermore, unlike slow manifolds, the finite-horizon PMs considered in this article are not invariant for the dynamics. To the contrary, they correspond to manifolds for which the dynamics {\mk wanders} around, within {\mk some} margin whose size (in a mean square sense) is {\mk strictly smaller than} the {\mkr energy} unexplained by the $\mathcal{H}^{\c}$-modes.}.} In statistical terms, a finite-horizon PM function $h^{\mathrm{pm}}$ {\mk can thus be thought of as a slaving relationship between the high modes and the low ones such}
that the fraction of {\mkr energy}\footnote{over the time interval $[0, t^\ast]$.} of $y_{\s}$  unexplained by $h^{\mathrm{pm}}(y_{\c})$ ({\mk {\it i.e.} {\mkr {\it via}} this slaving relationship}) is less than {\mk unity}.

{\mk In more precise terms, we are left with the following definition:}

\bd  \label{def:PM}

Let $t^\ast>0$ be fixed, $\mathcal{V}$ be an open set in $\mathcal{H}_\alpha$, and $\mathcal{U}$ an open set in $L^2(0,t^\ast; \mathcal{H})$.  A manifold $\mathfrak{M}$ of the form
\bea \label{eq:PM_def}
\mathfrak{M} := \{\xii + h^{\mathrm{pm}}(\xi) \mid \xi \in \mathcal{H}^{\c}\}
\eea
is called a {\it finite-horizon parameterizing manifold (PM)} over the time interval $[0, t^\ast]$ associated with the PDE \eqref{SEE} if the following conditions are satisfied:

\bi
\item[(i)] The function $h^{\mathrm{pm}}: \mathcal{H}^{\c} \rightarrow \mathcal{H}^{\s}_\alpha$ is continuous. 
\item[(ii)] The following inequality holds for any $y_0 \in \mathcal{V}$ and any $u \in \mathcal{U}$:
\bea \label{PM condition}
\int_0^{t^\ast}  \bigl \|y_{\s}(t,y_0; u) -  h^{\mathrm{pm}}(y_{\c}(t, y_0; u)) \bigr \|_\alpha^2 \, \d t < \int_0^{t^\ast}  \|y_{\s}(t, y_0; u)\|_\alpha^2 \, \d t,
\eea
where $y_{\c}(\cdot, y_0; u)$ and $y_{\s}(\cdot, y_0; u)$ are the projections to respectively the subspaces $\mathcal{H}^{\c}$ and $\mathcal{H}^{\s}_\alpha$ of the solution $y(\cdot, y_0; u)$ for the PDE ~\eqref{SEE} driven by $u$ emanating from $y_0$.  

\ei

For a given initial datum $y_0$, if $y_{\s}(\cdot, y_0; u)$ is not identically zero, the {\mk {\it parameterization defect} of $\mathfrak{M}$  over $[0,t^*]$, and 
associated with the control $u$,} is defined as the following ratio: 
\be\label{Eq_QualPM}
\boxed{
Q(t^\ast, y_0; u) :=  \frac{\int_0^{t^\ast}  \bigl \|y_{\s}(t, y_0; u) - h^{\mathrm{pm}}(y_{\c}(t, y_0; u)) \bigr \|_\alpha^2 \, \d t}{ \int_0^{t^\ast} \|y_{\s}(t, y_0; u)\|_\alpha^2 \, \d t}.}
\ee

\ed

{\mk Note that in Sections~\ref{Sect_Burgers}, \ref{Sect_Burgers_h2} and \ref{Sect_Burgers_local}, we will illustrate numerically that finite-horizon PMs can actually be obtained from the uncontrolled PDE \eqref{SEE-null}, with still possibly small parameterization defects when a controller $u$ is  applied. The procedure to build in practice such PMs from the uncontrolled PDE \eqref{SEE-null} is described in the next section; see also \cite[Section 8.5]{CLW13a} for the construction of PMs over arbitrarily (and sufficiently) large horizons.}



%
%
%
%
%

\subsection{Finite-horizon parameterizing manifolds as  pullback limits of  backward-forward systems: the leading-order case}  \label{ss:h1}
We consider now the important problem of the practical determination of finite-horizon PMs for PDEs of type \eqref{SEE}. 
As mentioned above, following  \cite{CLW13a},  the pullback {\mkr approximation} {\mk of the high modes in terms of the low ones {\it via} appropriate} auxiliary {\mk systems associated with} the {\it uncontrolled} PDE \eqref{SEE-null} will constitute  the key ingredient  to propose a solution to this problem; see also \cite[Section 8.5]{CLW13a}. In that respect, we consider first the following {\it backward-forward system} associated with the uncontrolled PDE \eqref{SEE-null}:
\begin{subequations}\label{LLL}
\begin{align}
& \frac{\mathrm{d} y^{(1)}_{\c}}{\mathrm{d} s} =  L_\lambda^{\c} y^{(1)}_{\c}, 
 && s \in [ -\tau, 0],   \quad \; y^{(1)}_{\c}(s)\vert_{s=0} = \xii,  \label{LLL1} \\
& \frac{\mathrm{d} y^{(1)}_{\s}}{\mathrm{d} s} =  L_\lambda^{\s} y^{(1)}_{\s}  +  P_{\s} F(y^{(1)}_{\c}), 
&&  s \in [-\tau, 0],  \qquad y^{(1)}_{\s}(s) \vert_{s=-\tau}= 0, \label{LLL2}
\end{align}
\end{subequations}
{\mk  where $L_\lambda^{\c} := P_{\c} L_\lambda$, $L_\lambda^{\s} := P_{\s} L_\lambda$, and $\xi \in \mathcal{H}^{\c}$.  We refer to Section \ref{Sect_Burgers_h2} for other {\it backward-forward systems} used in the construction of {\mk higher-order} finite-horizon PMs. }

In the system above, the  initial value of $y_{\c}^{(1)}$ is prescribed at $s = 0$, and the initial value of $y_{\s}^{(1)}$ at $s= - \tau$. {\mk The solution of this system is obtained by using a two-step {\it backward-forward integration procedure} \textemdash\, where Eq.~\eqref{LLL1}  is integrated first backward and Eq.~\eqref{LLL2} is then integrated forward  \textemdash\,}  made possible due to the partial coupling  present in \eqref{LLL} where $y^{(1)}_{\c}$ forces the evolution equation of $y^{(1)}_{\s}$ but not reciprocally. Due to {\mk this} forcing introduced by  $y_{\c}^{(1)}$ which emanates (backward) from $\xi$, the solution process $y_{\s}^{(1)}$ depends naturally on $\xi$. For that reason, we will emphasize this dependence as $y_{\s}^{(1)}[\xi]$ hereafter. 

It is clear that the solution to the above system is given by:
\bea
y^{(1)}_{\c}(s) & = e^{s L_\lambda^{\c}}\xi, \hspace{11em}   s \in [-\tau, 0], \; \xi \in \mathcal{H}^{\c}, \\
y_{\s}^{(1)}[\xi]{\mk(-\tau, s)} & = \int_{-\tau}^s e^{(s-\tau') L_\lambda^{\s}} P_{\s} F(e^{\tau' L_\lambda^{\c}}\xi) \d \tau', \hspace{2em}   s \in [-\tau, 0].
\eea

{\mk The dependence in $\tau$ and $s$ in $y_{\s}^{(1)}[\xi]$ is made apparent to emphasize the two-time description  employed for the description of the non-autonomous dynamics  {\mkr inherent to} \eqref{LLL2}; see {\it e.g.}~\cite{CLR13,CSG11}. Adopting the language of non-autonomous dynamical systems \cite{CLR13,CSG11}, we then define $h^{(1)}_\lambda(\xi)$ as the  following} {\it pullback limit} of the $y_{\s}^{(1)}$-component of the solution to the above system, {\it i.e.},
\begin{equation} \label{det_h1}
\boxed{{\mk h^{(1)}_\lambda(\xi)  :=\lim_{\tau \rightarrow +\infty} y^{(1)}_{\s}[\xi](-\tau, 0)} =  \int_{-\infty}^0 e^{-\tau' L^{\s}_\lambda} P_{\s} F(
e^{\tau' L^{\c}_\lambda}\xi) \,\mathrm{d} \tau', \quad \Forall \xi \in  \mathcal{H}^{\c},}
\end{equation}
when the latter limit exists.  We derive hereafter necessary and sufficient conditions for such a limit  to exist. 

In that respect, first note that since $L_\lambda$ is  self-adjoint, we have 
\be \label{uc expansion}
e^{\tau' L^{\c}_\lambda}\xii = \sum_{i = 1}^m e^{\tau' \beta_i(\lambda)} \xi_i e_i,   
\ee
where $\xi = \langle \xi, e_i \rangle$, $i \in \mathcal{I} :=\{1, \cdots, m\}$ with $m=\mathrm{dim}(\mathcal{H}^{\c})$, and $\langle \cdot, \cdot \rangle$ denoting  the inner-product in the ambient Hilbert space $\mathcal{H}$. 

Now for a fixed $\tau>0$, by projecting $y_{\s}^{(1)}[\xi](-\tau,0)$ against each eigenmode $e_n$ for $n > m$, we obtain, by using \eqref{uc expansion} and the $k$-linear property of $F$, 
\bea \label{h1 expansion}
y_{\s}^{(1)}[\xi](-\tau,0)
& = \sum_{n > m} \int_{-\tau}^0 e^{-\tau' \beta_n(\lambda)} \Bigl \langle F \Bigl(
\sum_{i = 1}^m e^{\tau' \beta_i(\lambda)} \xi_i e_i \Bigr), e_n \Bigr \rangle \,\mathrm{d}\tau' \, e_n \\ 
& = \sum_{n > m} \sum_{(i_1, \cdots, i_k) \in \mathcal{I}^k}  \int_{-\tau}^0 e^{- \beta_n(\lambda)\tau'  + \bigl( \sum_{j = 1}^k \beta_{i_j}(\lambda)\bigr) \tau'}  \,\mathrm{d}\tau'  \Bigl \langle  F(e_{i_1}, \cdots,  e_{i_k}), e_n \Bigr \rangle  e_n. 
\eea

From this identity, we infer that $h^{(1)}_\lambda$ is well defined if and only if each integral 
\bes
\int_{-\infty}^0 e^{- \beta_n(\lambda) \tau' + \bigl ( \sum_{j = 1}^k \beta_{i_j}(\lambda)\bigr) \tau'}  \,\mathrm{d}\tau'
\ees 
converges,  whenever the corresponding nonlinear interaction $F(e_{i_1}, \cdots,  e_{i_k})$ as projected against $e_n$, is non-zero. Namely, $h^{(1)}_\lambda$ exists if and only if the following (weak) {\it non-resonance condition} holds:
\begin{equation}  \label{NR} \tag{NR}
\begin{aligned} 
& \Forall \, (i_1, \cdots, i_k) \in \mathcal{I}^k,  \ n > m, \text{ it holds that} \\
 & \Bigl (\langle F(e_{i_1}, \cdots,  e_{i_k}), e_n \rangle \neq 0 \Bigr) \Longrightarrow \biggl ( \sum_{j = 1}^k \beta_{i_j}(\lambda) - \beta_n(\lambda) > 0 \biggr);
\end{aligned}
\end{equation}
see also \cite[Sect.~7]{CLW13a}.

Assuming the above \eqref{NR}-condition, it follows then from \eqref{det_h1} and \eqref{h1 expansion} that $h^{(1)}_\lambda$ takes the following form:
\be  \label{h1}
\boxed{h^{(1)}_\lambda(\xi)  = \sum_{n > m} \sum_{(i_1, \cdots, i_k) \in \mathcal{I}^k}  \frac{\xi_{i_1}\cdots \xi_{i_k}}{\sum_{j = 1}^k \beta_{i_j}(\lambda) - \beta_n(\lambda)} \Bigl \langle  F(e_{i_1}, \cdots,  e_{i_k}), e_n \Bigr \rangle  e_n.}
\ee

In particular under the \eqref{NR}-condition, each $e_n$-component of $h^{(1)}_\lambda(\xi)$ is \textemdash\, in the $\xi$-variable  \textemdash\, an homogeneous polynomial of order $k$, the order of the nonlinearity $F$. For that reason, $h^{(1)}_\lambda$ will be referred to as the {\it leading-order} finite-horizon PM when appropriate, that is when the latter provides a finite-horizon PM. We clarify in the remaining of this section, some (idealistic) conditions under which such a property is met by the manifold function $h^{(1)}_\lambda$ for the PDE \eqref{SEE}. In practice these conditions can be violated, while  the manifold function $h^{(1)}_\lambda$ defined by \eqref{h1} still constitutes a finite-horizon PM; see {\HL Sections \ref{Sec_numresults} and \ref{Sect_Burgers_local}} for numerical illustrations. 


To delineate conditions under which $h^{(1)}_\lambda$ is a finite-horizon PM is still valuable for the theory. This is the purpose of Lemma \ref{Lem:h1_FPM} below which relies on another key property of  $h^{(1)}_\lambda$ such as defined by \eqref{det_h1},  that can be explained using the language of invariant manifold theory for PDEs \cite{CLW13a,MW14}. The latter states that the manifold function $h^{(1)}_\lambda$ constitutes  \textemdash\, for the uncontrolled PDE \eqref{SEE-null}  \textemdash\,  the leading-order approximation of some local invariant manifold near the trivial steady state; {\HL see \cite[Appendix A]{MW14} and \cite[Sect. 7]{CLW13a}}.
Based on this result we formulate the following lemma about the existence of finite-horizon PMs.

\bl \label{Lem:h1_FPM}

Let $\lambda$ be fixed and $\mathcal{H}^{\c}$ be the subspace spanned by the first $m$ eigenmodes of the linear operator $L_\lambda$. Assume that the standing hypothesis of Section ~\ref{Sect_functfram} holds, and that
\be \label{eigen_relation}
\beta_m(\lambda) > 2k \beta_{m+1}(\lambda).
\ee
Assume furthermore that the non-resonance condition \eqref{NR} holds so that the pullback limit $h^{(1)}_\lambda$ defined by \eqref{det_h1} exists.

Assume that $h^{(1)}_\lambda$ is non-degenerate in the sense that there exists $C > 0$ such that
\be \label{non-degeneracy}
\|h^{(1)}_\lambda(\xi)\|_\alpha \ge C \|\xi\|_\alpha^k, \qquad \xi \in \mathcal{H}^{\c}.
\ee

Then, for any fixed $t^\ast>0$, there exist open neighborhoods $\mathcal{V} \subset \mathcal{H}^{\s}_\alpha$ and $\mathcal{U} \subset L^2(0, t^\ast; \mathcal{H})$ containing the origins of the respective spaces, such that $h^{(1)}_\lambda$ is a finite-horizon parameterizing manifold over the time interval $[0, t^\ast]$ for the PDE \eqref{SEE} driven by any {\mkr control} $u \in \mathcal{U}$ and with initial data taken from $\mathcal{V}$.

\el

\bp

Let us first recall some related elements from \cite{CLW13a}. Note that the PDE \eqref{SEE-null} fits into the framework of \cite[Cor.~7.1]{CLW13a}.\footnote{Eq.~\eqref{SEE-null} corresponds to a deterministic situation dealt with in \cite{CLW13a} by setting the noise amplitude to zero.} Since the nonlinearity $F$ is assumed to be $k$-linear for some $k\ge 2$, according to \cite[Cor.~7.1]{CLW13a}, under the assumption \eqref{eigen_relation}, there {\mk exists} of a local invariant manifold associated with the PDE \eqref{SEE-null} of the form,
\be
\mathfrak{M}_{\lambda}^{\mathrm{loc}} := \{\xi + h_{\lambda}^{\mathrm{loc}} (\xi) \mid \xi \in \mathfrak{B}\},
\ee
where $h_{\lambda}^{\mathrm{loc}}: \mathcal{H}^{\c} \rightarrow \mathcal{H}^\s_\alpha$ is the corresponding {\mk local manifold function}, $\mathfrak{B} \subset \mathcal{H}^{\c}$ is an open neighborhood of the {\mk origin in $ \mathcal{H}^{\c}$}, {\mk and $h_{\lambda}^{\mathrm{loc}}(0)=0$}.  {\mk Recall that the \eqref{NR}-condition ensures the pullback limit $h^{(1)}_\lambda$ given in \eqref{det_h1} to be well-defined.   According to \cite[Cor.~7.1]{CLW13a}, the manifold function $h^{(1)}_\lambda$ under its form \eqref{h1} provides then} the {\it leading order approximation} of the local invariant manifold function $h_{\lambda}^{\mathrm{loc}}$, {\mk \it{i.e.}}
\be \label{h1-error}
\|h_{\lambda}^{\mathrm{loc}}(\xi) - h^{(1)}_\lambda(\xi)\|_\alpha = o(\|\xi\|_\alpha^k).
\ee

It follows from \eqref{h1-error} that for all {\HL $\epsilon>0$}  sufficiently small, there  exists a neighborhood $\mathfrak{B}_1 \subset \mathfrak{B}$ such that
\be \label{h1-error-2}
\|h_{\lambda}^{\mathrm{loc}}(\xi) - h^{(1)}_\lambda(\xi)\|_\alpha \le  {\HL \epsilon} \|\xi\|_\alpha^{k+1}, \qquad \xi \in \mathfrak{B}_1.
\ee
This together with the non-degeneracy condition on $h^{(1)}_\lambda$ given by \eqref{non-degeneracy} implies that
\be 
\|h_{\lambda}^{\mathrm{loc}}(\xi)\|_\alpha \ge \|h^{(1)}_\lambda(\xi)\|_\alpha - \|h_{\lambda}^{\mathrm{loc}}(\xi) - h^{(1)}_\lambda(\xi)\|_\alpha \ge C \|\xi\|_\alpha^{k} - {\HL \epsilon} \|\xi\|_\alpha^{k+1}.
\ee
By possibly choosing ${\HL \epsilon}$ smaller, and $\mathfrak{B}_1$ to be a smaller neighborhood of the origin, we obtain
\be\label{h1-error-3}
\|h_{\lambda}^{\mathrm{loc}}(\xi)\|_\alpha \ge \frac{C}{2} \|\xi\|_\alpha^{k}, \qquad \xi \in \mathfrak{B}_1.
\ee

We show now that the condition \eqref{PM condition} required in Definition \ref{def:PM} holds for solutions of the uncontrolled PDE \eqref{SEE-null} emanating from sufficiently small initial data on the local invariant manifold $\mathfrak{M}^{\mathrm{loc}}_\lambda$.

For this purpose, we note that for any fixed $t^\ast >0$, by continuous dependence of the solutions to \eqref{SEE-null} on the initial data, given any sufficiently small initial datum on the local invariant manifold $\mathfrak{M}^{\mathrm{loc}}_\lambda$, the solution stays on $\mathfrak{M}^{\mathrm{loc}}_\lambda$ over $[0, t^\ast]$. Let $\mathfrak{B}_2 \subset \mathfrak{B}_1$ be a neighborhood of the origin in $\mathcal{H}^{\c}$ so that each initial datum of the form $y_0 :=\xi + h_{\lambda}^{\mathrm{loc}}(\xi)$, $\xi \in \mathfrak{B}_2$, satisfies the aforementioned property, and the corresponding solution $y(\cdot, y_0; 0)$ satisfies furthermore that
\be \label{yc-range}
y_{\c}(t, y_0; 0):=P_{\c} y(t, y_0; 0) \in \mathfrak{B}_1, \qquad \Forall t \in [0, t^\ast],
\ee
where the latter property can be guaranteed by choosing $\mathfrak{B}_2$ properly thanks again to the continuous dependence of the solution on the initial data.

By the local invariant property of $\mathfrak{M}^{\mathrm{loc}}_\lambda$, we have
\bes
y_{\s}(t, y_0; 0): = P_{\s} y(t, y_0; 0) = h_{\lambda}^{\mathrm{loc}}(y_{\c}(t, y_0; 0)), \qquad \Forall t \in [0, t^\ast].
\ees

Now, for each such chosen initial datum, thanks to \eqref{h1-error-2} and \eqref{yc-range}, we get
\bea
  \int_0^{t^\ast} \bigl \|y_{\s}(t,y_0; 0) - {\mk h^{(1)}_\lambda}(y_{\c}(t, y_0; 0)) \bigr \|_\alpha^2 \, \d t  & = \int_0^{t^\ast} \bigl \|h_{\lambda}^{\mathrm{loc}}(\xi (t)) - {\mk h^{(1)}_\lambda}(y_{\c}(t, y_0; 0)) \bigr \|_\alpha^2 \, \d t \\
& \le \int_0^{t^\ast} {\HL \epsilon} \|y_{\c}(t, y_0; 0)\|_\alpha^{2(k+1)} \, \d t \\
& \le {\HL \epsilon} \max_{t \in [0, t^\ast]} \|y_{\c}(t, y_0; 0)\|_\alpha^{2} \int_0^{t^\ast} \bigl \|y_{\c}(t, y_0; 0)\bigr \|_\alpha^{2k} \, \d t.
\eea
Besides, by \eqref{h1-error-3} we have
\bea
\int_0^{t^\ast} \|y_{\s}(t, y_0; 0)\|_\alpha^2 \, \d t = \int_0^t \|h_{\lambda}^{\mathrm{loc}}(y_{\c}(t, y_0; 0)) \|_\alpha^2 \, \d t \ge \frac{C}{2} \int_0^{t^\ast} \|y_{\c}(t, y_0; 0)\|_\alpha^{2k} \d t.
\eea
We obtain then for all $y_0 = \xi + h_{\lambda}^{\mathrm{loc}}(\xi)$ with $\xi \in \mathfrak{B}_2$ that
\be \label{PM-goal}
\frac{\int_0^{t^\ast} \bigl \|y_{\s}(t,y_0; 0) - h^{\mathrm{(1)}}_{\lambda}(y_{\c}(t, y_0; 0)) \bigr \|_\alpha^2 \, \d t }{\int_0^{t^\ast} \|y_{\s}(t, y_0; 0)\|_\alpha^2 \, \d t} \le \frac{2{\HL \epsilon}}{C} \max_{t \in [0, t^\ast]} \|y_{\c}(t, y_0; 0)\|_\alpha^{2}.
\ee
The RHS can be made less than one by again the continuity argument and by possibly choosing $\mathfrak{B}_2$ to be an even smaller neighborhood.

By appealing to the continuous dependences on initial data $y_0$ and the {\mkr control} $u$ of the solution $y(0,y_0; u)$ to the controlled PDE \eqref{SEE}, there exist an open set $\mathcal{V}$ in $\mathcal{H}_\alpha$ containing the set $\{ y_0 = \xi + h_{\lambda}^{\mathrm{loc}}(\xi) \mid \xi \in \mathfrak{B}_2\}$, and an open set $\mathcal{U}$ of the origin in $L^2(0, t^\ast; \mathcal{H})$, such that the solution $y(0,y_0; u)$ satisfies \eqref{PM-goal} with the RHS of \eqref{PM-goal} staying less than one as $y_0$ various in $\mathcal{V}$ and the control $u$ varies in $\mathcal{U}$. The proof is complete.

\ep

{\mk
We conclude this section by some remarks regarding possible ways of constructing more elaborated finite-horizon PMs as well as PMs relying on decompositions of the phase space $\mathcal{H}$  involving other bases than a standard eigenbasis.

\br \label{Rmk_POD}

\bi

\item[i)] More elaborated backward-forward systems than \eqref{LLL} can be imagined in order to design finite-horizon PMs of smaller parameterization defect than offered by $h^{(1)}_{\lambda}$; see \cite[Section 8.3]{CLW13a}. The idea remains however the same, namely to parameterize the high-modes as pullback limits of some approximation of the time-history of the dynamics of low modes. We refer to Section \ref{Sect_Burgers_h2} for such a parameterization leading in particular to finite-horizon PMs whose $e_n$-components are polynomials of higher order than for those constituting $h^{(1)}_{\lambda}$. As we will see in Section \ref{Sec_numresultsh2}, such higher-order PMs can give rise to a better design of suboptimal solutions to a given optimal control problem (including terminal payoff terms) than those accessible from the leading order finite-horizon PM $h^{(1)}_{\lambda}$; see also Remark \ref{Rem_h2better} below.

\item[ii)] Note also that the usage of the eigenbasis in the decomposition of the phase space $\mathcal{H}$ is not essential for
the definition of the finite-horizon PMs as well as for  the construction of PM candidates based on the backward-forward procedure presented in this section or discussed above. In practice, empirical bases such as the POD basis \cite{Holmes_al12} can be adopted to decompose the phase space into resolved low-mode part and its orthogonal complement (the high-mode part). By doing so, the resulting subspaces $\mathcal{H}^{\c}$ and $\mathcal{H}^{\s}$ are not invariant subspaces of the linear operator $L_\lambda$ anymore, and explicit formulas such as \eqref{h1} should be revised accordingly; this important  point for applications will be addressed elsewhere. 


\ei

\er
}

\section{Finite-Horizon Parameterizing Manifolds for Suboptimal Control of PDEs} \label{Sect_reduced model}
\subsection{Abstract results}
{\mk Given a finite-horizon PM, we present hereafter an abstract formulation of the corresponding reduced equations  from which we will see how suboptimal solutions to the  problem \eqref{P1} can be efficiently synthesized once an analytic formulation of such reduced equations is available; see Sections \ref{Sect_Burgers},  \ref{Sect_Burgers_h2} and \ref{Sect_Burgers_local}.}

The approach consists of reducing the PDE \eqref{SEE} governing the evolution of the state $y(t)$ to an ordinary differential equation (ODE) system which is aimed to model the evolution of the low modes $P_{\c} y(t)$, {\mk by substituting their interactions with the high modes $P_{\s}y(t)$, by means of the parameterizing function $h$ associated with a given PM.}

For simplicity, we assume that the nonlinearity $F$ is bilinear, {\mk denoted by $B$ hereafter so 
that 
\bes
B: \mathcal{H}_\alpha \times \mathcal{H}_\alpha \rightarrow \mathcal{H},
\ees
is thus a continuous bilinear mapping.}

 {\mk For the sake of readability,  the notations introduced in the previous sections are completed by those summarized in Table~\ref{tab:1} below. Note also that throughout this article, $B(v)$ will be sometimes used in place of $B(v,v)$ to simplify the presentation.}


\begin{table}[h]
\caption{Glossary of principal symbols used in Sections~\ref{Sect_reduced model} -- \ref{Sect_Burgers_h2}}
\label{tab:1}       
\centering
\begin{tabular}{ll}
\toprule\noalign{\smallskip}
symbol & terminology \\ 
\noalign{\smallskip}\hline\noalign{\smallskip}
$y_{\c}, y_{\s}$ & the low-mode and high-mode projections of a given PDE solution $y$: $y_{\c}:= P_{\c}y$ and $y_{\s}:= P_{\s}y$  \\ 
$(y^\ast, u^\ast)$ & an optimal pair for the original optimal control problem \eqref{P1} \\ 
$z$ &  state variable of  the {\it PM-based reduced system} \eqref{SEE2a} involved in \eqref{RP1}    \\ 
$(z_{R}^\ast, u_R^\ast)$ & an optimal pair for the reduced problem \eqref{RP1}; $u_R^\ast$ is the {\it PM-based  suboptimal control} for \eqref{P1} \\
$y_{R}^\ast$ &  the {\it suboptimal trajectory} of the underlying PDE driven by $\mathfrak{C} \ur^\ast$ \\
$z^\ast$ &  the trajectory of the PM-based reduced system driven by $P_{\c}  \mathfrak{C}P_{\c} u^\ast$ \\
$l_R$ &  the trajectory $z_{R}^\ast$ ``lifted'' onto the given parameterizing manifold:  $l_R:= z_{R}^\ast + h(z_{R}^\ast)$ \\
$l^\ast$ &  the trajectory $z^\ast$ ``lifted'' onto the given parameterizing manifold:  $l^\ast := z^\ast + h(z^\ast)$ \\
\noalign{\smallskip} \bottomrule 
\end{tabular}
\end{table}

{\mk Recall that the} subspace $\mathcal{H}^{\c}$ is spanned by the first $m$ {\mk dominant} eigenmodes associated with the linear operator $L_\lambda$ for some positive integer $m$. We denote as before its topological complements in $\mathcal{H}$ and $\mathcal{H}_\alpha$ by $\mathcal{H}^{\s}$ and $\mathcal{H}^{\s}_\alpha$, respectively. Let $h:\mathcal{H}^{\c} \rightarrow \mathcal{H}^{\s}_\alpha$  be a finite-horizon PM function associated with \eqref{SEE}; see Definition~\ref{def:PM}. The {\mk corresponding PM-based reduced optimal control problem \eqref{RP1} below, is then built from the following $m$-dimensional PM-based reduced system:
\begin{subequations} \label{SEE2_v1}
\begin{align}
&\frac{\d z}{\d t} = L^{\c}_\lambda z + P_{\c} B(z + h(z)) +  P_{\c} \mathfrak{C} P_{\c} u (t), \quad t \in (0, T],\label{red_sys}\\ 
&\hspace{-4.36cm}\mbox{ supplemented by }\nonumber\\
&z(0)  = P_{\c} y_0 \in \mathcal{H}^{\c};
\end{align}
\end{subequations}
the system \eqref{red_sys} being aimed to model the dynamics of the low modes $P_{\c}y(t)$ by $z(t)$, and the dynamics of the high modes $P_{\s} y(t)$ by $h(z(t))$.  To avoid pathological situations, we will assume throughout this article that $P_{\c} \mathfrak{C} P_{\c}$ is non-zero.

To simplify the presentation, we will assume furthermore  that  the PM function $h$ has been chosen so that 
for any $z(0)$ in $\mathcal{H}^\c$, the problem \eqref{SEE2_v1} admits a well-defined global ($\mathcal{H}^\c$-valued) solution that is continuous in time.} {\hl Such PM functions are identified in the case of a Burgers-type equation in Sections~\ref{Sect_Burgers}--\ref{Sect_Burgers_local}; see also Appendix~\ref{Sect_energy_est} for more details on the corresponding well-posedness problem for the associated reduced systems.}

Note that only the {\mk low-mode projection} of the controller $u$, $P_{\c} u$,  is kept in the above reduced model. In the following we denote {\mk by $\ur := P_{\c} u \in L^2(0,T; \mathcal{H}^{\c})$  this $m$-dimensional controller}. Then, {\mk the problem \eqref{SEE2_v1}} can be rewritten as:
\begin{subequations} \label{SEE2}
\begin{align}
& \frac{\d z}{\d t} = L^{\c}_\lambda z + P_{\c} B(z + h(z)) + P_{\c} \mathfrak{C} \ur (t), \quad t \in (0, T],  \label{SEE2a}\\
& z(0)  = P_{\c} y_0 \in \mathcal{H}^{\c},
\end{align}
\end{subequations}
{\mk and the cost functional \eqref{J_intro} is substituted by
\be  \label{J2}
J_R(z, \ur) := \int_0^T \bigl[ \mathcal{G}\bigl(z(t) + h(z(t))\bigr) + \mathcal{E}( \ur(t)) \bigr] \d t.
\ee
}

The finite-horizon PM-based reduced optimal control problem is then given by:
\be  \label{RP1}  \tag {$\mathcal{P}_{\mathrm{sub}}$}
\begin{aligned}
 \min J_R(z, \ur)  \quad \! \text{s.t.} \quad \!  (z, \ur) \in L^2(0,T; \mathcal{H}^{\c}) \times L^2(0,T; \mathcal{H}^{\c})  \quad \! \text{solves} \quad \! \eqref{SEE2}.
\end{aligned}
\ee

Throughout this section, we assume that the {\mk original problem \eqref{P1} as well as its reduced form \eqref{RP1}} admit {\mk each} an optimal control, denoted {\mk respectively} by $u^\ast$ and $\ur^\ast$. Theorem \ref{THM_controller} below provides then an important  {\it a priori} estimate for the theory. It gives indeed a measure on how far to the optimal control $u^*$ a suboptimal control $u^*_R$ built on a given PM is.  More precisely, under a second-order sufficient optimality condition on the cost functional $J$,  an {\it a priori} estimate of $\|\ur^\ast - u^\ast\|^2_{L^2(0,T; \mathcal{H})}$ is expressed in terms of key quantities associated with a given PM on one hand, and key quantities associated with the optimal control $u^*$, on the other; {\mk see \eqref{thm1:goal} below}. 
These quantities involve the parameterization defects associated with $u^\ast$ and $u_R^{*}$; the energy contained in the high modes of the optimal and suboptimal PDE trajectories associated with  $u^\ast$ and $u_R^{*}$, respectively; and the {\mk high-mode energy remainder $\| P_{\s}u^\ast\|_{L^2(0,T; \mathcal{H})}$ of $u^\ast$}. Our treatment  is here inspired by \cite{IK08} but differs however from the latter  by the use of PMs instead of {\mk AIMs}; the framework of PMs allowing for a natural interpretation of the  error estimate \eqref{thm1:goal} derived hereafter that as we will see in the applications, will help analyze the performances of a PM-based suboptimal controller; see Sections~\ref{Sect_Burgers}--\ref{Sect_Burgers_h2}, and Section~\ref{Sect_Burgers_local}.


\bt \label{THM_controller}
Assume that the optimal control problem \eqref{P1} admits an optimal controller $u^\ast$, where the cost functional $J$ defined in \eqref{J_intro} satisfies the assumptions of Section  \ref{Sect_functfram}.

Assume furthermore there exists $\sigma > 0$ such that the following (local) second order sufficient optimality condition holds:
\be \label{2nd optimality}
J(y(\cdot; v), v) - J(y^\ast, u^\ast)  \ge \sigma \|v - u^\ast\|^2_{L^2(0,T; \mathcal{H})}, 
\ee
where $v \in L^2(0, T; \mathcal{H})$ is chosen from some neighborhood $\mathcal{U}$ of $u^\ast$, and $y(\cdot; v)$ denotes the solution to \eqref{SEE} with $v$ in place of the controller $u$. 



Assume finally that the corresponding 
PM-based reduced optimal control problem \eqref{RP1} admits an optimal controller $\ur^\ast$, which is furthermore contained in $\mathcal{U}$, and that the underlying PM function $h: \mathcal{H}^\c \rightarrow \mathcal{H}_\alpha^\s$ is locally Lipschitz.

Then, the suboptimal controller $\ur^\ast$ satisfies the following error estimate
\bea  \label{thm1:goal}
\|\ur^\ast - u^\ast\|^2_{L^2(0,T; \mathcal{H})} & \le \frac{\mathcal{C}}{\sigma} 
\Bigl( \sqrt{Q(T, y_0; \ur^\ast)} \| y_{R,\s}^\ast\|_{L^2(0,T; \mathcal{H}_\alpha)} \\
& \hspace{1em} +  \sqrt{Q(T,y_0; u^\ast)} \| y_{\s}^\ast\|_{L^2(0,T; \mathcal{H}_\alpha)}  +   {\hh \|\mathfrak{C}\|\| P_{\s}u^\ast\|_{L^2(0,T; \mathcal{H})}} \Bigr),
\eea
where $Q(T, y_0; \ur^\ast)$ and $Q(T, y_0;  u^\ast)$ denote the parameterization defects of the finite-horizon PM function $h$ associated with the controllers in Eq.~\eqref{SEE} taken to be respectively $\ur^\ast$ and $u^\ast$; $ y_{R,\s}^\ast := P_{\s} y_R^\ast$ and $y_{\s}^\ast := P_{\s} y^\ast$ denote the  high-mode projections of the suboptimal trajectory $y_R^\ast$ and the optimal trajectory $y^\ast$ to Eq.~\eqref{SEE} driven respectively  by $\mathfrak{C}\ur^\ast$ and $\mathfrak{C}u^\ast$; and $\mathcal{C}$  denotes a positive constant depending in particular on $T$ and the local Lipschitz constant of $h$; see \eqref{constant_C} below.

\et

{\mk

Besides the suboptimal trajectory $y_R^\ast$, another trajectory of theoretical interest is the ``lifted'' trajectory by the PM function $h$, of the (low-dimensional) optimal trajectory $z_R^\ast:=z(\cdot, P_{\c}y_0; \ur^\ast)$ of the reduced optimal control problem \eqref{RP1}. This lifted trajectory is defined as
\bes
l_R(t):= z_{R}^\ast(t) + h(z_{R}^\ast(t)),
\ees
for which if $z_{R}^\ast$ constitutes a good approximation of the low-mode projection $P_\c y^\ast$ and $h$ has a small parameterization defect\footnote{so that $h(z_{R}^\ast)$ is a good approximation of the high-mode projection $P_\s y^\ast$.}, $l_R$ provides a good approximation of the optimal trajectory $y^\ast$, itself. 

This intuitive idea is made precise in  Corollary \ref{cor_y} below that provides  a general condition under which an error estimate regarding  the distance $\|y^\ast -l_R \|^2_{L^2(0,T; \mathcal{H})}$, between the lifted trajectory $l_R$ and the optimal trajectory $y^\ast$,  can be deduced from the error estimate \eqref{thm1:goal} about the distance between the respective controllers; see \eqref{cor_est2} below. This condition concerns the {\mkr  $L^2$-response} over the interval $[0,T]$ of the PM-based reduced system \eqref{SEE2a} with respect to perturbation of the {\mkr control term} $\mathfrak{C} P_{\c} u^\ast$.

}

\bc  \label{cor_y}

In addition to the assumptions of Theorem \ref{THM_controller}, assume that the  PM-based reduced system \eqref{SEE2a}
satisfies the following sublinear response  property:

\vspace{1ex}
There exist $\kappa> 0$ and a neighborhood $\mathcal{U} \subset L^2(0,T; \mathcal{H}^{\c})$ of $P_{\c} u^\ast$, such that the following inequality holds for all $\ur \in \mathcal{U}$:
\bea \label{sublinear response}
\|z(\cdot, P_{\c}y_0; \ur) - z^\ast(\cdot, P_{\c}y_0; P_{\c} u^\ast)\|_{L^2(0,T; \mathcal{H})}  \le  \kappa \|\ur - P_{\c} u^\ast\|_{L^2(0,T; \mathcal{H})},
\eea
where $z(\cdot, P_{\c}y_0; \ur)$ denotes the solution to \eqref{SEE2} emanating from $P_{\c}y_0$ and driven by $\mathfrak{C} \ur$.
\vspace{1ex}

Then, {\mk the following error estimate} between the optimal trajectory $z_R^\ast:=z(\cdot, P_{\c}y_0; \ur^\ast)$ for the reduced optimal control problem \eqref{RP1} and the low-mode projection  $y_{\c}^\ast := P_{\c}y^\ast$ of the optimal trajectory associated with \eqref{P1}, {\mk holds}:
\bea \label{cor_est1}
& \|y_{\c}^\ast - z_R^\ast\|^2_{L^2(0,T; \mathcal{H})} \\
& \le 2 T \Bigl( \widetilde{\mathcal{C}}_1 Q(T,y_0; u^\ast) \| y^\ast_{\s}\|_{L^2(0,T; \mathcal{H}_\alpha)}^2  + {\hh \widetilde{\mathcal{C}}_2   \|\mathfrak{C}\|^2 \| P_{\s}u^\ast\|^2_{L^2(0,T; \mathcal{H})}}  \Bigr)  \\
& + \frac{2 \kappa^2  \mathcal{C}}{\sigma} \Bigl( \sqrt{Q(T, y_0; \ur^\ast)} \| y_{R,{\s}}^\ast\|_{L^2(0,T; \mathcal{H}_\alpha)}  +  \sqrt{Q(T,y_0; u^\ast)} \| y_{\s}^\ast\|_{L^2(0,T; \mathcal{H}_\alpha)} + {\hh \|\mathfrak{C}\|\| P_{\s}u^\ast\|_{L^2(0,T; \mathcal{H})}} \Bigr),
\eea
where $\mathcal{C}$ is the same positive constant as given by \eqref{thm1:goal} in Theorem~\ref{THM_controller} and 
$\widetilde{\mathcal{C}}_1$, $\widetilde{\mathcal{C}}_2$ are given by \eqref{low_mode_est} in Lemma~\ref{lem:bilinear estimate} below.

Moreover, the following error estimate regarding  the distance $\|y^\ast -l_R \|^2_{L^2(0,T; \mathcal{H})}$, between the lifted trajectory $l_R$ and the optimal trajectory $y^\ast$, holds
\bea  \label{cor_est2}
& \|y^\ast -l_R \|^2_{L^2(0,T; \mathcal{H})} 
   \le 4\bigl[C_\alpha^2 + \widetilde{\mathcal{C}}_1 T \bigl (1+ 2 (C_1C_\alpha \mathrm{Lip}(h)\vert_{V_{\c}})^2\bigr) \bigr] Q(T,y_0; u^\ast) \| y^\ast_{\s}\|_{L^2(0,T; \mathcal{H}_\alpha)}^2 \\
&  \hspace{2em} + \frac{4 \kappa^2  \mathcal{C}}{\sigma} [1 + 2 (C_1C_\alpha \mathrm{Lip}(h)\vert_{V_{\c}})^2] \Bigl( \sqrt{Q(T, y_0; \ur^\ast)} \| y_{R,{\s}}^\ast\|_{L^2(0,T; \mathcal{H}_\alpha)} +   \sqrt{Q(T,y_0; u^\ast)} \| y_{\s}^\ast\|_{L^2(0,T; \mathcal{H}_\alpha)} \Bigr)   \\
&  \hspace{2em}  + {\hh 4 \bigl (1+ 2 (C_1C_\alpha \mathrm{Lip}(h)\vert_{V_{\c}})^2\bigr) \Bigl[ \widetilde{\mathcal{C}}_2 T  \|\mathfrak{C}\|^2 \| P_{\s}u^\ast\|^2_{L^2(0,T; \mathcal{H})} +  \frac{\kappa^2  \mathcal{C}}{\sigma} \|\mathfrak{C}\|\| P_{\s}u^\ast\|_{L^2(0,T; \mathcal{H})} \Bigr]},
\eea
where $C_1$ and $C_\alpha$ are some generic constants given by \eqref{equiv. norms} and \eqref{continuous embedding}, respectively; and  $\mathrm{Lip}(h)\vert_{V_{\c}}$ is the local Lipschitz constant of the PM function $h$ over some bounded set $V_{\c} \subset \mathcal{H}^{\c}$; see \eqref{Vc} and \eqref{Lip h}.

\ec

\medskip

Finally, the last corollary concerns a refinement  of the error estimate \eqref{thm1:goal} which consists of identifying conditions under which the contribution of {\mk the high-mode energy remainder $\| P_{\s}u^\ast\|_{L^2(0,T; \mathcal{H})}$ of the optimal control, can be removed in the upper bound of} $\|\ur^\ast - u^\ast\|^2_{L^2(0,T; \mathcal{H})}$.

\bc  \label{Cor_2}

Assume that the assumptions given in Theorem~\ref{THM_controller} hold. Assume furthermore that  the linear operator $\mathfrak{C}$ leaves stable the subspaces $\mathcal{H}^{\c}$ and $\mathcal{H}^{\s}$, {\it i.e.}
\be  \label{C_condition}
\mathfrak{C} \mathcal{H}^{\c} \subset \mathcal{H}^{\c}  \;  \mbox{ and }\; \; \mathfrak{C} \mathcal{H}^{\s} \subset \mathcal{H}^{\s}.
\ee
Then, {\hl the error estimate \eqref{thm1:goal} reduces to:}
\bea  \label{cor2:goal}
\|\ur^\ast - u^\ast\|^2_{L^2(0,T; \mathcal{H})} & \le \frac{\mathcal{C}}{\sigma} 
\Bigl( \sqrt{Q(T, y_0; \ur^\ast)} \| y_{R,\s}^\ast\|_{L^2(0,T; \mathcal{H}_\alpha)} +  \sqrt{Q(T,y_0; u^\ast)} \| y_{\s}^\ast\|_{L^2(0,T; \mathcal{H}_\alpha)}  \Bigr).
\eea

\ec

Similarly, the corresponding results of Corollary~\ref{cor_y} under the additional condition \eqref{C_condition} amounts to dropping 
the terms involving $P_{\s}u^\ast$ on the RHS of the estimates \eqref{cor_est1} and \eqref{cor_est2}.

\medskip

\subsection{Proofs of Theorem \ref{THM_controller} and Corollaries \ref{cor_y} and \ref{Cor_2}}
For the proofs of the above results, we will make use of the following preparatory lemma.

\bl \label{lem:bilinear estimate}
Given any control $u \in L^2(0, T; \mathcal{H})$, we denote by $y(t)$ the {\mk corresponding} solution to \eqref{SEE}.
{\mk Let $h: \mathcal{H}^\c \rightarrow \mathcal{H}_\alpha^\s$ be a PM function assumed to be locally Lipschitz, and $z(t)$ be the solution to the corresponding PM-based reduced system \eqref{SEE2a} driven by ${\hh P_{\c}} \mathfrak{C}P_{\c}u$  and emanating from $P_\c y(0)$}. 

Then, there exists $\widetilde{\mathcal{C}}_1,  \widetilde{\mathcal{C}}_2 > 0$ such that 
\be  \label{low_mode_est}
\|y_{\c}(t) -  z(t) \|^2 \le  \widetilde{\mathcal{C}}_1 \int_0^t \| y_{\s}(s) - h(y_{\c}(s))\|_\alpha^2 \,\d s  + 
\widetilde{\mathcal{C}}_2 \|\mathfrak{C}\|^2 \int_0^t \| P_{\s}u(s)\|^2 \,\d s, \qquad  t \in [0, T],
\ee
where $y_{\c}:=P_{\c}y$, $y_{\s}:=P_{\s}y$; and $\widetilde{\mathcal{C}}_1$, $\widetilde{\mathcal{C}}_2$  depend in particular on $T$ and the local Lipschitz constant of $h$; see \eqref{y-z est} below.

\el

\bp

Let us {\mk introduce} $w(t):=y_{\c}(t) -  z(t)$. By projecting \eqref{SEE} against the subspace $\mathcal{H}^{\c}$, we obtain
\beas
\frac{\d y_{\c}}{\d t} = L^{\c}_\lambda y_{\c} + P_{\c}B(y_{\c} + y_{\s}) +  P_{\c}\mathfrak{C} u(t),  \quad y_{\c}(0)  = P_{\c} y_0 \in \mathcal{H}^{\c}.
\eeas
This together with \eqref{SEE2_v1} implies that $w$ satisfies the following problem:
\be \label{eq:w}
\frac{\d w}{\d t} = L^{\c}_\lambda w + P_{\c} \bigl( B(y_{\c} + y_{\s}) -  B(z + h(z))\bigr)   +  {\hh P_{\c} \mathfrak{C} P_{\s}u}, \quad w(0)  = 0,
\ee
recalling that $u-P_\c u= P_\s u.$

By taking the $\mathcal{H}$-inner product on both sides of \eqref{eq:w} with $w$, we obtain:
\be  \label{energy est:1}
\frac{1}{2}\frac{\d \|w\|^2}{\d t} = \langle L^{\c}_\lambda w, w \rangle + \langle P_{\c} \bigl( B(y_{\c} + y_{\s}) -  B(z + h(z))\bigr), w \rangle  +  {\hh \langle P_{\c} \mathfrak{C} P_{\s}u, w \rangle }.
\ee

Since $B: \mathcal{H}_\alpha \times \mathcal{H}_\alpha \rightarrow \mathcal{H}$ is a continuous bilinear mapping, there exists $C_B > 0$ such that for any $v_1$ and $v_2$ in $\mathcal{H}_\alpha$, it holds that
\bea \label{bilinear est}
\|B(v_1) - B(v_2)\| & = \|B(v_1, v_1) - B(v_2, v_2)\| \\
& \le \|B(v_1, v_1) - B(v_1, v_2)\| + \|B(v_1, v_2) - B(v_2, v_2)\|  \\
& \le C_B\|v_1\|_\alpha \|v_1 - v_2\|_\alpha + C_B \|v_1 - v_2\|_\alpha \|v_2\|_\alpha \\
& \le C_B ( \|v_1\|_\alpha + \|v_2\|_\alpha)  \|v_1 - v_2\|_\alpha. 
\eea
Thanks to the above bilinear estimate, we get thus 
\be  \label{RHS_control}
 \langle P_{\c} \bigl( B(y_{\c} + y_{\s}) -  B(z + h(z))\bigr), w \rangle 
 \le C_B \bigl( \|y_{\c} + y_{\s}\|_\alpha +\| z + h(z)\|_\alpha \bigr) \, \|y_{\c} + y_{\s} - z - h(z)\|_\alpha \, \|w\|.
\ee

{\mk On the other hand, the assumptions made at the end of Section \ref{Sect_functfram} and in this section regarding the well-posedness problem associated respectively  with Eq.~\eqref{SEE} and the reduced system \eqref{SEE2a}, ensure the existence of a bounded set $V$  in $\mathcal{H}_\alpha$, such that $y(t)$ and $z(t) + h(z(t))$ stay in $V$ for all $t \in [0,T]$. As a consequence,} there exists a constant $C(V) > 0$, such that 
\be \label{sum_est}
C_B \bigl( \|y_{\c}(t) + y_{\s}(t)\|_\alpha + \|z(t) + h(z(t))\|_\alpha \bigr) \le C(V), \qquad t \in [0,T].
\ee
Note also that {\mk by using the local Lipschitz property of $h$, we get}  
\bea  \label{diff_est}
& \hspace{-2em} \|y_{\c}(t) + y_{\s}(t) - z(t) - h(z(t))\|_\alpha \\
& \le \| y_{\c}(t) - z(t)\|_\alpha + \|y_{\s}(t) - h(y_{\c}(t))\|_\alpha + \| h(y_{\c}(t)) - h(z(t))\|_\alpha \\
& \le (1 + \mathrm{Lip}(h)\vert_{V_{\c}}) \|y_{\c}(t) -  z(t) \|_\alpha +  \| y_{\s}(t) - h(y_{\c}(t))\|_\alpha \\
& \le C_1 (1 + \mathrm{Lip}(h)\vert_{V_{\c}}) \|w(t)\| +  \| y_{\s}(t) - h(y_{\c}(t))\|_\alpha,  \qquad t \in [0,T],
\eea
{\mk where  $V_\c=P_\c V$, and $C_1 $ in the  the last inequality denotes the generic positive constant for which 
\be  \label{equiv. norms}  
\|v\|_\alpha \le C_1 \|v\|, \qquad \Forall v \in \mathcal{H}^{\c},
\ee
due to the finite-dimensional nature of $\mathcal{H}^{\c}$. }

By using {\mk now the estimates} \eqref{sum_est} and \eqref{diff_est} in \eqref{RHS_control}, we get
\bea \label{RHS_control:2}
& \langle P_{\c} \bigl( B(y_{\c}(t) + y_{\s}(t)) -  B(z(t) + h(z(t)))\bigr), w(t) \rangle  \\
& \hspace{0.5em} \le  C_1 C(V) (1 + \mathrm{Lip}(h)\vert_{V_{\c}}) \|w(t)\|^2  + C(V)  \| y_{\s}(t) - h(y_{\c}(t))\|_\alpha \|w(t)\| \\
& \hspace{0.5em}  \le C_1 C(V) (1 + \mathrm{Lip}(h)\vert_{V_{\c}}) \|w(t)\|^2 + \frac{[C(V)]^2}{2} \| y_{\s}(t) - h(y_{\c}(t))\|_\alpha^2 + \frac{1}{2}\|w(t)\|^2,
\eea
where we have applied {\mk the standard} Young's inequality $ab < \frac{a^2}{2} +\frac{b^2}{2}$ to derive the last inequality.

Since $L_\lambda$ is assumed to be self-adjoint {\mk with dominant eigenvalue $\beta_1(\lambda)$}, we obtain
\bea \label{RHS_control:3}
\langle L^{\c}_\lambda w(t), w(t) \rangle {\mk = \sum_{i=1}^m  \beta_i(\lambda) |w_i(t)|^2}\le \beta_1(\lambda) \|w(t)\|^2. 
\eea

Note also that
\be  \label{RHS_control:4}
 \langle P_{\c} \mathfrak{C} P_{\s}u, w \rangle \le \|\mathfrak{C}\| \|P_{\s}u(t)\|\|w(t)\| \le \frac{1}{2} \|\mathfrak{C}\|^2 \|P_{\s}u(t)\|^2 + \frac{1}{2} \|w(t)\|^2.
\ee

Using \eqref{RHS_control:2}--\eqref{RHS_control:4} in \eqref{energy est:1}, we obtain
\bea\label{Eq_interm1}
\frac{1}{2}\frac{\d \|w(t)\|^2}{\d t} & \le \Bigl( 1 + \beta_1(\lambda) + C_1 C(V) (1 + \mathrm{Lip}(h)\vert_{V_{\c}}) \Bigr) \|w(t)\|^2 \\
&  \hspace{1.5em} + \frac{[C(V)]^2}{2} \| y_{\s}(t) - h(y_{\c}(t))\|_\alpha^2  + \frac{1}{2} \|\mathfrak{C}\|^2 \|P_{\s}u(t)\|^2.
\eea
Now, by a standard application of the Gronwall's inequality, we obtain for all $t \in [0, T]$, 
\bea \label{y-z est}
\hspace{-1em}\|w(t)\|^2  & = \|y_{\c}(t) - z(t)\|^2 \\
& \le \int_0^t e^{2[1 + \beta_1(\lambda) + C_1 C(V) (1 + \mathrm{Lip}(h)\vert_{V_{\c}})](t-s)}  \biggl( [C(V)]^2 \| y_{\s}(s) - h(y_{\c}(s))\|_\alpha^2 + \|\mathfrak{C}\|^2 \| P_{\s}u(s)\|^2 \biggr) \d s \\
& \le  e^{2[1 + \beta_1(\lambda) + C_1 C(V) (1 + \mathrm{Lip}(h)\vert_{V_{\c}})]T} \biggl( [C(V)]^2\int_0^t \| y_{\s}(s) - h(y_{\c}(s))\|_\alpha^2 \d s +  \|\mathfrak{C}\|^2 \int_0^t \| P_{\s}u(s)\|^2 \,\d s \biggr),
\eea
taking into account that $w(0) = y_{\c}(0) - z(0) = 0$, by assumption. The estimate \eqref{low_mode_est} is thus proved.

\ep



\medskip

We present now the proofs of Theorem~\ref{THM_controller} and Corollaries~\ref{cor_y} and \ref{Cor_2}.
 
\medskip
\noindent {\bf Proof of Theorem~\ref{THM_controller}.}  Let us denote by $y^\ast$ {\mk in $C^1([0,T]; \mathcal{H}) \cap C([0,T]; \mathcal{H}_\alpha)$ the optimal trajectory to the optimal control problem \eqref{P1}, and by $y^\ast_R$ (in the same functional space) the trajectory of Eq.~\eqref{SEE} corresponding to the control $u$ taken to be the optimal (low-dimensional) controller $\ur^\ast$ of the reduced optimal control problem \eqref{RP1}.}

Let us also introduce {\mk the lifted trajectories}
\be
l_R = z^\ast_R + h(z^\ast_R),  \mbox{ and } l^\ast = z^\ast + h(z^\ast), 
\ee
where $z^\ast_R$ and $z^\ast$ are {\mk the} solutions to \eqref{SEE2} driven respectively by $P_\c \mathfrak{C} \ur^\ast(t)$ and {\mk $P_\c \mathfrak{C} P_{\c} u^\ast(t)$}, $t\in[0,T]$.

Thanks to the second order optimality condition \eqref{2nd optimality}, the proof boils down to the derivation of a suitable upper bound for $\Delta:=J(y^\ast_R, \ur^\ast) - J(y^\ast, u^\ast)$, which is organized as follows.

In Step 1, we {\mk  reduce the control of $\Delta$ to the control of $J(y^\ast_R, \ur^\ast) - J(l_R, \ur^\ast)  + J(l^\ast, u^\ast)  - J(y^\ast, u^\ast)$ by using  the optimality property of the pair $(z^\ast_R, \ur^\ast)$ for the reduced problem \eqref{RP1}}. The main interest {\mk in doing so relies on the fact that only $\| y^\ast_R-l_R\|$ and $\| y^\ast- l^\ast\|$ are then determining in the control of $\Delta$; see Step 2.  This leads  in turn to an upper bound of $\Delta$ expressed in terms of key quantities for the design of suboptimal controller in our PM-based theory.}

{\mk  In that respect, the upper bound of $\Delta$ derived in \eqref{J est:2} involves $\| y^\ast_{R,\s} - h(y^\ast_{R,\c})\|_{L^2(0,T; \mathcal{H})}$ and $\| y_{\s}^\ast - h(y_{\c}^\ast)\|_{L^2(0,T; \mathcal{H})}$,  the {\mkr energy} (over the interval $[0,T]$) of the high modes unexplained by  the PM function when applied respectively to $ y^\ast_{R,\c}$ and $y_{\c}^\ast$; and  involves $\|y^\ast_{R,\c} -  z^\ast_R \|_{L^2(0,T; \mathcal{H})}$ and $\|y_{\c}^\ast-  z^\ast \|_{L^2(0,T; \mathcal{H})}$,  the errors associated with the modeling of the $y^\ast_{R,\c}$- and $y_{\c}^\ast$-dynamics by the reduced system \eqref{SEE2a}.

Thanks to Lemma~\ref{lem:bilinear estimate}, we can bound the {\mk two} latter {\mk quantities} by the former {\mk ones} together with a term involving the {\mk energy contained in the high modes of $u^\ast$}.  {\mk This is the purpose of Step 3}.  {\mk The desired result follows then by rewriting the relevant unexplained {\mkr energies} by using  the parameterization defects associated with the PM function $h$ and the controllers $u^\ast$ and $u^\ast_R$.}

\medskip

{\bf Step 1.} Since $(y^\ast, u^\ast)$ is an optimal pair for \eqref{P1}, we get
\bea \label{J est -0}
0 & \le J(y^\ast_R, \ur^\ast) - J(y^\ast, u^\ast) \\
& = J(y^\ast_R, \ur^\ast) - J(l_R, \ur^\ast)  + J(l_R, \ur^\ast)  - J(l^\ast, u^\ast)  + J(l^\ast, u^\ast)  - J(y^\ast, u^\ast).
\eea

Since $(z^\ast_R, \ur^\ast)$ is an optimal pair for the reduced problem \eqref{RP1}, we obtain
\be \label{J-tilde est}
J_R(z^\ast_R, \ur^\ast)  - J_R(z^\ast, P_{\c}u^\ast) \le 0.
\ee
Note also that 
\bes
J(l_R, \ur^\ast) = J_R(z^\ast_R, \ur^\ast), 
\ees
and that according to \eqref{C2}
\bes
J(l^\ast, u^\ast) \ge J_R(z^\ast, P_{\c}u^\ast),
\ees
since $\|P_{\c}u^\ast\|\leq \|u^\ast\|$.

Consequently, 
\be
J(l_R, \ur^\ast)  - J(l^\ast, u^\ast)  \le 0.
\ee 
We obtain then from \eqref{J est -0} that
\bea \label{J est -1}
0 \le J(y^\ast_R, \ur^\ast) - J(y^\ast, u^\ast) \le J(y^\ast_R, \ur^\ast) - J(l_R, \ur^\ast)  + J(l^\ast, u^\ast)  - J(y^\ast, u^\ast).
\eea

\medskip

{\bf Step 2.} Let $V \subset \mathcal{H}_\alpha$ be a bounded set such that 
\be
y^\ast_R(t), \quad l_R(t), \quad  y^\ast(t), \quad l^\ast(t) \in V \qquad \Forall t \in [0, T].
\ee 
Let also 
\be  \label{Vc}
V_{\c} = P_{\c} V.
\ee 
It is clear that $P_{\c}y^\ast_R(t)$, $P_{\c}y^\ast(t)$, $z^\ast_R(t)$ and $z^\ast(t)$ are contained in $V_{\c}$ for all $t \in [0, T]$.

Recalling \eqref{C1}, we denote by $\mathrm{Lip}(\mathcal{G})\vert_{V}$ the Lipschitz constant of $\mathcal{G}:\mathcal{H} \rightarrow \mathbb{R}^+$ restricted to the bounded set $V$. In \eqref{J est -1}, by applying Lipschitz estimates to the $\mathcal{G}$-{\mk part} of the cost functional $J$, we obtain
\bea \label{J est}
0 & \le J(y^\ast_R, \ur^\ast) - J(y^\ast, u^\ast) \\
& \le \mathrm{Lip}(\mathcal{G})\vert_{V} (\|y^\ast_R - l_R\|_{L^1(0,T; \mathcal{H})} +  \|l^\ast - y^\ast \|_{L^1(0,T; \mathcal{H})}) \\
& \le \sqrt{T} \mathrm{Lip}(\mathcal{G})\vert_{V} (\|y^\ast_R - l_R\|_{L^2(0,T; \mathcal{H})} +  \|l^\ast - y^\ast \|_{L^2(0,T; \mathcal{H})}),
\eea
where the last inequality follows from H\"older's inequality.

Recall that $l_R(t) = z^\ast_R(t) + h(z^\ast_R(t))$. Let us also rewrite $y^\ast_R(t)$ as $y^\ast_{R,\c}(t) + y^\ast_{R,\s}(t)$ with $y^\ast_{R,\c}(t)=P_{\c}y^\ast_R(t)$ and $y^\ast_{R,\s}(t)=P_{\s}y^\ast_R(t)$. We obtain then
\bea \label{y-est 1}
 \|y^\ast_R(t) - l_R(t)\|  & \le  \|y^\ast_{R,\c}(t) -  z^\ast_R(t) \| + \| y^\ast_{R,\s}(t) - h(z^\ast_R(t)) \|  \\
&  \le   \|y^\ast_{R,\c}(t) -  z^\ast_R(t) \| + \| y^\ast_{R,\s}(t) - h(y^\ast_{R,\c}(t))\| + \| h(y^\ast_{R,\c}(t))-  h(z^\ast_R(t)) \|.
\eea

Let us denote by $\mathrm{Lip}(h)\vert_{V_\c}$ the Lipschitz constant of $h: \mathcal{H}^{\c} \rightarrow \mathcal{H}^{\s}_\alpha$ restricted to the bounded set $V_{\c}$. We get 
\bea \label{Lip h}
\| h(y^\ast_{R,\c}(t))-  h(z^\ast_R(t)) \|_\alpha & \le \mathrm{Lip}(h)\vert_{V_\c} \|y^\ast_{R,\c}(t) -  z^\ast_R(t) \|_\alpha \\
&  \le C_1 \mathrm{Lip}(h)\vert_{V_\c} \|y^\ast_{R,\c}(t) -  z^\ast_R(t) \|, \quad t \in [0, T],
\eea
where we have used the equivalence between the norms on $\mathcal{H}^{\c}$; see \eqref{equiv. norms}. 

Since $\mathcal{H}_\alpha$ is continuously embedded into $\mathcal{H}$, there exists a generic positive constant $C_\alpha$, such that 
\be \label{continuous embedding}
\|v\| \le C_\alpha \|v\|_\alpha, \qquad \Forall v \in \mathcal{H}_\alpha.
\ee
We obtain then
\be  \label{Lip h-2}
\| h(y^\ast_{R,\c}(t))-  h(z^\ast_R(t)) \|  \le C_1 C_\alpha \mathrm{Lip}(h)\vert_{V_\c} \|y^\ast_{R,\c}(t) -  z^\ast_R(t) \|.
\ee
This together with \eqref{y-est 1} leads to 
\beas
 \|y^\ast_R(t) - l_R(t)\| \le  (1 + C_1 C_\alpha \mathrm{Lip}(h)\vert_{V_\c}) & \|y^\ast_{R,\c}(t) -  z^\ast_R(t) \| \\
& +  \| y^\ast_{R,\s}(t) - h(y^\ast_{R,\c}(t))\|, \quad t \in [0, T].
\eeas

Similarly,
\beas
\| l^\ast(t) - y^\ast(t)\| \le  (1 + C_1 C_\alpha \mathrm{Lip}(h)\vert_{V_{\c}}) & \|y_{\c}^\ast(t) -  z^\ast(t) \| \\
&  +  \| y_{\s}^\ast(t) - h(y_{\c}^\ast(t))\|, \quad t \in [0, T].
\eeas
Reporting the above two estimates into \eqref{J est}, we obtain
\bea  \label{J est:2}
0 & \le J(y^\ast_R, \ur^\ast) - J(y^\ast, u^\ast) \\
& \le 2 \sqrt{T} \mathrm{Lip}(\mathcal{G})\vert_{V}  \Bigl( \| y^\ast_{R,\s} - h(y^\ast_{R,\c})\|_{L^2(0,T; \mathcal{H})} +  \| y_{\s}^\ast - h(y_{\c}^\ast)\|_{L^2(0,T; \mathcal{H})}  \\
& \hspace{1.5em}+ (1 + C_1 C_\alpha \mathrm{Lip}(h)\vert_{V_{\c}})  \bigl ( \|y^\ast_{R,\c} -  z^\ast_R \|_{L^2(0,T; \mathcal{H})} +  \|y_{\c}^\ast-  z^\ast \|_{L^2(0,T; \mathcal{H})} \bigr)  \Bigr).
\eea

\medskip

{\bf Step 3.} By using Lemma~\ref{lem:bilinear estimate} (see \eqref{y-z est} above), we obtain:
\beas
{\hh \|y^\ast_{R,\c} -  z^\ast_R\|_{L^2(0,T; \mathcal{H})} \le  \sqrt{T} C(V) e^{[1 + \beta_1(\lambda) + C_1 C(V) (1 + \mathrm{Lip}(h)\vert_{V_{\c}})]T} \| y^\ast_{R,\s} - h(y^\ast_{R,\c})\|_{L^2(0,T; \mathcal{H}_\alpha)},} 
\eeas
{\hh where we have used $P_{\s}u_R^\ast = 0$ since $u_R^\ast$ lives in $L^2(0,T; \mathcal{H}^{\c})$; and the same lemma leads to}
\beas
\|y_{\c}^\ast & -  z^\ast\|_{L^2(0,T; \mathcal{H})} \\
& \le  \sqrt{T} e^{[1 + \beta_1(\lambda) + C_1 C(V) (1 + \mathrm{Lip}(h)\vert_{V_{\c}})]T}\Bigl(C(V) \| y_{\s}^\ast  - h(y_{\c}^\ast)\|_{L^2(0,T; \mathcal{H}_\alpha)} +  \|\mathfrak{C}\|  \| P_{\s}u^\ast\|_{L^2(0,T; \mathcal{H})} \Bigr).
\eeas

Now, by reporting these estimates in \eqref{J est:2} and using again the property of continuous embedding \eqref{continuous embedding}, we obtain:
\bea  \label{J est:3}
0 & \le J(y^\ast_R, \ur^\ast) - J(y^\ast, u^\ast) \\
&  \le \mathcal{C}(V, \mathrm{Lip}(h)\vert_{V_{\c}}, T) \Bigl( \| y^\ast_{R,\s} - h(y^\ast_{R,\c})\|_{L^2(0,T; \mathcal{H}_\alpha)} +  \| y_{\s}^\ast - h(y_{\c}^\ast)\|_{L^2(0,T; \mathcal{H}_\alpha)} +  {\hh \|\mathfrak{C}\|  \| P_{\s}u^\ast\|_{L^2(0,T; \mathcal{H})} } \Bigr),
\eea
where
\bea  \label{constant_C}
& \mathcal{C}(V, \mathrm{Lip}(h)\vert_{V_{\c}}, T) := 2 C_\alpha  \sqrt{T}\mathrm{Lip}(\mathcal{G})\vert_{V}  \\
& +  2 {\hh \max\{C(V), 1\}} T \mathrm{Lip}(\mathcal{G})\vert_{V} ( 1 + C_1 C_\alpha \mathrm{Lip}(h)\vert_{V_{\c}})  e^{[1 + \beta_1(\lambda) + C_1 C(V) (1 + \mathrm{Lip}(h)\vert_{V_{\c}})]T}.
\eea

In terms of parameterization defects defined in \eqref{Eq_QualPM}, the above estimate \eqref{J est:3} can be rewritten as:
\bea
0 & \le J(y^\ast_R, \ur^\ast) - J(y^\ast, u^\ast)  \\
& \le \mathcal{C}(V, \mathrm{Lip}(h)\vert_{V_{\c}}, T) \Bigl( \sqrt{Q(T, y_0; \ur^\ast)} \| y^\ast_{R,\s}\|_{L^2(0,T; \mathcal{H}_\alpha)}  \\
& \hspace{10em}+  \sqrt{Q(T,y_0; u^\ast)} \| y_{\s}^\ast\|_{L^2(0,T; \mathcal{H}_\alpha)}  +  {\hh \|\mathfrak{C}\|  \| P_{\s}u^\ast\|_{L^2(0,T; \mathcal{H})} } \Bigr),
\eea
where $Q(T, y_0; \ur^\ast)$ and $Q(T, y_0;  u^\ast)$ are the parameterization defects of the finite-horizon PM function $h$ when the control in \eqref{SEE} is taken to be $\ur^\ast$ and $u^\ast$, respectively.  

The proof is complete. 


\medskip
\noindent {\bf Proof of Corollary~\ref{cor_y}.} The estimate given by \eqref{cor_est1} can be derived directly from Theorem~\ref{THM_controller} and Lemma~\ref{lem:bilinear estimate} by noting that
\bes 
\|y_{\c}^\ast - z_R^\ast\|^2_{L^2(0,T; \mathcal{H})} \le 2 \|y_{\c}^\ast - z^\ast\|^2_{L^2(0,T; \mathcal{H})}
+ 2 \|z^\ast - z_R^\ast\|^2_{L^2(0,T; \mathcal{H})}.  
\ees

Indeed, the first term on the RHS above can be controlled as follows by Lemma~\ref{lem:bilinear estimate}:
\beas
\|y_{\c}^\ast - z^\ast\|^2_{L^2(0,T; \mathcal{H})} & \le \int_{0}^T\Bigl(  \widetilde{\mathcal{C}}_1 \int_0^t \| y^\ast_{\s}(s) - h(y^\ast_{\c}(s))\|_\alpha^2 \,\d s   + {\hh \widetilde{\mathcal{C}}_2  \|\mathfrak{C}\|^2 \int_0^t \| P_{\s}u(s)\|^2 \,\d s}   \Bigr) \d t \\
& \le   T \bigl( \widetilde{\mathcal{C}}_1  \| y_{\s}^\ast - h(y_{\c}^\ast)\|^2_{L^2(0,T; \mathcal{H}_\alpha)} + {\hh \widetilde{\mathcal{C}}_2   \|\mathfrak{C}\|^2 \| P_{\s}u\|^2_{L^2(0,T; \mathcal{H})}}  \bigr) \\
& \le T \bigl( \widetilde{\mathcal{C}}_1  Q(T,y_0; u^\ast) \| y^\ast_{\s}\|_{L^2(0,T; \mathcal{H}_\alpha)}^2  + {\hh \widetilde{\mathcal{C}}_2   \|\mathfrak{C}\|^2 \| P_{\s}u\|^2_{L^2(0,T; \mathcal{H})}}  \bigr).
\eeas

For the term $\|z^\ast - z_R^\ast\|^2_{L^2(0,T; \mathcal{H})}$,  according to the condition \eqref{sublinear response} on the sublinear response and Theorem~\ref{THM_controller}, we obtain
\beas
\|z^\ast - z_R^\ast\|^2_{L^2(0,T; \mathcal{H})} & \le  \kappa^2 \|\ur^\ast - P_{\c} u^\ast\|_{L^2(0,T; \mathcal{H})}^2   \le  \kappa^2 \|\ur^\ast - u^\ast\|_{L^2(0,T; \mathcal{H})}^2 \\
& \le  \frac{\mathcal{C}\kappa^2}{\sigma} 
\Bigl( \sqrt{Q(T, y_0; \ur^\ast)} \| y_{R,\s}^\ast\|_{L^2(0,T; \mathcal{H}_\alpha)} \\
& \hspace{10em} +  \sqrt{Q(T,y_0; u^\ast)} \| y_{\s}^\ast\|_{L^2(0,T; \mathcal{H}_\alpha)}  +  {\hh \|\mathfrak{C}\|  \| P_{\s}u^\ast\|_{L^2(0,T; \mathcal{H})} } \Bigr).
\eeas
We obtain then \eqref{cor_est1} by combining the above two estimates.

The estimate \eqref{cor_est2} follows from \eqref{cor_est1} by noting that
\beas
\|y^\ast -( z_R^\ast  & + h(z_R^\ast)) \|^2_{L^2(0,T; \mathcal{H})} \\
& \le 2 \|y_{\c}^\ast - z_R^\ast\|^2_{L^2(0,T; \mathcal{H})}  + 2 \| y_{\s}^\ast - h(z_R^\ast) \|^2_{L^2(0,T; \mathcal{H})} \\
& \le 2 \|y_{\c}^\ast - z_R^\ast\|^2_{L^2(0,T; \mathcal{H})}  + 4  \| y_{\s}^\ast - h(y_{\c}^\ast) \|^2_{L^2(0,T; \mathcal{H})} + 4  \|h(y_{\c}^\ast) - h(z_R^\ast) \|^2_{L^2(0,T; \mathcal{H})} \\
& \le 2 \|y_{\c}^\ast - z_R^\ast\|^2_{L^2(0,T; \mathcal{H})}  + 4 C_\alpha^2 \| y_{\s}^\ast - h(y_{\c}^\ast) \|^2_{L^2(0,T; \mathcal{H}_\alpha)} + 4  \|h(y_{\c}^\ast) - h(z_R^\ast) \|^2_{L^2(0,T; \mathcal{H})};
\eeas 
and that 
\bes
\|h(y_{\c}^\ast) - h(z_R^\ast) \|_{L^2(0,T; \mathcal{H})} \le C_1 C_\alpha \mathrm{Lip}(h)\vert_{V_\c} \|y_{\c}^\ast - z_R^\ast\|_{L^2(0,T; \mathcal{H})};
\ees
see \eqref{Lip h-2} for more details about the derivation of this last inequality {\hh (with $y_{R,\c}^\ast$ therein replaced by $y_{\c}^\ast$ here)}.


\medskip
\noindent {\bf Proof of Corollary~\ref{Cor_2}.} Note that if $\mathfrak{C}$ leaves stable the two subspaces $\mathcal{H}^{\c}$ and $\mathcal{H}^{\s}$, then in Lemma~\ref{lem:bilinear estimate}, the equation \eqref{eq:w} satisfied by the difference $w(t):=y_{\c}(t) -  z(t)$ is simplified into the following:
\bes
\frac{\d w}{\d t} = L^{\c}_\lambda w + P_{\c} \bigl( B(y_{\c} + y_{\s}) -  B(z + h(z))\bigr), \quad w(0)  = 0,
\ees
where the term $P_{\c} \mathfrak{C} P_{\s}u$ vanishes here. Consequently, the terms involving $P_s u$ in the subsequent estimates are dropped out, leading then to the the estimate given in \eqref{cor2:goal}.


\section{2D-Suboptimal Controller Synthesis Based on the Leading-Order Finite-Horizon PM: Application to a Burgers-type Equation} \label{Sect_Burgers}
{\mk We apply in this section and the next, the PM-based reduction approach introduced above for the design of suboptimal solutions to an optimal control  problem of a Burgers-type equation, in the case of globally distributed control {\mkr laws}. The more challenging case of locally distributed control {\mkr laws,} is addressed in Section   \ref{Sect_Burgers_local}.} 
\subsection{{\mk Cost functional of terminal payoff type} for a Burgers-type equation, and existence of optimal solution}\label{Sect_TP_existence}

The model considered here takes the following form, which is posed on the interval $(0, l)$ driven by a {\mk globally} distributed control {\mkr term} {\HL $\mathfrak{C} u(x, t)$}:
\be \label{eq:Burgers}
\frac{\mathrm{d} y}{\d t} = \nu y_{xx}  + \lambda y  - \gamma y  y_x + \mathfrak{C} u(x, t),  \qquad (x,t) \in (0,l) \times (0, T],
\ee
where $\nu, \lambda$ and $\gamma$ are positive parameters, the final time $T > 0$ is fixed, {\mk and conditions on the linear operator $\mathfrak{C}$ are specified in Section~\ref{ss:reduction-h1} below}. 

The equation is supplemented with the Dirichlet boundary condition 
\be\label{Dirc_cond}
y(0,t;u) = y(l,t;u) = 0, \qquad t \in[0, T];
\ee
and appropriate initial condition 
\be  \label{initial:Burgers}
y(x, 0) = y_0(x), \qquad x\in (0,l).
\ee

{\mk The classical Burgers equation (with $\lambda = 0$ in \eqref{eq:Burgers}) has widely served as a theoretical laboratory  to test various methodologies  devoted to the design of optimal/suboptimal controllers {\mkr of} nonlinear distributed-parameter systems; see {\it e.g.} \cite{baker2000,Choi_al93,krstic09,kunisch2004hjb,Volkwein01} and references therein. The inclusion of the term $\lambda y$ here allows for the presence of linearly unstable modes, which lead in turn to the existence of non-trivial (and nonlinearly) stable steady states for the uncontrolled version  of \eqref{eq:Burgers} provided that  $\lambda$ is large enough; see \cite{HW06}. The latter {\mkr property} will be used in the choices of initial data and targets for the associated optimal control problems analyzed hereafter. From a physical perspective, we mention that  \eqref{eq:Burgers} arises in the modeling of flame front propagation \cite{berestycki2001meta}. This model will serve us here to demonstrate the effectiveness of the PM approach introduced above in the design of suboptimal solutions to optimal control problems.} 

{\mk In that respect, we consider the following cost functional  associated with \eqref{eq:Burgers}--\eqref{initial:Burgers}, }
\be  \label{JJ}
J(y, u) = \int_0^T \bigl( \frac{1}{2}\|{\HL y(\cdot, t; y_0, u)}\|^2  + \frac{\mu_1}{2}\|u(\cdot, t)\|^2 \bigr) \d t + \frac{\mu_2}{2} \|y(\cdot, T; y_0, u) - Y\|^2,
\ee
constituted by a {\it running cost} {\hl along the controlled trajectory} and a {\it terminal payoff} term {\hl defining a penalty on the final state}; here $\mu_1$ and $\mu_2$ are some positive constants, $Y \in L^2(0,l)$ is some given target profile, and $\|\cdot\|$ denotes the $L^2(0,l)$-norm.

Compared to the cost functional \eqref{J_intro} associated with the optimal control problem \eqref{P1} given in Section~\ref{Sect_functfram}, we have added here a terminal payoff term $\frac{\mu_2}{2} \|{\HL y(\cdot, T; y_0, u) - Y}\|^2$ to the running cost $\int_0^T \bigl( \frac{1}{2}\|{\HL y(\cdot,t; y_0, u)}\|^2  + \frac{\mu_1}{2}\|{\HL u(\cdot,t)}\|^2 \bigr) \d t$. In Section~\ref{Sect_reduced model}, the optimal control problem \eqref{P1} involving only the latter type of running cost, has served to identify the determining quantities controlling the distance to an optimal control of a suboptimal solution to  \eqref{P1} built from a PM-reduced system; see Theorem \ref{THM_controller} {\HL and Corollary \ref{Cor_2}}.  For a functional cost of type \eqref{JJ}, error estimates similar to \eqref{thm1:goal} {\HL and \eqref{cor2:goal}} can be derived by controlling appropriately  the contribution of the terminal payoff term to  $J(y^\ast_R, \ur^\ast) - J(y^\ast, u^\ast)$ in the estimate \eqref{J est}.  For instance, the error estimate \eqref{cor2:goal} becomes
\bea 
 \|  \ur^\ast  - u^\ast\|^2_{L^2(0,T; \mathcal{H})}  \le \frac{\mathcal{C}}{\sigma}
\Bigl( \sqrt{Q(T, y_0; \ur^\ast)} \| y_{R,\s}^\ast\|_{L^2(0,T; \mathcal{H}_\alpha)} & + \sqrt{Q(T,y_0; u^\ast)} \| y_{\s}^\ast\|_{L^2(0,T; \mathcal{H}_\alpha)} \Bigr) \\
&  + \frac{|C_T(y^*_{R,T}, Y) - C_T(y_T^*, Y)|}{\sigma},
\eea
where $C_T(v, Y) := \frac{\mu_2}{2} \|v - Y\|^2$,  $y^*_{R,T} = y^*_{R}(T)$ and  $y^*_{T} = y^*(T)$. We dealt with the simpler situation of a single running cost type functional in Section~\ref{Sect_reduced model} in order not to overburden the presentation. Furthermore, as we will see in  this section and the forthcoming ones, the error estimates derived in Section \ref{Sect_reduced model}  are sufficient enough to provide useful (and computable) insights to help analyze the performances of a PM-based suboptimal controller.\footnote{Note that in practice, although the second order optimality condition  \eqref{2nd optimality} is difficult to check, the error estimates such as \eqref{cor2:goal} will still demonstrate their relevance for the performance analysis; see Section \ref{Sec_numresults}.}


{\mk The interest of cost functionals such as \eqref{JJ} is that they} arise naturally when the goal is to drive the state $y(\cdot; u)$ of \eqref{eq:Burgers} as close as possible to a target profile $Y$ at the final time $T$, while keeping the cost of the control, expressed by $\frac{\mu_1}{2} \int_0^T \|u(t)\|^2 \d t$, as low as possible. Here, the terminal payoff term gives a measurement of the ``proximity'' to the target $Y$ at the final-time SPDE profile. If one can make $\mu_2 = +\infty$, it means the problem is exactly controllable, if not the system is approximately controllable \cite{Lions88}.

{\mk We turn now to the precise description of the optimal control problem considered in this section and the next.
Adopting the notations of Section \ref{Sect_functfram}, the functional spaces are
\be  \label{spaces}
\mathcal{H}:=L^2(0,l),  \qquad \mathcal{H}_1:=H^2(0,l)\cap H_0^1(0,l),  \qquad \mathcal{H}_{1/2}:= H_0^1(0,l),
\ee 
the linear operator $L_\lambda: \mathcal{H}_1 \rightarrow \mathcal{H}$ is given by
\be  \label{L Burgers}
L_\lambda y := \nu \partial_{xx}^2 y + \lambda y,
\ee 
and the nonlinearity $F$ is expressed by the bilinear term
\bea\label{defB}
B: & \mathcal{H}_{1/2} \times  \mathcal{H}_{1/2}  \rightarrow \mathcal{H}\\
& (y,y) \mapsto B(y,y) := - \gamma y \partial_x y,
\eea
with slight abuse of notations, understanding \eqref{L Burgers} and $ y \partial_x y$ in \eqref{defB} within  the appropriate weak sense.

The optimal control problem {\hh for which we will propose suboptimal solutions takes here the following form}:
\be  \label{BP}  
\begin{aligned}
 & \hspace{5em} \min J(y, u)  \quad \text{with $J$ defined in \eqref{JJ}} \qquad \text{s.t.} \\
&  (y, u) \in L^2(0,T; \mathcal{H}) \times L^2(0,T; \mathcal{H}) \text{ solves the problem \eqref{eq:Burgers}--\eqref{initial:Burgers}}. 
\end{aligned}
\ee
}


It can be checked by standard energy estimates that for any given controller $u \in L^2(0,T; \mathcal{H})$, initial datum $y_0 \in \mathcal{H}$ and any finite $T>0$, there exists a unique weak solution\footnote{in the sense recalled in \eqref{Eq_weaksense} below.}  $y(\cdot; y_0, u)$ for the problem \eqref{eq:Burgers}--\eqref{initial:Burgers} such that $y(\cdot; y_0, u) \in L^2(0,T; \mathcal{H}_{1/2})$ and
$y'(\cdot; y_0, u) \in L^2(0,T; (\mathcal{H}_{1/2})^{-1})$, where $(\mathcal{H}_{1/2})^{-1} = H^{-1}(0,l)$ is the dual of $\mathcal{H}_{1/2} = H_0^1(0,l)$; see {\it e.g.}~\cite{Volkwein01} for the standard Burgers equation subject to affine control.

Note also that $y(\cdot; y_0, u) \in C([0,T]; \mathcal{H}) $ thanks to the continuous embedding
\bes
\mathcal{W}:=\{y\; \vert y \in L^2(0,T; \mathcal{H}_{1/2}) \text{ and } \frac{\d y}{\d t} \in L^2(0,T; (\mathcal{H}_{1/2})^{-1}) \} \subset C([0,T]; \mathcal{H});
\ees
see {\it e.g.} \cite[Sect.~5.9 Thm.~3]{Evans10} for more details. {\mk This last property implies thus that the cost functional $J$ given by \eqref{JJ} is well defined for any pair $(y,u) \in \mathcal{W} \times L^2(0,T; \mathcal{H}) $  that satisfies the problem \eqref{eq:Burgers}--\eqref{initial:Burgers} in the weak sense \eqref{Eq_weaksense}.}

Within this functional setting, the existence of an optimal pair to \eqref{BP}  in $\mathcal{W} \times L^2(0,T; \mathcal{H})$, can be achieved by application of the {\it direct method of calculus of variations} \cite{dacorogna2007direct}.  The closest application of such a method that serves our purpose can be found in the proof of \cite[Prop.~4]{Volkwein01} for the standard Burgers equation where the author considered cost functional of tracking type; the arguments being easily adaptable to cost functional of the form \eqref{JJ}. We provide below a sketch of such arguments.

First note that given a minimizing sequence $\{(y^n, u^n)\} \in (\mathcal{W} \times L^2(0,T; \mathcal{H}))^{\mathbb{N}}$, since the cost functional $J$ defined by \eqref{JJ} is positive (and thus bounded from below) and satisfies 
\bes
J(y,u) \rightarrow \infty \quad \text{ if } \quad \|y\|_{L^2(0,T; \mathcal{H})} \rightarrow \infty \quad \text{ or } \quad \|u\|_{L^2(0,T; \mathcal{H})} \rightarrow \infty,
\ees
the minimizing sequence lives  in a bounded subset of the functional space $\mathcal{W} \times L^2(0,T; \mathcal{H})$. We can then extract a subsequence, say $\{(y^{n_j}, u^{n_j})\}$, which converges weakly to some element $(y^\ast, u^\ast) \in \mathcal{W} \times L^2(0,T; \mathcal{H})$; see {\it e.g.}~\cite[Thm.~3.18]{Brezis2011}. By using the fact that $\mathcal{W}$ is compactly embedded in $L^2(0,T;L^{\infty}(0,l))$ \cite{temam1984navier}, standard energy estimates on the nonlinear term allow to show that actually $(y^\ast, u^\ast)$ satisfies  \eqref{eq:Burgers}--\eqref{initial:Burgers} in the following weak sense, {\it i.e.} for any $\varphi \in L^2(0,T; \mathcal{H}_{1/2})$ and any $T>0$,
\be\label{Eq_weaksense}
\int_0^T \Big(\big\langle \frac{\mathrm{d} y^\ast}{\d t},\varphi \big\rangle_{\mathcal{H}_{1/2}^{-1}; \mathcal{H}_{1/2}}  - \langle B(y^\ast,y^\ast), \varphi \big \rangle_{\mathcal{H}}+ \nu \big\langle y^\ast , \varphi  \big\rangle_{ \mathcal{H}_{1/2}}  - \big\langle \lambda  y^\ast   + \mathfrak{C} u^\ast,\varphi \big\rangle_{\mathcal{H}} \Big)\d t =0,
\ee 
with $y^\ast(0)=y_0.$


Invoking now the lower semi-continuity property of the norm in Banach space (see {\it e.g.} \cite[Prop.~3.5 (iii)]{Brezis2011}) with respect to the convergence in the weak topology, from the functional form of $J$ given in \eqref{JJ} we conclude that $(y^\ast, u^\ast)$ is an optimal pair for the optimal control problem \eqref{BP}. Having ensured the existence of an optimal pair to \eqref{BP}, we turn now to the design of low-dimensional suboptimal  pairs based on the  
(leading-order) parameterizing manifold introduced in Section \ref{ss:h1}.

\subsection{Analytic derivation of the $h^{(1)}_\lambda$-based  2D reduced system for the design of suboptimal controllers}  \label{ss:reduction-h1}

{\mk We present in this section the analytic derivation of the $h^{(1)}_\lambda$-based reduced system on which we will rely to design suboptimal solutions to  problem \eqref{BP}.}
In this respect, we consider {\mk the} particular case where the subspace $\mathcal{H}^{\c}$ of the low-modes is chosen to be the subspace spanned by the first two eigenmodes of the linear operator $L_\lambda$ defined {\mk in} \eqref{L Burgers}.  {\mk Recall} that the eigenvalues of $L_\lambda$ are given by 
\be \label{eq_eigenvalues}
\beta_n(\lambda) := \lambda - \frac{\nu n^2\pi^2}{l^2}, \qquad \qquad n \in \mathbb{N}, 
\ee
and the corresponding eigenvectors are 
\be \label{Eq_eigenmodes}
e_n(x) := \sqrt{\frac{2}{l}}\sin\Bigl(\frac{n\pi x}{l} \Bigr), \qquad \qquad x\in(0,l).
\ee
{\mk Throughout the numerical applications presented hereafter, we will choose $\lambda$ to be bigger than the critical value $\lambda_c:= \frac{\nu \pi^2}{l^2}$ {\HL such that} $L_\lambda$ admits one and only one unstable eigenmode.}  The subspace $\mathcal{H}^{\c}$ {\mk given by}
\be \label{Hc-2modes}
\mathcal{H}^{\c} := \mathrm{span}\{e_1, e_2\},
\ee
{\mk is thus spanned by one unstable and one stable mode.}

 {\mk For the regimes considered hereafter, it can be checked that the \eqref{NR}-condition is satisfied, leading in particular to a well-defined $h^{(1)}_\lambda$}. {\mk We take as a finite-horizon PM candidate, the manifold function $h^{(1)}_\lambda$ provided by the explicit formula \eqref{h1} that we apply to the PDE \eqref{eq:Burgers}. Recall that according to Lemma \ref{Lem:h1_FPM}, the manifold function $h^{(1)}_\lambda$ provides a natural theoretical PM candidate. Numerical results reported in Fig.~\ref{fig:PM_Sect5} will support that this choice is in fact relevant for the regimes analyzed hereafter for the PDE \eqref{eq:Burgers} leading in particular to manifold functions with parameterization defect less than unity as required in Definition \ref{def:PM}.}

{\mk To analyze the performances achieved by the $h^{(1)}_\lambda$-based reduced system in the design of suboptimal solutions to \eqref{BP}, we place ourselves within the conditions of Corollary \ref{Cor_2}. In particular, we assume that the continuous linear operator $\mathfrak{C}: \mathcal{H} \rightarrow \mathcal{H}$ leaves stable  the subspaces $\mathcal{H}^{\c}$ and $\mathcal{H}^{\s}$:
\be
\mathfrak{C} \mathcal{H}^{\c} \subset \mathcal{H}^{\c},   \qquad \mathfrak{C} \mathcal{H}^{\s} \subset \mathcal{H}^{\s}.
\ee

Recall that under such assumptions, the high-mode energy remainder $\| P_{\s}u^\ast\|_{L^2(0,T; \mathcal{H})}$ of the (unknown) optimal controller $u^*$, does not contribute to the estimate of $\|\ur^\ast - u^\ast\|^2_{L^2(0,T; \mathcal{H})}$; leaving the parameterization defect as a key determining parameter in the control of the latter. In particular we will see in Section \ref{Sect_Burgers_h2} that other manifold functions with a smaller parameterization defect than  the one associated with  $h^{(1)}_\lambda$, lead to a design of better suboptimal solutions to  \eqref{BP} than those based on $h^{(1)}_\lambda$.}

{\mk To be more specific,} the operator $\mathfrak{C}$ when restricted to $\mathcal{H}^{\c}$ takes the following form 
\be  \label{C-def}
\mathfrak{C} e_1 = a_{11}e_1 + a_{12} e_2, \qquad \mathfrak{C} e_2 = a_{21}e_1 + a_{22} e_2,
\ee
where the coefficient matrix 
\be \label{M}
M := \begin{pmatrix}
a_{11} & a_{12}\\
a_{21}& a_{22}
\end{pmatrix}
\ee
 is {\mk chosen to be non-trivial to avoid pathological situations}.

Corresponding to the cost functional \eqref{JJ},  the cost associated {\mk with the $h^{(1)}_\lambda$-based reduced system takes the following form}: 
\be  \label{J2_Burgers}
J_R(z, \ur) = \int_0^T \bigl(  \frac{1}{2}\|z(t) + h^{(1)}_\lambda(z(t; P_{\c}y_0, u_R))\|^2  + \frac{\mu_1}{2} \|\ur(t)\|^2 \bigr) \d t + \frac{\mu_2}{2} \|z(T; P_{\c}y_0, u_R) -P_{\c}Y\|^2,
\ee
where $Y \in \mathcal{H}$ is {\mk some} prescribed target.

{\mk Recall that following \eqref{SEE2}, the $h^{(1)}_\lambda$-based reduced system intended to model the dynamics of the low modes $P_{\c}y$, takes the following abstract form:
\bea \label{reduced_Burgers}
& \frac{\d z}{\d t} = L^{\c}_\lambda z + P_{\c} B\Bigl(z + h^{(1)}_\lambda(z), z + h^{(1)}_\lambda(z)\Bigr) + P_{\c} \mathfrak{C} \ur(t), \qquad {\hl t \in (0, T]}, \\
 & z(0)  = P_{\c} y_0 \in \mathcal{H}^{\c},
\eea
where $y_0$ is the initial datum of the original PDE \eqref{eq:Burgers}, and $\ur \in L^2(0,T; \mathcal{H}^{\c})$ is a given control of the reduced system. }

We are thus left with the following reduced optimal control problem {\mk associated with} \eqref{BP}:
\be  \label{RBP} 
\begin{aligned}
 \hspace{-1em} \min J_R(z, \ur)  \quad \! \text{ s.t. } \quad\!  (z, \ur) \in L^2(0,T; \mathcal{H}^{\c}) \times L^2(0,T; \mathcal{H}^{\c})  \quad\! \text{solves} \quad \! \eqref{reduced_Burgers}.
\end{aligned}
\ee

{\mk We turn now to the description of the analytic form of \eqref{RBP}.}

\medskip
{\bf Analytic form of \eqref{RBP}.} {\mk We proceed with the explicit expression of $h^{(1)}_\lambda$ provided by \eqref{h1} that we apply to the Burgers-type equation \eqref{eq:Burgers}. In that respect the 
 nonlinear interactions between the $\mathcal{H}^\c$-modes as projected onto the $\mathcal{H}^\s$-modes given by
 $$B_{i_1i_2}^n:=\langle B(e_{i_1}, e_{i_2}), e_n \rangle,$$
 constitute key quantities to determine. 
In the case of the Burgers-type equation \eqref{eq:Burgers}, they take the following form: 
\bea \label{nonlinear_interaction}
B_{i_1i_2}^n  & = - \gamma \langle e_{i_1} ( e_{i_2})_x, e_n \rangle 
 = \begin{cases}
-  \alpha i_2, & n = i_1 + i_2, \\
- \alpha i_2 \mathrm{sgn}(i_1-i_2), & n = |i_1-i_2|, \\
0, & \text{otherwise},
\end{cases} 
\eea
where 
\be  \label{alpha}
\alpha := \frac{\gamma \pi}{\sqrt{2}l^{3/2}}. 
\ee
In particular, we have  
$$\langle e_{i_1} ( e_{i_2})_x, e_n \rangle = 0,$$ 
for any $n \ge 5$ and $i_1,i_2 \in\{1,2\}$.}

By using the above nonlinear interaction relations in \eqref{h1}, we obtain {\mk thus} the following  expression of $h^{(1)}_\lambda$:
\bea \label{h1_Burgers}
\boxed{h^{(1)}_\lambda(z_1 e_1 + z_2 e_2)  = \alpha_1(\lambda) z_1 z_2 e_3 + \alpha_2(\lambda) (z_2)^2 e_4 , \qquad (z_1, z_2) \in \mathbb{R}^2,}
\eea
where
\bea \label{alpha1-2}
\alpha_1(\lambda) & := - \frac{ 3 \gamma \pi}{\sqrt{2}l^{3/2} (\beta_1(\lambda) + \beta_2(\lambda) - \beta_3(\lambda))}, \\ \alpha_2(\lambda) & := - \frac{\sqrt{2} \gamma \pi}{l^{3/2}(2\beta_2(\lambda) - \beta_4(\lambda))},
\eea
{\mk with the $\beta_i(\lambda)$ given such as given by  \eqref{eq_eigenvalues}.} {\mk Note that this set of eigenvalues obey the \eqref{NR}-condition for any $\lambda$-value of interest here ({\it i.e.}~$\lambda > \lambda_c$). Note also that $\alpha_1(\lambda) < 0$ and $\alpha_2(\lambda) < 0$ for any such $\lambda$.}

Now, {\mk by using \eqref{h1_Burgers}}, we can rewrite \eqref{J2_Burgers} into the following {\mk explicit} form:
\be \label{J2_Burgers-c}
J_R(z, \ur) = \int_0^T [\mathcal{G}(z(t))  + \mathcal{E}(\ur (t))] \d t +  C_T(z(T), P_{\c}Y),
\ee
where
\bea
& \mathcal{G}(z) = \frac{1}{2} \|z + h^{(1)}_\lambda(z)\|^2 = \frac{1}{2} [(z_1)^2 + (z_2)^2 + (\alpha_1(\lambda) z_1 z_2 )^2 + (\alpha_2(\lambda)z_2^2)^2], \\
& \mathcal{E}(\ur)  = \frac{\mu_1}{2} \|\ur\|^2  = \frac{\mu_1}{2}[(\urc{1})^2 +  (\urc{2})^2], 
\eea
and 
\be  \label{C_T_sect5}
 C_T(z(T), P_{\c} Y)  := \frac{\mu_2}{2} \sum_{i=1}^m |z_i(T) - Y_i|^2,
\ee
with $z_i := \langle z, e_i\rangle$, $\urc{i} := \langle \ur, e_i\rangle$, and $Y_i := \langle Y, e_i\rangle$, $i = 1, 2$.

By using furthermore the expression of $h^{(1)}_\lambda$  {\mk given in \eqref{h1_Burgers} into \eqref{reduced_Burgers}, we obtain finally after projection onto $\mathcal{H}^\c$, the following  analytic formulation of  \eqref{reduced_Burgers}}:
\begin{equation}\label{eq:Burgers reduced}
\boxed{
\begin{aligned}
& \frac{\d z_1}{\d t} = \beta_1(\lambda) z_1 + \alpha \Big( z_1z_2 + \alpha_1(\lambda) z_1z_2^2  + \alpha_1(\lambda) \alpha_2(\lambda) z_1 z_2^3 \Big)  + a_{11}\urc{1}(t) + a_{21}\urc{2}(t), \\
& 
\frac{\d z_2}{\d t} = \beta_2(\lambda) z_2 + \alpha\Big( - z_1^2   + 2 \alpha_1(\lambda) z_1^2z_2  + 2 \alpha_2(\lambda)  z_2^3\Big) + a_{12}\urc{1}(t) + a_{22}\urc{2}(t),
\end{aligned}
}
\end{equation}
{\mk where} $\alpha_1(\lambda)$ and $\alpha_2(\lambda)$ are defined in \eqref{alpha1-2}, and $\alpha =   \frac{\gamma \pi}{\sqrt{2}l^{3/2}}.$

Note that for any given initial datum $(z_{1,0}, z_{2,0})$ and any $T>0$, the {\mk $h^{(1)}_\lambda$-based} reduced system \eqref{eq:Burgers reduced} admits a unique solution in $C([0,T]; \mathbb{R}^2)$; this is carried out through some simple but specific energy estimates that are provided in Appendix~\ref{Sect_energy_est} for the sake of clarity.

\subsection{{\mk Synthesis of suboptimal controllers by a Pontryagin-maximum-principle approach}} \label{ss:PMP}

The {\mk analytic form \eqref{eq:Burgers reduced} of the $h^{(1)}_\lambda$-based reduced system \eqref{reduced_Burgers} allows for the use of standard techniques {\mkr from} finite-dimensional optimal control theory to solve the related reduced optimal control problem \eqref{RBP} \cite{bonnard2003singular,Bryson_al75,Kirk12,Kno81,SL12}.}  {\mk We follow below} an indirect {\mk approach relying on} the Pontryagin maximum principle (PMP); see {\it e.g.} \cite{bonnard2003singular,boscain2004optimal,Kirk12,Kno81,PBGM64,SL12}.   {\mk Usually, the use of the Pontryagin maximum principle allows to identify a set of necessary conditions to be satisfied by an optimal solution. However, as we will see, due to the particular form of the cost functionals considered here and the nature of the reduced control system \eqref{eq:Burgers reduced}, these conditions will turn out to be sufficient to ensure the existence of a (unique) optimal control for the reduced problem.} {\mk  Relying on a PMP approach allows also for theoretical insights that can be gained on the  reduced optimal control problem \eqref{RBP} from the (costate-based) explicit formula of the (reduced) optimal controller reachable by such an approach; see \eqref{uc-p_v2} and Lemmas \ref{Lemma3} and \ref{Lem:M} below.}

{\mk In that perspective}, let us denote the {\mk $h^{(1)}_\lambda$-based reduced vector field involved in} \eqref{eq:Burgers reduced}, by 
\bes
f(z, \ur):= (f_1(z,\ur),  f_2(z, \ur))^{\mathrm{tr}}.
\ees
We introduce now the following Hamiltonian associated with the {\mk reduced} optimal control problem \eqref{RBP}:
\be \label{H}
H(z, p, \ur)  := \mathcal{G}(z) + \mathcal{E}(\ur) + p_1 f_1(z, \ur) + p_2 f_2(z, \ur),
\ee
where $p:=(p_1, p_2)^{\mathrm{tr}}$ is the  {\it costate} (or {\it adjoint} state) associated with the state $z =(z_1, z_2)^{\mathrm{tr}}$.

It follows from the Pontryagin maximum principle that for a given pair
\bes
(z_R^\ast, \ur^\ast)  \in L^2(0,T; \mathcal{H}^{\c}) \times L^2(0,T; \mathcal{H}^{\c}) 
\ees
to be optimal for the {\hl reduced} problem \eqref{RBP}, it must satisfy {\mk the following constrained Hamiltonian system:}
\begin{subequations}  \label{Pontryagin relation}
\begin{align}
& \begin{rcases}
\frac{\displaystyle \d z^\ast_{R}}{\displaystyle \d t}  = \nabla_{p}H(z^\ast_{R}, p^\ast_{R}, u^\ast_{R}) = f(z^\ast_{R}, \ur^\ast),\\ 
\frac{\displaystyle \d p^\ast_{R}}{\displaystyle \d t} =  - \nabla_{z}H(z^\ast_{R}, p^\ast_{R}, u^\ast_{R})=  g(z^\ast_{R}, p^\ast_{R}),
\end{rcases} &  (\text{Hamiltonian system for $(z^\ast_R, p^\ast_R)$}) \\
& \nabla_{u_R} H(z^\ast_{R}, p^\ast_{R}, u^\ast_{R})   = 0,  &  (\text{$1^{\mathrm{st}}$-order optimality condition})  \label{Pontryagin-c} \\
 &  p_{R}^\ast(T) =  \nabla_z  C_T(z_R^\ast(T), P_{\c}Y), &  (\text{terminal condition})    \label{Pontryagin-d}
\end{align}
\end{subequations}
where {\hl $\nabla_x$ stands for the gradient operator along the $x$-direction}, $p_R^\ast = p_{R,1}^\ast e_1 + p_{R,2}^\ast e_2$ is the costate associated with $z_R^\ast$, and the vector field $g = (g_1, g_2)^{\mathrm{tr}}$ {\mk has the following expression} 
\bea  \label{eq:g}
g_1(z, p) & := - z_1    - \beta_1(\lambda) p_1  -  \alpha  p_1 z_2 + 2 \alpha p_2  z_1 - \alpha \alpha_1(\lambda) p_1  (z_2 )^2     \\
&  \hspace{2em} - 4 \alpha \alpha_1(\lambda) p_2  z_1 z_2 - (\alpha_1(\lambda))^2 z_1  (z_2 )^2    -  \alpha \alpha_1(\lambda) \alpha_2(\lambda) p_1  (z_2 )^3 ,   \\
g_2(z , p) & := - z_2   - \beta_2(\lambda) p_2   -  \alpha  p_1  z_1   - 2 \alpha \alpha_1(\lambda) p_1  z_1  z_2  - 2 \alpha \alpha_1(\lambda) p_2  (z_1 )^2 \\
& \hspace{2em} + 6 \alpha \alpha_2(\lambda) p_2  (z_2 )^2 - (\alpha_1(\lambda))^2 (z_1 )^2 z_2  \\
& \hspace{2em}  -  3 \alpha \alpha_1(\lambda) \alpha_2(\lambda) p_1  z_1 (z_2 )^2   - 2 (\alpha_2(\lambda))^2  (z_2 )^3.
\eea

Note also that 
\bes
\nabla_{u_R} H(z^\ast_{R}, p^\ast_{R}, u^\ast_{R}) = \Bigl (\mu_1\urc{1}^\ast + a_{11} p_{R,1}^\ast + a_{12}  p_{R,2}^\ast \, , \,  \mu_2\urc{2}^\ast + a_{21} p_{R,1}^\ast +  a_{22} p_{R,2}^\ast \Bigr)^{\mathrm{tr}}.
\ees
{\mk The \text{$1^{\mathrm{st}}$-order optimality condition} \eqref{Pontryagin-c} reduces then to}
\bea  \label{uc-p}
(\urc{1}^\ast, \urc{2}^\ast)  = - \Bigl( \frac{a_{11} p_{R,1}^\ast + a_{12} p_{R,2}^\ast}{\mu_1},  \frac{a_{21} p_{R,1}^\ast +  a_{22} p_{R,2}^\ast}{\mu_1} \Bigr ),
\eea
{\mk which  written into  a compact form, gives}
\be  \label{uc-p_v2}
\boxed{u_R^\ast =  - \frac{1}{\mu_1} M p_{R}^\ast,}
\ee
{\mk where $M$ is the matrix introduced in \eqref{M}.}

Thanks to the relation \eqref{uc-p} between $\ur^\ast$ and the costate $p_{R}^\ast$, we get
\bea \label{f3-f4}
a_{11} \urc{1}^\ast + a_{21} \urc{2}^\ast & = - \frac{1}{\mu_1}\bigl( (a_{11})^2 + (a_{21})^2\bigr) p_{R,1}^\ast - \frac{1}{\mu_1} (a_{11}a_{12} + a_{21}a_{22}) p_{R,2}^\ast \\
& =: f_3(p_{R,1}^\ast,p_{R,2}^\ast), \\
a_{12} \urc{1}^\ast + a_{22} \urc{2}^\ast & = - \frac{1}{\mu_1} (a_{11}a_{12} + a_{21}a_{22}) p_{R,1}^\ast  - \frac{1}{\mu_1}\bigl( (a_{12})^2 + (a_{21})^2\bigr) p_{R,2}^\ast \\
& =: f_4(p_{R,1}^\ast,p_{R,2}^\ast).
\eea

Finally, the {\hl terminal} condition \eqref{Pontryagin-d} leads to 
\be \label{p-condition}
p_{R,i}^\ast(T) = \mu_2 (z_{R,i}^\ast(T) -  Y_i), \qquad i = 1, 2. 
\ee

By using the above relations, we can reformulate the {\hl set of necessary conditions} \eqref{Pontryagin relation} as {\hl the following} boundary-value problem (BVP) {\hl to be satisfied by} $z_{R}^\ast$ and $p_{R}^\ast$:
\bea \label{eq:bvp}
&  \frac{\d z_1}{\d t} = \beta_1(\lambda) z_1 + \alpha z_1 z_2 + \alpha \alpha_1(\lambda) z_1 (z_2)^2  + \alpha \alpha_1(\lambda) \alpha_2(\lambda) z_1 (z_2)^3 + f_3(p_1, p_2), \\
& \frac{\d z_2}{\d t} = \beta_2(\lambda) z_2 - \alpha (z_1)^2   + 2 \alpha \alpha_1(\lambda) (z_1)^2z_2  + 2 \alpha \alpha_2(\lambda)  (z_2)^3  + f_4(p_1, p_2),\\
& \frac{\d p_1}{\d t} =   g_1(z, p),  \\
& \frac{\d p_2}{\d t} =   g_2(z, p),
\eea
subject to the boundary conditions
\be \label{eq:bvp bc}
z_1(0) = \langle y_0, e_1 \rangle, \quad z_2(0) = \langle y_0, e_2 \rangle, \quad  p_1(T) = \mu_2 (z_{1}(T) -  Y_1), \quad p_2(T) = \mu_2 (z_2(T) -  Y_2),
\ee
where $f_3$ and $f_4$ are given by \eqref{f3-f4}, and  $g_1(z,p)$ and $ g_2(z,p)$ are given by \eqref{eq:g}.

Once this BVP is solved,  the corresponding controller $\ur^\ast$ determined by \eqref{uc-p_v2} constitutes then a {\mk natural candidate to solve the $h^{(1)}_\lambda$-based reduced optimal control problem \eqref{RBP}}. For the problem at hand, since the cost functional \eqref{J2_Burgers} is quadratic in $u_R$ and the {\mk dependence on the controller is affine for the system of equations \eqref{eq:Burgers reduced}}, it is known that the controller $\ur^\ast$ so obtained is actually the unique optimal controller of the reduced problem \eqref{RBP}; see {\it e.g.}~\cite[Sect.~5.3]{Kirk12} and \cite{Trelat2012}. This observation also holds for the other reduced optimal control problems derived in later sections.

It is worth mentioning that the solution of the above BVP depends on the coefficient matrix $M$ defined in \eqref{M} associated with the linear operator $\mathfrak{C}$ through {\mk the expressions of} $f_3$ and $f_4$ {\mk given in} \eqref{f3-f4}. {\mk However,} due to the specific form of $f_3$ and $f_4$, different choices of $M$ can lead to the same solution of the BVP.   {\mk More precisely, the solutions of \eqref{eq:bvp}--\eqref{eq:bvp bc} remain unchanged as long as  $M$  stays in the group of $2\times2$ orthogonal matrices. The following lemma summarizes this result}. 

\bl \label{Lemma3}

The solution of \eqref{eq:bvp}--\eqref{eq:bvp bc} is the same for any $M \in O(2)$. 

\el

\bp

The result follows trivially by noting that given any $M \in O(2)$, it holds that $M^{\mathrm{tr}}M = I$. In particular, the following {\mk basic} identities hold:
\bes
(a_{11})^2 + (a_{21})^2 = (a_{12})^2 + (a_{22})^2 = 1, \qquad a_{11}a_{12} + a_{21}a_{22} = 0.
\ees

By using the above {\mk identities} in \eqref{f3-f4}, we obtain for any $M \in O(2)$ that 
\bes
f_3(p_{R,1}^\ast,p_{R,2}^\ast) = - \frac{1}{\mu_1}p_{R,1}^\ast  \qquad f_4(p_{R,1}^\ast,p_{R,2}^\ast) = - \frac{1}{\mu_1}p_{R,2}^\ast,
\ees 
which is independent of $M$. The desired result follows.  
\ep

In connection to the above lemma, let us make {\mk finally} the following basic observation, which {\mk will be of some interest in the} numerical experiments. 

\bl   \label{Lem:M}

For any two bounded linear operators $\mathfrak{C}_i: \mathcal{H} \rightarrow \mathcal{H}$ (i = 1,2), if they leave invariant the subspaces $\mathcal{H}^\c$ and $\mathcal{H}^\s$, and their actions on the low modes differs only by an orthogonal transformation, {\it i.e.,}
\bes
\mathfrak{C}_i \mathcal{H}^\c \subset \mathcal{H}^\c, \qquad  \mathfrak{C}_i \mathcal{H}^\s \subset \mathcal{H}^\s, \qquad P_{\c} \mathfrak{C}_1 =  M P_{\c} \mathfrak{C}_2 \quad \text{with $M \in O(2)$},
\ees
then the optimal pairs $(z_R^{\ast}, u_R^{\ast})$ and $(\overline{z}_R^{\ast}, \overline{u}_R^{\ast})$, corresponding to the reduced optimal control problem \eqref{RBP} with $\mathfrak{C}$ in \eqref{reduced_Burgers} taken to be $\mathfrak{C}_1$ and $\mathfrak{C}_2$ respectively, satisfy the following relation:
\bes
z_R^\ast = \overline{z}_R^\ast, \qquad u_R^\ast = M^{-1} \overline{u}_R^\ast, \qquad J_R(z_R^\ast, z_R^\ast) = J_R(\overline{z}_R^\ast, \overline{u}_R^\ast).
\ees
If we assume furthermore that $P_{\s} \mathfrak{C}_1= P_{\s} \mathfrak{C}_2$, then analogous results hold for the original optimal control problem \eqref{BP}.

\el

\br

{\mk The above result is not limited to the  two-dimensional nature of} $\mathcal{H}^\c$ given by \eqref{Hc-2modes}, and can be generalized to {\mk a higher dimension $m$, as long as} $\mathcal{H}^\c$ is spanned by the first $m$ eigenmodes, {\mk and $M$ lives in $O(m)$.}

\er

\subsection{{\mk Suboptimal pair $(y_R^\ast,u_R^\ast)$ to  \eqref{BP} based on $h^{(1)}_\lambda$: Numerical aspects}}\label{Sec_numaspects}

{\mk  The method used to solve the reduced optimal control problem \eqref{RBP} being clarified in the previous section, we turn now to the practical aspects concerning the synthesis of an $h^{(1)}_\lambda$-based suboptimal pair $(y_R^\ast,u_R^\ast)$ to the optimal control problem \eqref{BP} associated with the Burgers-type equation \eqref{eq:Burgers}. This synthesis is organized in two steps. First, the BVP problem \eqref{eq:bvp}--\eqref{eq:bvp bc} is solved to get the $h^{(1)}_\lambda$-based suboptimal controller $u_R^\ast$ according to the costate-based explicit expression \eqref{uc-p_v2}.  Second, this suboptimal controller is then used in \eqref{eq:Burgers} to get the suboptimal trajectory $y_R^\ast$ driven by $\mathfrak{C} u_R^\ast$. We explain below how these steps are numerically carried out.}

{\mk Recall that the uncontrolled {\mk Burgers-type equation} admits two locally stable steady states $y^{\pm}$ (emerging from a pitch-fork bifurcation) when $\lambda$ is above the critical value $\lambda_c = \frac{\nu \pi^2}{l^2}$ at which the leading eigenmode $e_1$ loses its linear stability \cite{HW06}.  In the experiments below we take $y^+$ as initial data $y_0$, the target $Y$ being specified in Section \ref{Sec_numresults}. }

{\mk Shooting and collocation methods are commonly used  to solve two-point boundary value problems \cite{Asher_al95,Bonnard_al06,Bryson_al75,keller1976,RS72}. A convenient collocation code is the Matlab built-in solver \texttt{bvp4c.m}\footnote{See \cite{KS01} for more details about \texttt{bvp4c}. We also mention that all the numerical experiments performed in this article have been carried out by using the Matlab version \texttt{7.13.0.564} (R2011b).}, which is used to solve the aforementioned BVP \eqref{eq:bvp}--\eqref{eq:bvp bc} as well as other BVPs encountered in later sections.}

{\mk The simulation of the Burgers equation \eqref{eq:Burgers}  as driven by the 2D suboptimal controller $\ur^\ast$ is then performed by means of a semi-implicit Euler scheme where at each time step the nonlinear term $yy_x = (y^2)_x/2$ and the controller $\ur^\ast(x,t)$ are treated explicitly, while the linear term is treated implicitly.} The Laplacian operator is discretized using a standard second-order central difference approximation. The resulting semi-implicit scheme now reads as follows:
\bea \label{Burgers discrete-1}
y_{j}^{n+1} - y_{j}^{n} = 
 \Big( \nu \Delta_d y_{j}^{n+1} + \lambda  y_{j}^{n+1} - \frac{\gamma}{2} \nabla_d\big((y_{j}^{n})^2 \bigr)  + u^{R,n}_j\Big)\delta t,  \quad j \in \{1, \cdots, N_x - 1\},
\eea
where $y_{j}^{n}$ {\mk denotes} the discrete approximation of $y(j\delta x, n\delta t)$; $u^{R,n}_j$, the discrete approximation of $\ur^\ast(j\delta x, n\delta t)$; $\delta x$, the mesh size of the spatial discretization; $\delta t$, the time step;   {\mk  while $\Delta_d$ and $\nabla_d$ denote  the discrete Laplacian and  discrete first-order derivative given respectively by}
\bes
\Delta_d y_{j}^{n}= \frac{y_{j-1}^{n} - 2y_{j}^{n} + y_{j+1}^{n}}{(\delta x)^2}; \qquad \nabla_d \big( (y_j^{n})^2\big) = \frac{ (y_{j+1}^{n})^2 - (y_j^{n})^2 }{\delta x}, \quad j \in \{1, \cdots, N_x - 1\}.
\ees  
The {\mk Dirichlet boundary condition \eqref{Dirc_cond} becomes} 
$$
y_0^n=y_{N_x}^n=0, 
$$
where $N_x +1$ is the  number of grid points used for the discretization of the spatial domain $[0, l]$.

{\mk The time-dependent $(N_x -1)$-dimensional  vector solution to \eqref{Burgers discrete-1}  is denoted by $\mathbf{Y}^n$, and is intended to be an approximation of the suboptimal trajectory $y_R^\ast$ at time $t = n \delta t$.} Let us also denote by $\mathbf{U}^n$ the spatial discretization of $\ur^\ast(x,n \delta t)$ for $x\in [\delta x, l - \delta x]$, given by
\bes
\mathbf{U}^n := \bigl( \ur^\ast(\delta x,n \delta t), \cdots, \ur^\ast((N_x-1)\delta x,n \delta t) \bigr)^{\mathrm{tr}}.
\ees
Then after rearranging the terms, equation \eqref{Burgers discrete-1} can be rewritten into the following algebraic system:
\bea \label{Burgers discrete}
\bigl( (1- \lambda \delta t) \mathbf{I} - \nu  \delta t \mathbf{A}  \bigr)\mathbf{Y}^{n+1} = \mathbf{Y}^{n} -  \frac{\gamma}{2} \delta t \mathbf{B}[\mathbf{S}(\mathbf{Y}^{n})]  + \delta t \mathbf{U}^{n},
\eea
where $\mathbf{I}$ is the $(N_x-1) \times (N_x-1)$ identity matrix, $\mathbf{A}$ is the tridiagonal matrix associated with the discrete Laplacian $\Delta_d$, $\mathbf{B}$ is the matrix associated with the discrete spatial derivative $\nabla_d$, and $\mathbf{S}(\mathbf{Y}^{n})$ denotes the vector whose entries are the square of the corresponding entries of $\mathbf{Y}^{n}$. 

{\mk Since the eigenvalues of $\mathbf{A}$ are given by $\frac{2}{(\delta x)^2} \Big( \cos( \frac{j \pi \delta x}{l}) - 1 \Big)$ ($j = 1, \cdots, N_x - 1$) and the corresponding eigenvectors are the discretized version of the first $N_x-1$ sine modes $e_1, \cdots, e_{N_x-1}$ given in \eqref{Eq_eigenmodes},  the eigenvalues of the matrix $\mathbf{M}:=(1- \lambda \delta t) \mathbf{I} - \nu  \delta t \mathbf{A}$ of the LHS of \eqref{Burgers discrete} can be obtained easily, and the corresponding eigenvectors are still the discretized sine functions}. At each time step, the algebraic system \eqref{Burgers discrete} can thus be solved {\mk efficiently} using the {\it discrete sine transform}. {\mk To do so,} we first compute the discrete sine transform of the RHS {\mk and then divide the elements of the transformed vector by the eigenvalues of $\mathbf{M}$ to which  the inverse discrete sine transform is applied} to find $\mathbf{Y}^{n+1}$; see {\it e.g.}~\cite[{\att Sect.~3.2}]{Eyre98} for more details. {\mk In the numerical results that follow, the} discrete sine transform has been handled by using the {\mk Matlab} built-in function \texttt{dst.m}.

Finally, it is worthwhile mentioning that we have used a uniform time mesh for the integration of the PDE, whereas the $\ur^\ast$ is defined on a non-uniform mesh due to the adaptive mesh feature of the \texttt{bvp4c} solver. This discrepancy is resolved by using linear interpolation to obtain the value of $\ur^\ast$ at the uniform mesh used in the PDE scheme.

{\mk For the sake of comparison, the synthesis of a suboptimal controller based on a two-mode Galerkin approximation has been carried out following the same steps and the same numerical treatment described above. The corresponding suboptimal controller $u_{G}^\ast$ {\mkr associated with} the 2D Galerkin-based reduced optimal problem \eqref{GBP} is also obtained {\it via} a PMP approach which leads to  solving a BVP   described in Appendix \ref{ss:2D-Galerkin}; see  \eqref{BVP-Galerkin}. The same procedure is applied to higher-dimensional Galerkin-based reduced optimal control problems \eqref{GBP'} derived in Appendix \ref{ss:ND-Galerkin}. }

\subsection{2D-suboptimal controller synthesis based on $h^{(1)}_\lambda$, and control performances: Numerical results}\label{Sec_numresults}

We assess in this section the control performances achieved by the $h^{(1)}_\lambda$-based suboptimal pair $(y_R^\ast,u_R^\ast)$ of the optimal control problem \eqref{BP}  such as synthesized according to the procedure described above. These performances are compared with those achieved by a suboptimal solution computed from the 2D Galerkin-based reduced optimal control problem \eqref{GBP}.  In that respect,  the cost \eqref{JJ} evaluated at the suboptimal pair $(y(\cdot; y_0, u_R^\ast), u_R^\ast)$ will be compared with the cost evaluated at the suboptimal pair $(y(\cdot; y_0, u_G^\ast), u_G^\ast)$, where $u_G^\ast$ is the suboptimal controller synthesized from \eqref{GBP}.

We also set the coefficient  $\mu_2$ weighting the terminal payoff part of the cost functional \eqref{JJ} to be sufficiently large so that the  comparison of  the solution profile at the final time $T$ of \eqref{Burgers discrete-1} \textemdash\, driven by the corresponding synthesized controller \textemdash\, with the prescribed target profile $Y$,  provides a  way to visualize the performance of the synthesized suboptimal controller.

The simulations reported below, are performed for $\delta t=0.001$ and $N_x = 251$ with $l = 1.3\pi$ so that $\delta x \approx 0.02$. The system parameters are taken to be $\nu = 1$, $\gamma = 2.5$, and $\lambda = 3 \lambda_c \approx 1.78$. The parameters $\mu_1$ and $\mu_2$ in the cost functional \eqref{JJ} are taken to be $\mu_1 = 1$ and $\mu_2 = 20$. For all the simulations conducted in this article, the relative tolerance  for the \texttt{bvp4c} has been set to $10^{-8}$ and the BVP mesh size parameter has been set to 1.6E$4$. The linear operator $\mathfrak{C}: \mathcal{H} \rightarrow \mathcal{H}$ is taken to be the identity mapping for the sake of simplicity. According to  Lemma~\ref{Lem:M}, any operator $\mathfrak{C}$ such that $P_{\c} \mathfrak{C} \in O(2)$ and $P_{\s}\mathfrak{C} = \mathrm{Id}_{\mathcal{H}^{\s}}$ can be reduced to this case.

 The numerical results at the final time $T=3$ are reported in Fig.~\ref{fig:state_Sect5}. The left panel of this figure presents for this final time, the solution profile  to \eqref{Burgers discrete-1} as  driven by $\ur^\ast$ and $\u_G^\ast$,  respectively.  For these simulations, the target profile has been chosen to be given by 
\be\label{Y_target}
Y = -0.1\langle y^-, e_1\rangle e_1 + 1.6 \langle y^-, e_2 \rangle e_2.
\ee
The right panel of Fig.~\ref{fig:state_Sect5} shows the two components of the synthesized suboptimal controllers $\ur^\ast$ and $\u_G^\ast$.

As can be observed, the {\mk (approximate) PDE final state} $y(T; \ur^\ast)$ associated with the controller $\ur^\ast$ captures the main qualitative feature of the target, while $y(T; u_{G}^\ast)$ associated with the controller $u_{G}^\ast$ {\mk fails in this task}.  {\mk At a more quantitative level,  the relative $L^2$-errors between the respective driven PDE final states  and the target $Y$ are  given by 
\bes
\frac{\|y(T; y_0, u_R^\ast) - Y\|}{\|Y\|} = 22.81\%, \mbox{ and }   \frac{\|y(T; y_0, u_G^\ast) - Y\|}{\|Y\|} = 76.28\%.
\ees
}

 This discrepancy in the control performance {\mk as revealed on the above relative $L^2$-errors,} goes with a noticeable discrepancy between the respective numerical values of the cost, namely 
\bes
J(y(\cdot; y_0, u_R^\ast), u_R^\ast) = 9.75, \mbox{ and }   J(y(\cdot; y_0, u_G^\ast), u_G^\ast) = 30.77.
\ees
These preliminary results clearly indicate that given a decomposition {\mk $\mathcal{H}^\c \oplus \mathcal{H}^\s$ of $\mathcal{H}$},  the slaving relationships between the $\mathcal{H}^\s$-modes and the $\mathcal{H}^\c$-modes such as parameterized by $h^{(1)}_{\lambda}$, participate in improving the control performance of the suboptimal solutions synthesized from a reduced system involving only the (partial) interactions between the  $\mathcal{H}^\c$-modes as modeled by a low-dimensional Galerkin approximation.

\begin{figure}[!hbtp]
   \centering
   \includegraphics[height=0.3\textwidth, width=\textwidth]{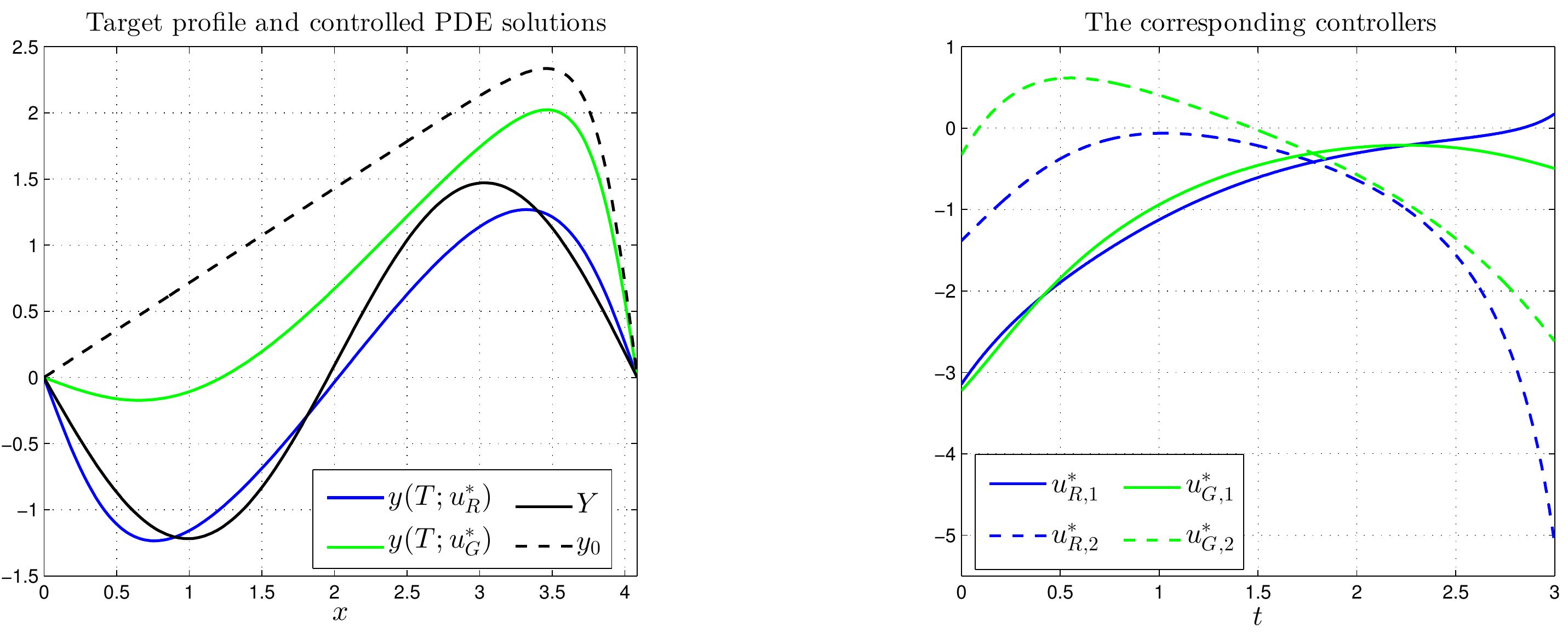}
\vspace{-0.5em}
  \caption{{\footnotesize {\bf Left panel}: PDE solution profiles at the final time $T = 3$ driven respectively by the suboptimal controllers $\ur^\ast$ and $u^\ast_{G}$ with initial profile $y_0$ taken to be $y^+$ (the locally stable positive steady state of the {\mkr uncontrolled} PDE); the target $Y$ (in solid black) is taken to be $-0.1\langle y^-, e_1\rangle e_1 + 1.6 \langle y^-, e_2 \rangle e_2$. {\bf Right panel}: The controller $\ur^\ast= \urc{1}^\ast e_1 + \urc{2}^\ast e_2$ synthesized by the finite-horizon PM-based reduced optimal control problem \eqref{RBP}; and the controller $u_{G}^\ast=u_{G,1}^\ast e_1 + u_{G,2}^\ast e_2$ synthesized by the Galerkin-based reduced optimal control problem \eqref{GBP}. Here, the system parameters are taken to be $l = 1.3\pi$, $\nu = 1$, $\gamma = 2.5$, $\lambda = 3 \lambda_c$. The time step in the PDE solver is $\delta t = 0.001$ and spatial mesh size $\delta x \approx 0.02$. The parameters $\mu_1$ and $\mu_2$ in the cost functional \eqref{JJ} are taken to be $\mu_1 = 1$ and $\mu_2 = 20$.  The corresponding costs are $J(y(\cdot; y_0, u_R^\ast), u_R^\ast) = 9.75$ and  $J(y(\cdot; y_0, u_G^\ast), u_G^\ast) = 30.77$.
}}   \label{fig:state_Sect5}
\end{figure}

To better assess the control performance achieved by the $h^{(1)}_\lambda$-based suboptimal pair $(y_R^\ast,u_R^\ast)$, we compared with the performance achieved by a (suboptimal) solution to  \eqref{BP} based on a high-dimensional Galerkin approximation of \eqref{eq:Burgers}. In that respect, we checked that the cost associated with a suboptimal pair $(y(\cdot; y_0, \widetilde{u}_G^\ast), \widetilde{u}_G^\ast)$, where $\widetilde{u}_G^\ast$ is a 
controller synthesized by solving the BVP \eqref{bvp-Galerkin-m}  associated with an $m$-dimensional Galerkin-based reduced  optimal problem \eqref{GBP'}, can serve as good estimate of the cost associated with the (genuine) optimal solution to the problem \eqref{BP} provided that $m$ is sufficiently large.  We indeed observed that increasing the dimension beyond $m=16$ does not result in significant change of the cost {\mk value} (up to six significant digits) and we thus retained  the results obtained for $m=16$ as reference for providing a good approximation of the optimal solution to \eqref{BP}.
For $m=16$, the corresponding values of the cost \eqref{JJ}, and the relative $L^2$-error for the final time solution profile are given by
\bes
J(y(\cdot; y_0, \widetilde{u}_G^\ast), \widetilde{u}_G^\ast) = 8.41, \mbox{ and }  \frac{\|y(T; y_0, \widetilde{u}_G^\ast) - Y\|}{\|Y\|} = 13.75\%.
\ees
These values when compared with those obtained for the two-dimensional $h^{(1)}_\lambda$-based reduced problem \eqref{RBP} indicates that the two-dimensional controller $u_R^\ast$ already provides a fairly good control performance {\mk but at a much cheaper expense. 
}

On the other hand, the quantitative discrepancy observed {\mk on the cost values and relative $L^2$-errors} between the results based on \eqref{RBP} and those for the original optimal control problem (as indicated by the results based on the high-dimensional Galerkin reduced problem) can be attributed to two main factors according to the {\mk theoretical results of Section~\ref{Sect_reduced model}; see Corollary~\ref{Cor_2} and in particular the {\mk error estimate \eqref{cor2:goal}}. The first factor is related to the  parameterization defect associated with the finite-horizon  PM used here, namely $h^{(1)}_{\lambda}$}; and the second concerns the energy kept in the high modes of the solution {\mk either driven by the suboptimal controller $u_R^\ast$ or the optimal controller $u^*$ itself.}

For the remaining part of this section, we {\mk report on detailed numerical results which further emphasize the practical relevance of the aforementioned theoretic results provided by Corollary~\ref{Cor_2}}. These numerical results {\mk shown} in Figs.~\ref{fig:PM_Sect5} and \ref{fig:Energy_Sect5} are carried out by varying the final time $T$ in the range $[0.1, 5]$ while keeping other parameters the same as used in Fig.~\ref{fig:state_Sect5}.

Panel (a) of Fig.~\ref{fig:PM_Sect5} shows the cost {\mkr values}, when $T$ is varied, associated with the suboptimal pairs $(y_R^\ast,u_R^\ast)$ on one hand (blue curve), and associated with the suboptimal pairs $(\widetilde{y}_G^\ast, \widetilde{u}_G^\ast)$, on the other hand (black curve).  As one can observe up to  $T=3$, the suboptimal controllers $u_R^\ast$ synthesized from the $h^{(1)}_\lambda$-based reduced problem \eqref{RBP} gives access to suboptimal solutions whose cost values are close to those achieved by the optimal ones\footnote{As approximated from the 16-dimensional Galerkin-based reduced  optimal problem \eqref{GBP'}.}. Such good performances starts however to noticeably deteriorate as $T$ increases from $T=3$. 

The reasons of this deterioration are actually rich of teaching, as we explain now. If the error estimate \eqref{cor2:goal} is meaningful, analyzing its main constitutive elements should help understand what causes this deterioration. In that respect, we computed (i) the  corresponding parameterization defects\footnote{{\mk Note that, given a suboptimal controller,  the computation of the parameterization defects here and in latter sections, has been performed by integrating   the  discrete form \eqref{Burgers discrete-1} of  \eqref{eq:Burgers}, and by using the formula \eqref{Eq_QualPM}, where the $H^1$-norm has been used in place of the  $\|\cdot\|_{\alpha}$-norm; see Definition \ref{def:PM} and Section \ref{Sect_TP_existence} for the functional spaces defined in \eqref{spaces}.}} associated with $h^{(1)}_\lambda$ and a given suboptimal controller $u_R^\ast$, and (ii) the energy contained in the high modes of the PDE solution either driven by the suboptimal controller $u_R^\ast$ (leading to the suboptimal trajectory $y_R^\ast$) or the (sub)optimal controller $\widetilde{u}_G^*$ (leading to the (sub)optimal trajectory $\widetilde{y}_G^\ast$).

As a first result, the panels  (b)--(f) of Fig.~\ref{fig:PM_Sect5} show that $h^{(1)}_\lambda$ provides a finite-horizon PM for the whole range of $T$ analyzed here. The parameterization defects of $h^{(1)}_\lambda$ is furthermore robust with respect to variations of $T$, reaching a (nearly) constant value of about 0.57 for $T \ge 1$. At the same time, a substantial growth of the energy contained in the high modes of the suboptimal trajectories $y_R^\ast$ ({\it i.e.} $\|P_{\s} y_R^\ast (t)\|_{H^1(0,l)}$), is observed from $T=3$ to $T=5$ while $\|P_{\s} \widetilde{y}_G^\ast(t)\|_{H^1(0,l)}$ does not change significantly; see Fig.~\ref{fig:Energy_Sect5}.  A closer look at the numbers reveals that 
\begin{align*}
\begin{rcases}
& Q(T, y_0; u_R^\ast) = 0.57, \qquad \|P_{\s} y_R^\ast\|_{L^2(0,T; H^1(0,l))} = 2.26, \\
& Q(T, y_0; \widetilde{u}_G^\ast) = 0.63, \qquad \|P_{\s} \widetilde{y}_G^\ast\|_{L^2(0,T; H^1(0,l))} = 2.15, \\
\end{rcases}
\qquad \text{ for } \quad T = 3, \\
\begin{rcases}
& Q(T, y_0; u_R^\ast) = 0.59, \qquad \|P_{\s} y_R^\ast\|_{L^2(0,T; H^1(0,l))} = 3.0, \\
& Q(T, y_0; \widetilde{u}_G^\ast) = 0.57, \qquad \|P_{\s} \widetilde{y}_G^\ast\|_{L^2(0,T; H^1(0,l))} = 2.13,
\end{rcases}
\qquad \text{ for } \quad T = 5,
\end{align*}
which clearly shows that the RHS of the error estimate \eqref{cor2:goal} experiences a growth of about $15 \%$ when $T$ increases from $T=3$ to $T=5$.
This growth of the RHS of \eqref{cor2:goal}  comes with a growth related to the low-mode part of the LHS  of \eqref{cor2:goal}, {\it i.e.} $\|P_\c (u_R^\ast-\widetilde{u}_G^\ast)\|_{L^2(0,T;L^2(0,l))}^2$, of  about $10 \%$.   This deviation from $\widetilde{u}_G^\ast$, observed on its low-mode part, is consistent with the substantial growth observed on the cost value $J(y_R^\ast, u_R^\ast)$ as shown in Fig.~\ref{fig:PM_Sect5} (a).

To summarize, the error estimate \eqref{cor2:goal} given in Corollary~\ref{Cor_2} provides useful (and computable) insights that can be used to guide the design of PM-based suboptimal controllers with good control performance.    In particular, it addresses the importance of constructing PMs with small parameterization defects on one hand, while keeping small the energy contained in the high-modes, on the other. While the latter factor can be conceivably alleviated by increasing the dimension of the reduced phase space $\mathcal{H}^{\c}$, finite-horizon PMs with smaller parameterization defects than proposed by $h^{(1)}_{\lambda}$ can be thus expected to be even more useful for the design of low-dimensional suboptimal controllers with good performances. The next section addresses the construction of such finite-horizon PMs.

\begin{figure}[hbtp]
\centering
\includegraphics[height=0.55\textwidth, width=0.8\textwidth]{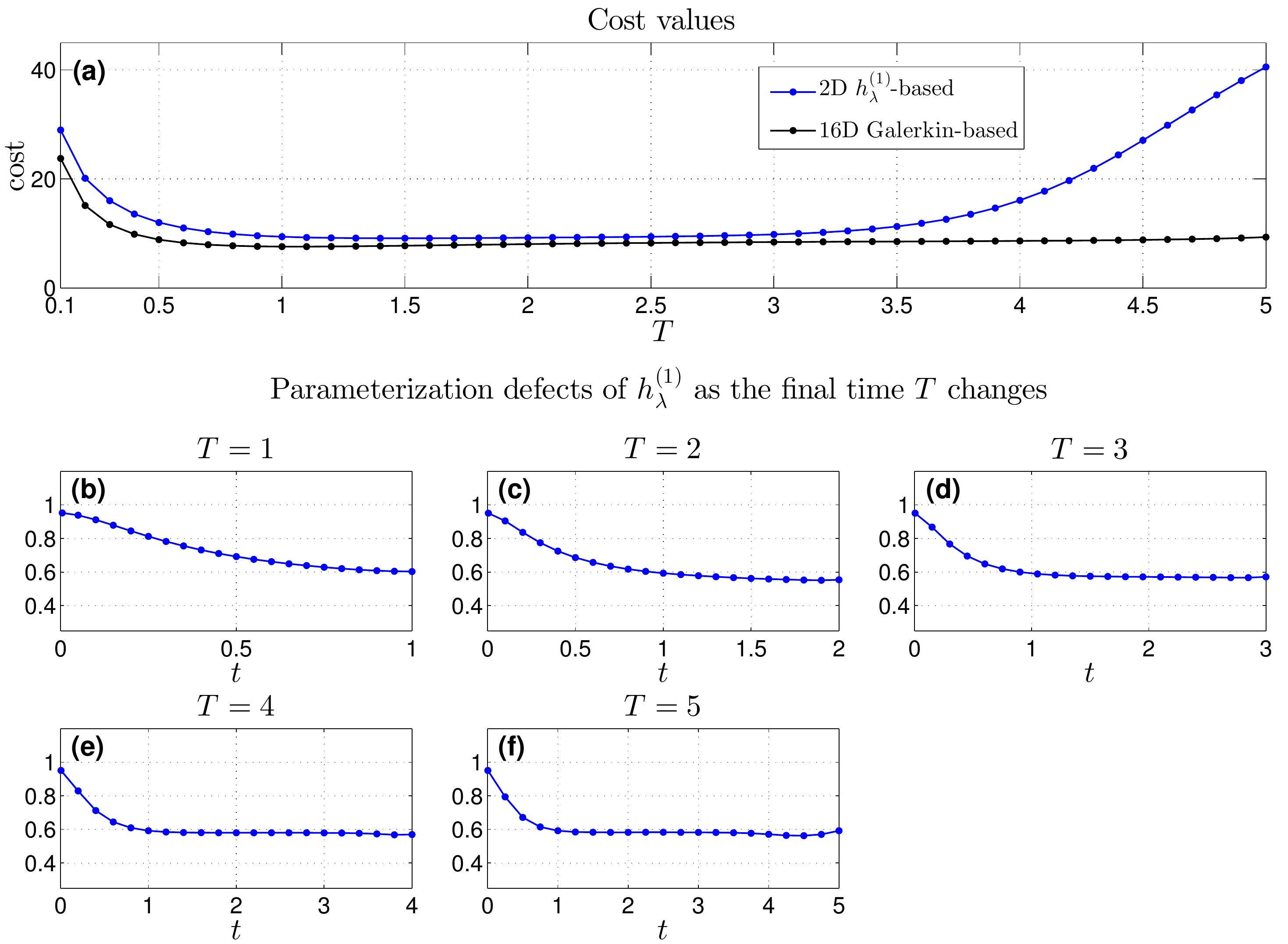}
\vspace{-0.5em}
\caption{{\footnotesize {\bf(a)}: The values of the corresponding cost functional $J$ defined by \eqref{JJ} associated with the suboptimal pair $(y_R^\ast,u_R^\ast)$ as well as the suboptimal pair $(\widetilde{y}_G^\ast, \widetilde{u}_G^\ast)$ as the final time $T$ various in $[0.1, 6]$, where $u_R^\ast$ denotes the suboptimal controller synthesized by the $h^{(1)}_\lambda$-based reduced problem and $\widetilde{u}_G^\ast$ the $16$-dimensional Galerkin based one; {\bf (b)-(f)}: The parameterization defect associated with the finite-horizon PM $h^{(1)}_\lambda$ over the time interval $[0, T]$ for various values of $T$. The parameters are the same as given in Fig.~\ref{fig:state_Sect5}.}} \label{fig:PM_Sect5}
\end{figure}

\begin{figure}[hbtp]
\centering
\includegraphics[height=0.3\textwidth, width=0.8\textwidth]{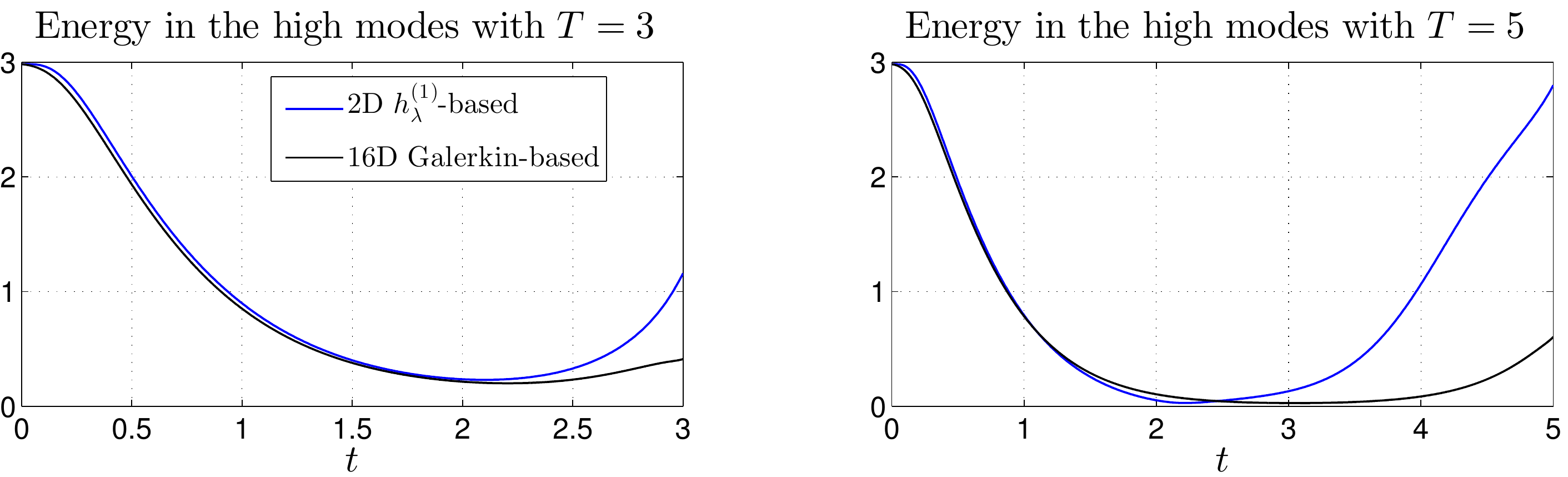}
\vspace{-0.5em}
\caption{{\footnotesize Energy contained in the high modes of the suboptimal trajectories $y_R^\ast$ and $\widetilde{y}_G^\ast$ for $T =3$ (left panel) and $T = 5$ (right panel). The plotted curves are $\|P_{\s} y_R^\ast(t)\|_{H^1(0,l)}$ (blue) and $\|P_{\s} \widetilde{y}_G^\ast(t)\|_{H^1(0,l)}$ (black). The parameters are the same as given in Fig.~\ref{fig:state_Sect5}.}} \label{fig:Energy_Sect5}
\end{figure}
\br \label{Rmk_bocop}
We mention that the numerical results reported in Fig.~\ref{fig:state_Sect5}  have been compared  with those obtained by solving the reduced optimal control problem \eqref{RBP} with the \texttt{BOCOP} toolbox \cite{Bocop}\footnote{In contrast to the indirect method adopted above, \texttt{BOCOP}  uses a direct method combining discretization and interior-point methods to solve the reduced optimal control problem  \eqref{RBP}  as implemented in the solver \texttt{IPOPT} \cite{wachter2006}; see the webpage \url{http://bocop.org} for more information.}. For the parameters used, the relative error under the $L^2$-norm between the controllers  numerically obtained by this toolbox and by our calculations has been observed to be within a margin of $0.1\%$. {\mk For the sake of reproducibility of the results for \eqref{RBP}, we provide the following numerical values of the components of $Y$ used in \eqref{Y_target}:  $\langle Y, e_1\rangle = 0.2561$  and $\langle Y,e_2\rangle= -1.9193$.}
\er

\section{2D-Suboptimal Controller Synthesis Based on Higher-Order Finite-Horizon PMs}  \label{Sect_Burgers_h2}

As illustrated in the previous section {\hh in the context of  a Burgers-type equation},  the finite-horizon PM $h^{(1)}_\lambda$ based on the simple one-layer backward forward system \eqref{LLL}, {\mk can be used efficiently to} obtain low-dimensional {\mk suboptimal controllers with relatively good performances} for certain cases.  Figures \ref{fig:PM_Sect5} and \ref{fig:Energy_Sect5} indicate {\mk that these performances can be altered when the parameterization defects associated with $h^{(1)}_\lambda$ is not specially small, while the energy contained in the high modes of the solution \textemdash\, either driven by the suboptimal controller $u_R^\ast$ or the optimal controller $u^*$ itself \textemdash\, get large, in agreement with the theoretical predictions of Corollary~\ref{Cor_2}.  The error estimate \eqref{cor2:goal} suggests that other finite-horizon PMs with smaller parameterization defects than $h^{(1)}_\lambda$  should help in the synthesis of suboptimal controllers with better performances. The main purpose of this section is to build effectively such PMs that {\mkr in particular} add higher-order terms to $h^{(1)}_{\lambda}$ (Theorem \ref{THM_h2} below) {\mkr which} will turn out to play a crucial role {\mkr to improve the performances} of the $h^{(1)}_\lambda$-based suboptimal controllers encountered so far; see Remark \ref{Rem_h2better} below.}

\subsection{Higher-order finite-horizon PMs based on two-layer backward-forward system: Analytic derivation}  \label{ss:h2}

{\mk We} follow \cite[Sect.~11]{CLW13a} and consider the following {\it two-layer backward-forward system} associated with the {\mk uncontrolled version of \eqref{eq:Burgers}}:
\begin{subequations}\label{L2}
\begin{align}
& \frac{\mathrm{d} y^{(1)}_{\c}}{\mathrm{d} s} =  L_\lambda^{\c} y^{(1)}_{\c}, 
 && s \in [ -\tau, 0],   \quad \; y^{(1)}_{\c}(s)\vert_{s=0} = \xii,  \label{L2-uc1} \\
& \frac{\mathrm{d} y^{(2)}_{\c}}{\mathrm{d} s} =  L_\lambda^{\c} y^{(2)}_{\c} + P_{\c} B( y^{(1)}_{\c},  y^{(1)}_{\c}), 
 && s \in [ -\tau, 0],   \quad \; y^{(2)}_{\c}(s)\vert_{s=0} = \xii,  \label{L2-uc2} \\
& \frac{\mathrm{d} y^{(2)}_{\s}}{\mathrm{d} s} =  L_\lambda^{\s} y^{(2)}_{\s}  +  P_{\s} B(y^{(2)}_{\c}, y^{(2)}_{\c}), 
&&  s \in [-\tau, 0],  \quad y^{(2)}_{\s}(s) \vert_{s=-\tau}= 0, \label{L2-us}
\end{align}
\end{subequations}
where  $L_\lambda^{\c} := P_{\c} L_\lambda$, $L_\lambda^{\s} := P_{\s} L_\lambda$, and $\xi \in \mathcal{H}^{
\c}$. 

{\mk Similar to the one-layer backward-forward system \eqref{LLL}}, the above system is integrated using {\mk a two-step {\it backward-forward integration procedure} where Eqns.~\eqref{L2-uc1}-\eqref{L2-uc2}  are integrated first backward, and Eq.~\eqref{L2-us} is then integrated forward.  }
We will emphasize the dependence {\mk on $\xi$  of the high-mode component $y_{\s}^{(2)}$ of this system as $y_{\s}^{(2)}[\xi]$.}

Theorem \ref{THM_h2} below identifies {\it non-resonance  conditions}  \eqref{NR2} under which the {\it pullback limit} of  $y_{\s}^{(2)}[\xi]$ exists {\mk as $\tau\rightarrow \infty$}. In particular, it provides an analytical {\mk expression} of this pullback limit. {\mk As it will be supported by the numerical results of Section \ref{Sec_numresultsh2}, this pullback limit will turn out to give access to finite-horizon PMs for a broad class of targets.}

\bt \label{THM_h2}

Consider the two-layer backward-forward system \eqref{L2} associated with the {\hl uncontrolled} Burgers-type equation \eqref{eq:Burgers}, {\hl {\it i.e.}~with $\mathfrak{C} = 0$}. Let {\hl $\mathcal{H}^{\c}$ be the subspace} spanned by the first two eigenmodes $e_1$ and $e_2$ of the corresponding linear operator $L_\lambda$ defined in \eqref{L Burgers}. Assume that the eigenvalues of $L_\lambda$ satisfy the following non-resonance conditions:
\begin{equation}  \label{NR2} \tag{NR2}
\begin{aligned} 
& \beta_1(\lambda) + \beta_2(\lambda) - \beta_3(\lambda) > 0, \qquad && \beta_1(\lambda) + 2 \beta_2(\lambda) - \beta_3(\lambda)  > 0, \\
&3 \beta_1(\lambda)  - \beta_3(\lambda)  > 0 , \qquad && 3 \beta_1(\lambda) + \beta_2(\lambda) - \beta_3(\lambda) > 0, \\
& 2 \beta_1(\lambda)  + \beta_2(\lambda) - \beta_4(\lambda) > 0,  && 4 \beta_1(\lambda) - \beta_4(\lambda) > 0,\\
&  2 \beta_2(\lambda) - \beta_4(\lambda) > 0.
\end{aligned}
\end{equation}

Then the pullback limit of the solution $y_{\s}^{(2)}[\xii]$ to \eqref{L2} exists and is given by:
\be \label{eq:h2}
\boxed{h^{(2)}_\lambda(\xi)  :=\lim_{\tau \rightarrow +\infty} y^{(2)}_{\s}[\xi]{\mk(-\tau, 0)}= 
 \int_{-\infty}^0 e^{-\tau' L^{\s}_\lambda} P_{\s} B\bigl(y^{(2)}_{\c}(\tau'), y^{(2)}_{\c}(\tau')\bigr) \d \tau', \quad \Forall \xii \in  \mathcal{H}^{\c}.}
\ee

Under the above conditions,  $h^{(2)}_\lambda$ has {\hl furthermore} the following analytic expression:
\bea \label{h2 expansion}
 h^{(2)}_\lambda(\xi_1e_1+\xi_2e_2) & = h^{(2),3}_\lambda(\xi_1,\xi_2) e_3 + h^{(2),4}_\lambda(\xi_1,\xi_2) e_4,  \qquad (\xi_1, \xi_2) \in \mathbb{R}^{2},
\eea
where
\begin{subequations}\label{coeffs_h2_a}
\begin{align}
 h^{(2),3}_\lambda(\xi_1,\xi_2) & := \langle  h^{(2)}_\lambda(\xi_1e_1+\xi_2e_2), e_3 \rangle \nonumber \\
& =  \mathbi{A} \xi_1 \xi_2 + \mathbi{B} (\xi_1)^3 + \mathbi{C} \xi_1(\xi_2)^2 + \mathbi{D} (\xi_1)^3 \xi_2,  \label{h23} \\
 h^{(2),4}_\lambda(\xi_1,\xi_2) & := \langle  h^{(2)}_\lambda(\xi_1e_1+\xi_2e_2), e_4 \rangle =  \mathbi{E} (\xi_2)^2 + \mathbi{F} (\xi_1)^2 \xi_2 + \mathbi{G} (\xi_1)^4,   \label{h24}
\end{align}
\end{subequations}
with
\bea\label{coeffs_h2}
\mathbi{A} & =  - \frac{3\alpha}{\beta_{1}(\lambda) + \beta_{2}(\lambda) -\beta_{3}(\lambda)}, \\
 \mathbi{B} & = - \frac{3\alpha^2}{(3 \beta_{1}(\lambda)  -\beta_{3}(\lambda))  (\beta_{1}(\lambda) + \beta_{2}(\lambda) -\beta_{3}(\lambda))}, \\
\mathbi{C} & = \frac{3\alpha}{  (\beta_{1}(\lambda) + 2 \beta_{2}(\lambda) -\beta_{3}(\lambda) ) (\beta_{1}(\lambda) + \beta_{2}(\lambda) -\beta_{3}(\lambda) )}, \\
\mathbi{D} & = \frac{3\alpha^3}{ (3 \beta_{1}(\lambda) - \beta_{3}(\lambda)) (\beta_{1}(\lambda) + \beta_{2}(\lambda) -\beta_{3}(\lambda) ) (\beta_{1}(\lambda) + 2 \beta_{2}(\lambda) -\beta_{3}(\lambda))} \\
& + \frac{3\alpha^3 }{ (3 \beta_{1}(\lambda) - \beta_{3}(\lambda)) ( 3\beta_{1}(\lambda) + \beta_{2}(\lambda) -\beta_{3}(\lambda) )(\beta_{1}(\lambda) + 2 \beta_{2}(\lambda) -\beta_{3}(\lambda))}, \\
\mathbi{E} & = - \frac{2 \alpha}{\beta_{2}(\lambda)  -\beta_{4}(\lambda)}, \hspace{0.5em} \mathbi{F} = - \frac{4 \alpha^2}{ (2 \beta_{1}(\lambda) + \beta_{2}(\lambda) -\beta_{4}(\lambda) ) (2 \beta_{2}(\lambda) -\beta_{4}(\lambda))}, \\
\mathbi{G} & = - \frac{4 \alpha^3}{ (4 \beta_{1}(\lambda) - \beta_{4}(\lambda) ) (2\beta_{1}(\lambda) + \beta_{2}(\lambda) -\beta_{4}(\lambda)) (2 \beta_{2}(\lambda) -\beta_{4}(\lambda))},
\eea
and
\bea
\alpha = \frac{\gamma \pi}{\sqrt{2} l^{3/2}}.
\eea

\et

\br\label{Rem_h2better}
{\mk Note that the analytic expression of $h^{(2)}_\lambda$  given in \eqref{h2 expansion} can be written as the sum of  $h^{(1)}_\lambda$ given by \eqref{h1_Burgers}\footnote{Using the symbols introduced here, $h^{(1)}_\lambda(\xi_1,\xi_2) = \mathbi{A} \xi_1 \xi_2 e_3 + \mathbi{E} (\xi_2)^2 e_4$ from \eqref{h1_Burgers}.} associated with the one-layer backward-forward system \eqref{LLL}, and some other higher-order terms. It is worth noting that the extra five terms contained in the expression of $h^{(2)}_\lambda$ result from the nonlinear self-interactions between the low modes  as brought by $P_{\c} B \bigl( y^{(1)}_{\c},  y^{(1)}_{\c}\bigr)$ in \eqref{L2-uc2}. Numerical results of Section \ref{Sec_numresultsh2} below, support the fact that these extra terms can be interpreted as corrective terms to $h^{(1)}_\lambda$.  Indeed,  as we will illustrate for the optimal control problem \eqref{BP}, these terms can help design suboptimal low-dimensional controller of better performances than those built from $h^{(1)}_{\lambda}$-based reduced system; the $h^{(2)}_\lambda$-based reduced system bringing extra higher-order terms corresponding to ``low-high'' and ``high-high'' interactions absent from the $h^{(1)}_{\lambda}$-based reduced system. {\mkr This last point can be observed by comparing  \eqref{eq:Burgers reduced} with \eqref{eq:Burgers reduced-h2} below, where both reduced systems are derived from the abstract formulation \eqref{SEE2} by setting the PM function $h$ to be  $h^{(1)}_{\lambda}$ or $h^{(2)}_{\lambda}$,  respectively.}
}
\er

\bp

{\mk A simple integration of \eqref{L2} shows that} for any $\tau > 0$ and $\xi \in \mathcal{H}^{\c}$  the solution to the backward-forward system \eqref{L2} is given by:
\begin{subequations}
\begin{align}
y^{(1)}_{\c}(s) & = e^{s L_\lambda^{\c}}\xi,  \\
y^{(2)}_{\c}(s) & = e^{s L_\lambda^{\c}}\xi - \int_{s}^0 e^{(s-\tau') L_\lambda^{\c}} P_{\c} B\bigl(y^{(1)}_{\c}(\tau'), y^{(1)}_{\c}(\tau')\bigr) \d \tau',   \\
y_{\s}^{(2)}[\xi]{\mk(-\tau, s)}  & = \int_{-\tau}^s e^{(s-\tau') L_\lambda^{\s}} P_{\s} B\bigl(y^{(2)}_{\c}(\tau'), y^{(2)}_{\c}(\tau')\bigr) \d \tau', \label{ys-2}
\end{align}
\end{subequations}
for all  $s \in [-\tau, 0]$.

{\mk Due to \eqref{ys-2}, the pullback limit of $y_{\s}^{(2)}[\xi](-\tau,0)$ takes the form given in \eqref{eq:h2} provided that the concerned integral exists}. We show below that the \eqref{NR2}-condition is necessary and sufficient for  {\mk such an integral} to exist.  {\mk In that respect, the fact that $\mathcal{H}^{\c}$ is spanned by the first two eigenmodes facilitate some of the manipulations as described below}.

{\mk First, note that the projections of $y^{(1)}_{\c}$ onto $e_1$ and $e_2$, give respectively,}
\be   \label{yc1}
y^{(1)}_{1}(s)  := \langle y^{(1)}_{\c}(s), e_1 \rangle  = e^{\beta_1(\lambda) s}\xi_1, \qquad y^{(1)}_{2}(s)  := \langle y^{(1)}_{\c}(s), e_2 \rangle  = e^{\beta_2(\lambda) s}\xi_2,
\ee
where $\xi_i := \langle \xi, e_i \rangle$, $i = 1, 2$.

{\mk To determine the projection of $y^{(2)}_{\c}$ against $e_1$ and $e_2$, we need to recall that  the nonlinear interaction laws \eqref{nonlinear_interaction}, give here}
\be
B_{11}^1 = 0, \qquad B_{12}^1 =   2\alpha, \quad  B_{21}^1 = - \alpha, \quad   B_{11}^2 = - \alpha, \quad  B_{12}^2 = B_{21}^2 = 0,
\ee
which leads to 
\beas
\langle B(y^{(1)}_{\c}, y^{(1)}_{\c}), e_1 \rangle &  = \bigl \langle B\bigl(y^{(1)}_{1}e_1 + y^{(1)}_{2} e_2, y^{(1)}_{1}e_1 + y^{(1)}_{2} e_2 \bigr), e_1 \bigr \rangle \\
&   = y^{(1)}_{1} y^{(1)}_{2}B_{12}^1 - y^{(1)}_{1}y^{(1)}_{2} B_{21}^1 = \alpha y^{(1)}_{1}y^{(1)}_{2}, \\
\langle B(y^{(1)}_{\c}, y^{(1)}_{\c}), e_1 \rangle & = \bigl( y^{(1)}_{1} \bigr)^2 B_{11}^2 = - \alpha \bigl( y^{(1)}_{1} \bigr)^2.
\eeas

 {\mk The projection of $y^{(2)}_{\c}$ against $e_1$ and $e_2$ are then given by}
\bea  \label{yc2}
y^{(2)}_{1}(s)  & := \langle y^{(2)}_{\c}(s), e_1 \rangle  = e^{\beta_1(\lambda) s}\xi_1 - \alpha \int_s^0 e^{\beta_1(\lambda)(s-\tau')}y^{(1)}_{1}(\tau') y^{(1)}_{2}(\tau') \d \tau', \\
y^{(2)}_{2}(s)  & := \langle y^{(2)}_{\c}(s), e_2 \rangle  = e^{\beta_2(\lambda) s}\xi_2 + \alpha \int_s^0 e^{\beta_2(\lambda)(s-\tau')} (y^{(1)}_{1}(\tau'))^2 \d \tau'.
\eea 

{\mk Relying again on to the nonlinear interaction laws \eqref{nonlinear_interaction}}, we have 
\bea
B_{11}^3 &= 0, \qquad B_{12}^3 =  - 2\alpha, \qquad  B_{21}^3 = - \alpha, \qquad B_{22}^3 = 0, \\
B_{11}^4  &= B_{12}^4 = B_{21}^4 = 0, \hspace{3.6em} B_{22}^4 =  - 2 \alpha, \\
B_{ij}^n &= 0, \qquad \Forall i, j \in \{1,2\}, \,n \ge 5,
\eea
which leads to
\bea
y^{(2)}_{3}[\xi]{\mk(-\tau, s)}  & := \langle y^{(2)}_{\s}[\xi]{\mk(-\tau, s)} , e_3 \rangle  =  - 3 \alpha \int_{-\tau}^s e^{\beta_3(\lambda)(s-\tau')} y^{(2)}_{1}(\tau') y^{(2)}_{2}(\tau') \d \tau', \\
y^{(2)}_{4}[\xi]{\mk(-\tau, s)} & := \langle y^{(2)}_{\s}[\xi]{\mk(-\tau, s)} , e_4 \rangle  = - 2\alpha \int_{-\tau}^s e^{\beta_4(\lambda)(s-\tau')} (y^{(2)}_{2}(\tau'))^2 \d \tau'.
\eea

By using the expressions of $y^{(2)}_{1}$ and $y^{(2)}_{2}$ given in \eqref {yc2}  (and using also \eqref{yc1}), it can be checked that the limit $h^{(2),3}_\lambda := \lim_{\tau \rightarrow +\infty} y^{(2)}_{3}[\xi]{\mk(-\tau, 0)} $ exists if and only if the first four inequalities in the \eqref{NR2}-condition hold, {\mk while}  $h^{(2),3}_\lambda$ is given by \eqref{h23} under these conditions. Similarly, the limit $h^{(2),4}_\lambda := \lim_{\tau \rightarrow +\infty} y^{(2)}_{4}[\xi]{\mk(-\tau, 0)} $ exists if and only if the last three inequalities in the \eqref{NR2}-condition hold, and  $h^{(2),4}_\lambda$ is given by \eqref{h24} under these conditions. The theorem is proved.

\ep

\subsection{Controller synthesis based on $h^{(2)}_\lambda$, and control performances: {\mk Analytic derivation and  numerical results}}\label{Sec_numresultsh2}

{\mk  {\bf Analytic derivation of the $h^{(2)}_\lambda$-based reduced optimal control problem.} Following \eqref{SEE2}, the $h^{(2)}_\lambda$-based reduced system intended to model the dynamics of the low modes $P_{\c}y$ of \eqref{eq:Burgers}, takes the following abstract form:}
\bea \label{SEE3_Burgers}
& \frac{\d z}{\d t} = L^{\c}_\lambda z + P_{\c} B\Bigl(z + h^{(2)}_\lambda(z), z + h^{(2)}_\lambda(z)\Bigr) +  \mathfrak{C} \ur,  \qquad {\hl t \in (0, T]}, \\ & z(0)  = P_{\c} y_0 \in \mathcal{H}^{\c}, 
\eea
where $y_0$ is the initial datum for the original PDE \eqref{eq:Burgers}. 

{\hl Analogous} to \eqref{J2_Burgers}, the cost functional {\hl associated with the reduced system \eqref{SEE3_Burgers}} is given by
\be \label{J3_Burgers}
\widehat{J}_{R}(z, \ur) = \int_0^T  \bigl( \frac{1}{2}\|z(t) + h^{(2)}_\lambda(z(t))\|^2  + \frac{\mu_1}{2} \|\ur(t)\|^2 \bigr) \d t +  C_T(z(T), P_{\c}Y),
\ee
where $C_T(z(T), P_{\c} Y)  := \frac{\mu_2}{2} \sum_{i=1}^m |z_i(T) - Y_i|^2$ is the terminal payoff term as defined in \eqref{C_T_sect5}, with $Y$ being {\mk some} prescribed target for \eqref{eq:Burgers}.

By using the analytic expression of $h^{(2)}_\lambda$ given  {\mk in \eqref{h2 expansion}-\eqref{coeffs_h2},} the cost functional \eqref{J3_Burgers} can be written {\mk into the following explicit form}:
\be \label{J3_Burgers-b}
\widehat{J}_R(z, \ur) = \int_0^T \Bigl[ \frac{1}{2} \mathcal{G}(z(t))  +  \frac{\mu_1}{2} \mathcal{E}(\ur(t)) \Bigr] \d t + C_T(z(T), P_{\c} Y),
\ee
where
\bea
\mathcal{G}(z) &=  \frac{1}{2} \|z + h^{(2)}_\lambda(z)\|^2 = \frac{1}{2} \Bigl[ (z_1)^2 + (z_2)^2 + (h^{(2),3}_\lambda(\xi_1,\xi_2))^2 + (h^{(2),4}_\lambda(\xi_1,\xi_2))^2 \Bigr], \\
\mathcal{E}(\ur) & =  \frac{\mu_1}{2} \|\ur\|^2 =  \frac{\mu_1}{2} [(\urc{1})^2 + (\urc{2})^2],
\eea
with $z_i := \langle z, e_i\rangle$ and  $\urc{i} := \langle \ur, e_i\rangle$, $i = 1, 2$. 

Now, by using again the analytic expression 
\bes
h^{(2)}_\lambda(\xi_1e_1+\xi_2e_2) = h^{(2),3}_\lambda(\xi_1,\xi_2) e_3 + h^{(2),4}_\lambda(\xi_1,\xi_2) e_4
\ees
 in \eqref{SEE3_Burgers} and projecting this equation against $e_1$ and $e_2$ respectively, we obtain, after simplification {\hl by} using the nonlinear interaction {\mk laws} \eqref{nonlinear_interaction}, {\mk the following analytic formulation of the $h^{(2)}_\lambda$-based reduced system \eqref{SEE3_Burgers}}:
\begin{equation}\label{eq:Burgers reduced-h2}
\boxed{
\begin{aligned}
 &  \frac{\d z_1}{\d t}  = \beta_1(\lambda) z_1 + \alpha \Bigl( z_1z_2 +  z_2 h^{(2),3}_\lambda(z_1,z_2)  + h^{(2),3}_\lambda(z_1,z_2) h^{(2),4}_\lambda(z_1,z_2)\Bigr)  + a_{11}\urc{1}(t) + a_{21}\urc{2}(t), \\
&  \frac{\d z_2}{\d t} = \beta_2(\lambda) z_2 - \alpha z_1^2   + 2 \alpha \Bigl( z_1 h^{(2),3}_\lambda(z_1,z_2) + z_2 h^{(2),4}_\lambda(z_1,z_2) \Bigr)  + a_{12}\urc{1}(t) + a_{22}\urc{2}(t),
\end{aligned}
}
\end{equation}
{\mk with $h^{(2),3}_\lambda(z_1,z_2)$ and $h^{(2),4}_\lambda(z_1,z_2)$ given by \eqref{coeffs_h2_a}-\eqref{coeffs_h2}.\footnote{{\mk Using this analytic formulation, we mention that the Cauchy problem for \eqref{eq:Burgers reduced-h2} 
can be dealt with by carrying out similar (but more tedious) energy estimates  as presented in Appendix~\ref{Sect_energy_est} for the two-dimensional $h^{(1)}_\lambda$-based reduced system \eqref{eq:Burgers reduced}}.}}

The resulting reduced optimal control problem based on $h^{(2)}_\lambda$ is {\hl thus}:
\be  \label{R2BP} 
\begin{aligned}
 \hspace{-1em} \min \widehat{J}_R(z, \ur)  \quad \! \text{ s.t. } \quad\!  (z, \ur) \in L^2(0,T; \mathcal{H}^{\c}) \times L^2(0,T; \mathcal{H}^{\c})  \quad\! \text{solves  \eqref{eq:Burgers reduced-h2}}.
\end{aligned}
\ee

By {\mk following similar arguments as provided in Section \ref{ss:reduction-h1}  and applying the Pontryagin maximum Principle, we can conclude that 
for a given pair
\bes
(\widehat{z}_R^\ast, \widehat{u}_{R}^\ast)  \in L^2(0,T; \mathcal{H}^{\c}) \times L^2(0,T; \mathcal{H}^{\c}) 
\ees
to be optimal for the $h^{(2)}_\lambda$-reduced optimal problem \eqref{R2BP}}, it is necessary and sufficient\footnote{{\mk The sufficient part is again due to the fact that the cost functional \eqref{J3_Burgers} is quadratic in $u_R$ and the dependence on the controller is affine for the system of equations \eqref{eq:Burgers reduced-h2}; see {\it e.g.}~\cite[Sect.~5.3]{Kirk12} and \cite{Trelat2012}.}} to satisfy the following {\hl set of} conditions:
\bea  \label{uc-p-h2}
\boxed{(\widehat{u}^\ast_{R,1}, \widehat{u}^\ast_{R,2})  = - \Bigl( \frac{a_{11} \widehat{p}_{R,1}^\ast + a_{12} \widehat{p}_{R,2}^\ast}{\mu_1},  \frac{a_{21} \widehat{p}_{R,1}^\ast +  a_{22} \widehat{p}_{R,2}^\ast}{\mu_1} \Bigr ),}
\eea
where $(\widehat{p}_{R,1}^\ast, \widehat{p}_{R,2}^\ast)$ is the costate associated with $(\widehat{z}_{R,1}^\ast, \widehat{z}_{R,1}^\ast)$, {\mk both determined by solving the following BVP}:
\bea \label{eq:bvp-h2}
&  \frac{\d z_1}{\d t} =  \beta_1(\lambda) z_1 + \alpha \Bigl( z_1z_2 +  z_2 h^{(2),3}_\lambda(z_1,z_2)  + h^{(2),3}_\lambda(z_1,z_2) h^{(2),4}_\lambda(z_1,z_2)\Bigr)   - \frac{1}{2} p_1, \\
& \frac{\d z_2}{\d t} = \beta_2(\lambda) z_2 - \alpha (z_1)^2   + 2 \alpha \Bigl( z_1 h^{(2),3}_\lambda(z_1,z_2) + z_2 h^{(2),4}_\lambda(z_1,z_2) \Bigr)- \frac{1}{2} p_2, \\
& \frac{\d p_1}{\d t} =   g_1(z, p),  \\
& \frac{\d p_2}{\d t} =   g_2(z, p),
\eea
subject to the boundary condition
\be \label{eq:bvp h2}
z_1(0) = \langle y_0, e_1 \rangle, \quad z_2(0) = \langle y_0, e_2 \rangle, \quad p_1(T) = \mu_2 (z_{1}(T) -  Y_1), \quad p_2(T) = \mu_2 (z_2(T) -  Y_2),
\ee
where
\beas
g_1(z, p) & := - z_1 -  h^{(2),3}_\lambda(z_1,z_2) \frac{\partial h^{(2),3}_\lambda(z_1,z_2)}{\partial z_1} -  h^{(2),4}_\lambda(z_1,z_2) \frac{\partial h^{(2),4}_\lambda(z_1,z_2)}{\partial z_1} \\
& \hspace{2em} -  p_1 \biggl( \beta_1(\lambda)   + \alpha z_2 +  \alpha z_2 \frac{\partial h^{(2),3}_\lambda(z_1,z_2)}{\partial z_1} + \alpha \frac{\partial h^{(2),3}_\lambda(z_1,z_2)}{\partial z_1}   h^{(2),4}_\lambda(z_1,z_2) \\
& \hspace{5em}  + \alpha  h^{(2),3}_\lambda(z_1,z_2) \frac{\partial h^{(2),4}_\lambda(z_1,z_2)}{\partial z_1} \biggr) \\
& \hspace{2em}  -  2  \alpha p_2 \biggl( -  z_1 + h^{(2),3}_\lambda(z_1,z_2)  + z_1 \frac{\partial h^{(2),3}_\lambda(z_1,z_2)}{\partial z_1} + z_2 \frac{\partial h^{(2),4}_\lambda(z_1,z_2)}{\partial z_1} \biggr), \\
g_2(z, p) & := -  z_2 -  h^{(2),3}_\lambda(z_1,z_2) \frac{\partial h^{(2),3}_\lambda(z_1,z_2)}{\partial z_2} - h^{(2),4}_\lambda(z_1,z_2) \frac{\partial h^{(2),4}_\lambda(z_1,z_2)}{\partial z_2} \\
& \hspace{0em} - \alpha  p_1 \biggl( z_1 + h^{(2),3}_\lambda(z_1,z_2) + z_2 \frac{\partial h^{(2),3}_\lambda(z_1,z_2)}{\partial z_2} + \frac{\partial h^{(2),3}_\lambda(z_1,z_2)}{\partial z_2}   h^{(2),4}_\lambda(z_1,z_2) \\
& \hspace{5em}  +  h^{(2),3}_\lambda(z_1,z_2) \frac{\partial h^{(2),4}_\lambda(z_1,z_2)}{\partial z_2} \biggr) \\
& \hspace{0em}  -  p_2 \biggl(  \beta_2(\lambda) + 2 \alpha z_1 \frac{\partial h^{(2),3}_\lambda(z_1,z_2)}{\partial z_2}
+ 2 \alpha h^{(2),4}_\lambda(z_1,z_2)  + 2 \alpha z_2 \frac{\partial h^{(2),4}_\lambda(z_1,z_2)}{\partial z_2}\biggr).
\eeas

{\mk The vector  field $(g_1,g_2)$ given above has been determined by evaluating $-\nabla_z \widehat{H}(z,p,u)$, with the following  Hamiltonian $\widehat{H}$, formed by application of the PMP to \eqref{R2BP}
$$\widehat{H}(z, p, u)  := \mathcal{G}(z) + \mathcal{E}(u) + p_1 \widehat{f_1}(z, u) + p_2 \widehat{f_2}(z, u),$$
where $(\widehat{f_1}, \widehat{f_2})$ denotes the vector field constituting the RHS of  the $z$-equations in \eqref{eq:bvp-h2}.
}

\medskip

 {\mk {\bf Numerical results.}} The above BVP is solved again using \texttt{bvp4c},  and the {\mk resulting} two-dimensional {\mk suboptimal controller $\widehat{u}^\ast_{R}$  is obtained according to \eqref{uc-p-h2}}. {\mk As before, the corresponding suboptimal trajectory $\widehat{y}_R^\ast$ of the PDE \eqref{eq:Burgers} is computed by driving \eqref{eq:Burgers} with $\widehat{u}^\ast_{R}$, following the numerical procedure described in Section \ref{Sec_numaspects}.}

\begin{figure}[!hbtp]
\centering
\includegraphics[height=0.5\textwidth, width=0.8\textwidth]{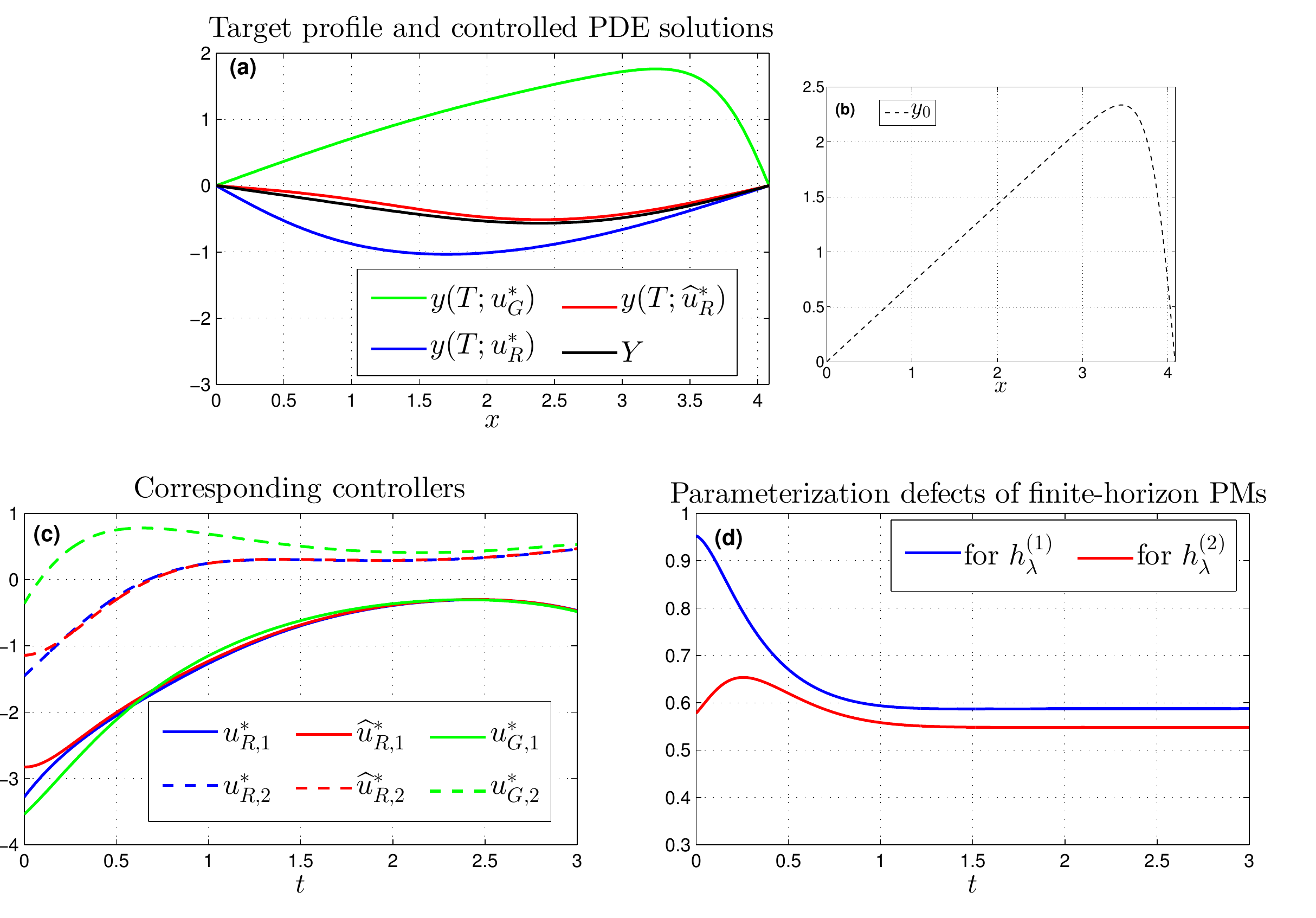}
\vspace{-0.5em}
\caption{{\footnotesize {\bf (a)}: Final state at $T = 3$ of the PDE solution profiles driven respectively by the suboptimal controllers $u^\ast_{G}$, $\ur^\ast$ and $\widehat{u}^\ast_{R}$ with initial profile taken to be $y^+$ as shown in {\bf (b)}; also shown in {\bf (a)} is the target state {\mk $Y$ given by \eqref{Y_target2}}.  {\bf (c)}: The suboptimal controllers $u^\ast_{G}= u^\ast_{G,1} e_1 + u^\ast_{G,2} e_2$ synthesized from the Galerkin-based reduced optimal control problem \eqref{GBP}; $\ur^\ast= \urc{1}^\ast e_1 + \urc{2}^\ast e_2$ synthesized from the $h^{(1)}_\lambda$-based reduced optimal control problem \eqref{RBP}; and $\widehat{u}^\ast_R=\widehat{u}^\ast_{R,1} e_1 + \widehat{u}^\ast_{R,2} e_2$ synthesized from the $h^{(2)}_\lambda$-based reduced optimal control problem \eqref{R2BP}. {\bf (d)}: Finite-horizon parameterization defects of $h^{(1)}_\lambda$ and $h^{(2)}_\lambda$ associated with the PDE \eqref{eq:Burgers} driven respectively by $\ur^\ast$ and $\widehat{u}^\ast_R$ over the time interval $[0, 3]$. The system parameters are the same as in Section~\ref{Sect_Burgers}; see caption of Fig.~\ref{fig:state_Sect5}.
}} \label{fig:state_Sect6}
\end{figure}
 
The corresponding control performance is shown in Fig.~\ref{fig:state_Sect6}, where the performance of the suboptimal controllers $u^\ast_{R}$ and $u^\ast_{G}$ associated with respectively the two-dimensional $h^{(1)}_\lambda$-based reduced optimal control problem \eqref{RBP} and the {\mk two-dimensional Galerkin-based one} \eqref{GBP} are also reported for comparison. In panel {\bf (a)} of Fig.~\ref{fig:state_Sect6}, {\mk we present the PDE final time solution profile $y(T,\widehat{u}_R^\ast)$, $y(T,\ur^\ast)$, and $y(T,\u_G^\ast)$ driven  respectively by $\widehat{u}_R^\ast$, $\ur^\ast$ and $\u_G^\ast$, for $T=3$.}
{\mk For these simulations, the target profile $Y$ has been  chosen to be again spanned by the first two eigenfunctions, but given this time by
\be\label{Y_target2}
Y = -0.3\langle y^+, e_1\rangle e_1  - 0.1 \langle y^+, e_2 \rangle e_2;
\ee
the initial profile is taken to be the positive steady state $y^+$ for the uncontrolled PDE as used in Section \ref{Sec_numresults}, see panel {\bf (b)}.}
{\mk The two components of the synthesized suboptimal controllers are shown in panel {\bf (c)},  and the parameterization defects associated with respectively $h^{(1)}_\lambda$ and $h^{(2)}_\lambda$ are shown in panel {\bf (d)}. The corresponding cost {\mkr values} and final-time relative $L^2$-errors are given in Table~\ref{tab:cost_Sect6} below.}

\begin{table} [!h]   
\center
\caption[c]{Cost {\mkr values} and final-time relative $L^2$-errors associated with the suboptimal controllers}   \label{tab:cost_Sect6}
\begin{tabular}{ |c | c  c c  c| }  
  \hline                       
      &  $u_G^\ast$   &  $u_R^\ast$   & $\widehat{u}_R^\ast$ &   $\widetilde{u}_G^\ast$ (with $m=16$)\\
\hline 
$J(y(\cdot; y^+, u), u)$  &    108.65  &  12.48  & 5.07   & 5.02  \\
\hline 
Relative $L^2$-error: $\|y(T; y^+, u) - Y\|/\|Y\|$ &  $405.60\%$ &   $107.23\%$   &  $15.07\%$ & $11.41\%$ \\
\hline
\end{tabular}
\vspace{1ex}
\begin{quote}
The cost $J$ is {\mk  the one defined in \eqref{JJ} associated with the optimal control problem \eqref{BP}.}  {\mk This cost is assessed} at the suboptimal pairs $(y(\cdot; y^+, u), u)$ with $u$ taken to be {\mk either} $u_G^\ast$, $u_R^\ast$, $\widehat{u}_R^\ast$, or $\widetilde{u}_G^\ast$. {\mk  The target $Y$ is given by \eqref{Y_target2}.} The suboptimal controller $u_G^\ast$ is synthesized from the {\mk 2D Galerkin-based reduced optimal control problem \eqref{GBP}; $u_R^\ast$ from the $h^{(1)}_\lambda$-based \eqref{RBP}; $\widehat{u}_R^\ast$ from the $h^{(2)}_\lambda$-based \eqref{R2BP}; and $\widetilde{u}_G^\ast$ from the $m$-dimensional Galerkin-based one \eqref{GBP'} with $m=16$.}  The latter  serves as a benchmark here. {\mk The model parameters are those used for Fig.~\ref{fig:state_Sect5}.}
\end{quote}
\end{table}

{\mk The results of Fig.~\ref{fig:state_Sect6} {\bf (a)} and Table~\ref{tab:cost_Sect6} illustrate that for a given reduced phase space \textemdash\, here the two-dimensional vector space $\mathcal{H}^\c$ \textemdash\, the slaving relationship of the high-modes (not in $\mathcal{H}^\c$) by the low modes  (in $\mathcal{H}^\c$) as parameterized by $h^{(2)}_\lambda$ can turn out to be superior than the one proposed by $h^{(1)}_\lambda$ for the synthesis of suboptimal solutions to \eqref{BP}, and can turn out to be clearly advantageous compared to suboptimal solutions for which no  slaving relationship whatsoever is involved such as for  those built from the 2D Galerkin-based reduced optimal control problem \eqref{GBP}. Again, Corollary~\ref{Cor_2} and the error estimate \eqref{cor2:goal} provide theoretical insights  that help understand why improving the quality of such a slaving relationship participates to improve the performance of a suboptimal controller. For instance, the improvement in getting closer to the prescribed target $Y$ (Fig.~\ref{fig:state_Sect6} {\bf (d)}) \textemdash\, accompanied with a noticeable reduction of the cost {\mkr values} (Table \ref{tab:cost_Sect6}) \textemdash\,  occurs when the PDE \eqref{eq:Burgers} is driven  by the $h^{(2)}_\lambda$-based suboptimal controller $\widehat{u}^\ast_{R}$ instead of the $h^{(1)}_\lambda$-based one $\ur^\ast$, and goes with a parameterization defect (overall) smaller for $h^{(2)}_\lambda$ than for $h^{(1)}_\lambda$ (Fig.~\ref{fig:state_Sect6} {\bf (d)}). Interestingly, this reduction of the parameterization defect {\mkr comes} with {\mkr the} higher-order terms contained in $h^{(2)}_{\lambda}$ (see Theorem \ref{THM_h2}) that can be thus {\mkr reasonably} interpreted as correction terms to the parameterization proposed  by  $h^{(1)}_{\lambda}$;  see {\mkr also} Remark~\ref{Rem_h2better}.}

However, such a statement has to be nuanced and an $h^{(2)}_{\lambda}$-based reduced system does not always lead to the significant advantages in the design of suboptimal solutions such as illustrated in Fig.~\ref{fig:state_Sect6}. The caveat relies on the fact that the parameterization defect associated with $h^{(2)}_\lambda$ also depends on the target profile. For instance, with the sign-changing target \eqref{Y_target} used in the experiments of Section \ref{Sec_numresults}, the suboptimal solutions designed from \eqref{R2BP} achieve {\HL comparable performances to} those designed from \eqref{RBP}.

{\mk {\mkr These remarks} motivate further analysis to arbitrate whether the success achieved for the target prescribed in \eqref{Y_target2} are pathological or robust, to some extent. For that purpose, we considered {\mkr deformations} of the target \eqref{Y_target2} taken to be of the form 
\be\label{Y_target3}
Y_{\sigma_1,\sigma_2}= - \sigma_1 \langle y^+, e_1 \rangle e_1 - \sigma_2 \langle y^+, e_2 \rangle e_2,
\ee
 with $\sigma_1 \in [0.2, 0.7]$ and $\sigma_2 \in [0.01, 0.5]$, and we {\mkr solved} the corresponding $h^{(2)}_{\lambda}$-based (resp. $h^{(1)}_{\lambda}$-based) reduced optimal problem to provide  the corresponding $h^{(2)}_{\lambda}$-based (resp. $h^{(1)}_{\lambda}$-based) suboptimal solutions. As a benchmark\footnote{{\mkr Here, $4$ significant digits of the cost $J$ are ensured with $m=16$ by comparing with cost values associated with higher-dimensional suboptimal controller synthesized from \eqref{GBP'}}.}, these solutions are compared with those obtained  from the $m$-dimensional Galerkin-based reduced optimal problem \eqref{GBP'} with $m=16$.
 The results are reported in Fig.~\ref{fig:contour} and in Fig.~\ref{fig:J} below.  Figure \ref{fig:contour}  shows for each $(\sigma_1,\sigma_2)$ the corresponding relative $L^2$-errors at the final-time solution profiles compared with the target $Y_{\sigma_1,\sigma_2}$; and Figure \ref{fig:J} shows the cost {\mkr values} associated with the suboptimal controllers $u_R^\ast$ and $\widehat{u}_R^\ast$, on one hand, and $\widetilde{u}_G^\ast$ obtained from the $m$-dimensional Galerkin-based reduced problem, on the other.}

\begin{figure}[!hbtp]
\centering
\includegraphics[height=0.6\textwidth, width=.8\textwidth]{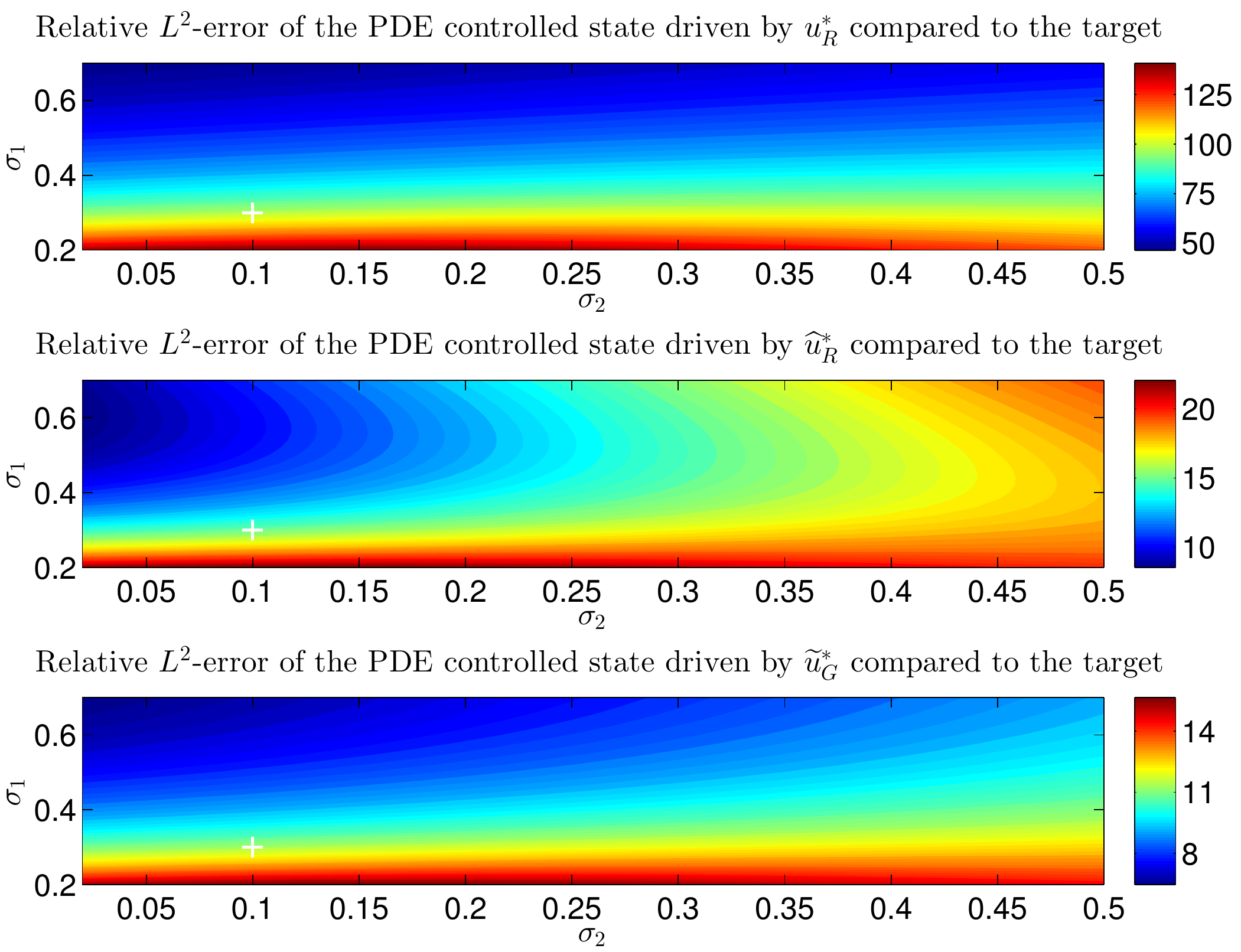}
\caption{{\footnotesize {\mk $(\sigma_1,\sigma_2)$-dependence of the relative $L^2$-error of the PDE final state $y(T,y^+; u)$ compared to the target $Y_{\sigma_1,\sigma_2}$ given by \eqref{Y_target3}}. Here the controller $u$ is taken to be either $u_{R}^\ast$ ({\bf upper panel}), or $\widehat{\u}_R^\ast$ ({\bf middle panel}), or $\widetilde{u}^\ast_G$ ({\bf lower panel}); the parameters $\sigma_1$ and $\sigma_2$ are taken to be $\sigma_1 \in [0.2, 0.7]$ and $\sigma_2 \in [0, 0.5]$; and the final time is $T = 3$. The markers ``+'' in the plots correspond {\mk to the results shown in Fig.~\ref{fig:state_Sect6}, for  $(\sigma_1, \sigma_2) = (0.3,0.1)$.}}} \label{fig:contour}
\end{figure}

{\mk
Figures \ref{fig:contour} and \ref{fig:J} show that the good performance achieved by the $h^{(2)}_\lambda$-based suboptimal controller shown in Fig.~\ref{fig:state_Sect6}~{\bf (a)}, is not isolated and can be even further improved  within a broad region of the $(\sigma_1,\sigma_2)$-parameter space when $Y_{\sigma_1,\sigma_2}$ is changed accordingly.  Compared to the bad performances observed on Fig.~\ref{fig:contour} (top panel) for the $h^{(1)}_\lambda$-based suboptimal controllers, these $h^{(2)}_\lambda$-based results provide strong evidence that the higher-order terms brought by $h^{(2)}_\lambda$ with respect to  $h^{(1)}_\lambda$, act as corrective terms in the {\mkr high-mode parametrization proposed by $h^{(1)}_\lambda$.}}  

These numerical results together with the theoretic results of Corollary~\ref{Cor_2} suggest that in order to design reduced problems {\mkr whose solutions would provide even better control performance than those reported here}, one can try to construct finite-horizon PMs with smaller {\mk parameterization} defects {\mkr than those achieved by $h^{2)}_\lambda$. In that}} respect, the {\mkr discussions and results of \cite[Sect.~8.3-8.5]{CLW13a}, presented in the context of asymptotic PMs}, can be valuable. 
In connection {\mk to the discussion concerning Figs.~\ref{fig:PM_Sect5} and \ref{fig:Energy_Sect5} in} Section~\ref{Sec_numresults}, the searching for better slaving relationships between the $\mathcal{H}^\s$-modes and the $\mathcal{H}^\c$-modes can be combined with the usage of higher dimensional {\mk reduced phase spaces} $\mathcal{H}^\c$ so that the energy kept in the high modes gets reduced. {\mk The next section shows that a moderate increase of $\mbox{dim}(\mathcal{H}^\c)$ can actually already help improve the performances based on $h^{(1)}_\lambda$, in the case of locally distributed control laws.}

\begin{figure}[!hbtp]
\centering
\includegraphics[height=0.5\textwidth, width=0.5\textwidth]{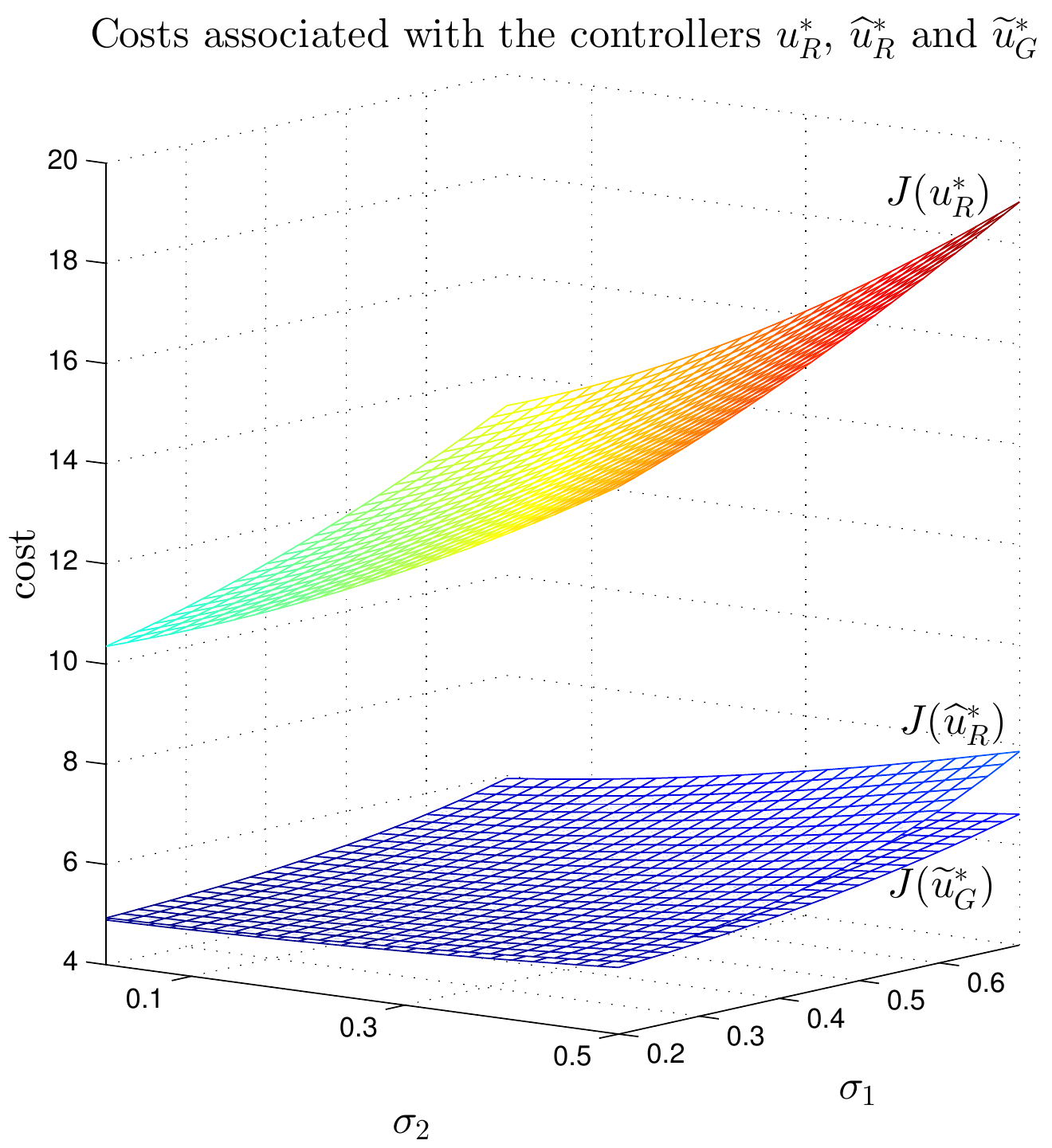}
\caption{{\footnotesize {\mk $(\sigma_1,\sigma_2)$-dependence of the cost values $J(y,u)$ given by \eqref{JJ} when $u=u_R^\ast$, $u=\widehat{u}^\ast_R$ and $u=\widetilde{u}^\ast_G$, respectively. The parameters $\sigma_1$ and $\sigma_2$ vary  in $[0.2, 0.7]$ and in $[0, 0.5]$, respectively. The final time is still $T = 3$.}}} \label{fig:J}
\end{figure}

\section{Synthesis of $m$-Dimensional Locally Distributed Suboptimal Controllers} \label{Sect_Burgers_local}

In this last section, we consider the more challenging case of optimal {\mkr locally distributed  control} problems associated with the {\mk Burgers-type equation \eqref{eq:Burgers}.} {\mk This situation corresponds to the case where the linear operator $\mathfrak{C}$ is associated with the characteristic function $\chi_\Omega$ of a subdomain $\Omega \subset [0, l]$, such that for any $u \in \mathcal{H} = L^2(0,l)$, the action of $\mathfrak{C}$ on $u$ is defined  by:}
\be
\mathfrak{C} u(x) = \chi_\Omega(x) u(x), \quad \Forall  x \in [0,l].
\ee

{\mk As used in the fully distributed case in the previous sections, we will consider for some prescribed (time-independent) target $Y$,} cost functionals of {\it terminal payoff type} {\mk such as:}
\be \label{JJ_terminal}
J^{\mathrm{TP}}(y, u) = \int_0^T \bigl( \frac{1}{2} \|y(t; y_0, u)\|^2  + \frac{\mu_1}{2}\|u(t)\|^2 \bigr) \d t + \frac{\mu_2}{2}\|y(T; y_0, u) - Y\|^2,
\ee
but {\mk also cost functionals} of {\it tracking type}:
\be \label{JJ_tracking}
J^{\mathrm{track}}(y, u) = \int_0^T \bigl( \frac{1}{2} \|y(t; y_0, u) - Y\|^2  + \frac{\mu_1}{2}\|u(t)\|^2 \bigr) \d t,
\ee
where {\hl in both cases}, $\mu_1$ and $\mu_2$ are some positive {\mk parameters}.

The optimal control problem {\hl takes thus one of the following forms}:
\be  \label{LBP}   
\begin{aligned}
 & \hspace{5em} \min J^{\mathrm{TP}}(y, u)  \quad \text{with $J^{\mathrm{TP}}$ defined in \eqref{JJ_terminal}} \qquad \text{s.t.} \\
&  (y, u) \in L^2(0,T; \mathcal{H}) \times L^2(0,T; \mathcal{H}) \text{ solves the problem \eqref{eq:Burgers}--\eqref{initial:Burgers}}. 
\end{aligned}
\ee
or
\be  \label{LBP'} 
\begin{aligned}
 & \hspace{5em} \min J^{\mathrm{track}}(y, u)  \quad \text{with $J^{\mathrm{track}}$ defined in \eqref{JJ_tracking}} \qquad \text{s.t.} \\
&  (y, u) \in L^2(0,T; \mathcal{H}) \times L^2(0,T; \mathcal{H}) \text{ solves the problem \eqref{eq:Burgers}--\eqref{initial:Burgers}}. 
\end{aligned}
\ee

{\mk The goal of this last section is to show that the PM-approach introduced above provides an efficient way to design suboptimal solutions for such optimal control problems associated with locally distributed control laws. For simplicity, we will focus on the performance achieved by the $h^{(1)}_\lambda$-based reduced system for the design of such suboptimal solutions, that is the following $m$-dimensional reduced system
\bea \label{reduced_Burgers_local}
& \frac{\d z}{\d t} = L^{\c}_\lambda z + P_{\c} B\Bigl(z + h^{(1)}_\lambda(z), z + h^{(1)}_\lambda(z)\Bigr) + P_{\c} \chi_{\Omega} \ur(t),  \qquad {\hl t \in (0, T],}
\eea
will be at the core of our synthesis of suboptimal controllers.

It is worthwhile to note that in general, the choice of the reduced dimension, $m$, depends typically on the system parameters such as the viscosity $\nu$, the domain size $l$ and the control parameter $\lambda$; and $m$ is chosen so that the resolved modes explain a sufficient large portion of the energy contained in the PDE solution. For the particular case of locally distributed control laws, the size and the location of the subdomain $\Omega$ plays also a determining role in sizing ``a good'' $m$. For instance, the smaller the subdomain $\Omega$ will be, the larger the dimension $m$ will need to be in order to obtain a reduced system useful for the design of good suboptimal controllers. Intuitively, this is related to the fact that {\mk further} eigenmodes are needed in order to obtain a reasonably good approximation of the characteristic function $\chi_\Omega$ when the size of the support $\Omega$ is {\mk further} reduced.  This intuition will be numerically {\mk confirmed} in Section \ref{Sec_local_num} below, where a reduction of 40 percent of the domain compared to the globally distributed case analyzed in Section \ref{Sec_numresults}, {\mkr led to a choice of $m=4$ for a} design of suboptimal controllers with comparable performances {\mk than those achieved in Section \ref{Sec_numresults}, from two-dimensional reduced systems}.
}

{\mk We now describe the $h^{(1)}_\lambda$-based reduced optimal control that will serve us to design the corresponding suboptimal controllers. First, note that the} cost functional associated with \eqref{reduced_Burgers_local} {\mk takes one of the following forms} 
\be \label{JJ_reduced_terminal}
J^{\mathrm{TP}}_R(z, \ur) = \int_0^T \Big(\frac{1}{2} \|z + h^{(1)}_\lambda(z)\|^2  +  \frac{\mu_1}{2} \|\ur\|^2  \Big) \d t + \frac{\mu_2}{2}\|z(T; z_0, \ur) - P_{\c}Y\|,
\ee
or
\be \label{JJ_reduced_tracking}
J^{\mathrm{track}}_R(z, \ur) = \int_0^T \Big(\frac{1}{2} \|z + h^{(1)}_\lambda(z) - Y\|^2  +  \frac{\mu_1}{2} \|\ur\|^2  \Big)\d t,
\ee
{\hl depending on whether \eqref{JJ_terminal} or \eqref{JJ_tracking} is considered.}

The reduced optimal control problem for \eqref{LBP} reads then as follows:
\be  \label{LRBP}   
\begin{aligned}
 \hspace{-1em} \min J^{\mathrm{TP}}_R(z, \ur)  \quad \! \text{ s.t. } \quad\!  (z, \ur) \in L^2(0,T; \mathcal{H}^{\c}) \times L^2(0,T; \mathcal{H}^{\c})  \quad\! \text{solves} \quad \!  \text{\eqref{reduced_Burgers_local}}.
\end{aligned}
\ee
Accordingly,  the reduced optimal control problem for \eqref{LBP'} reads:
\be  \label{LRBP'} 
\begin{aligned}
 \hspace{-1em} \min J^{\mathrm{track}}_R(z, \ur)  \quad \! \text{ s.t. } \quad\!  (z, \ur) \in L^2(0,T; \mathcal{H}^{\c}) \times L^2(0,T; \mathcal{H}^{\c})  \quad\! \text{solves} \quad\!  \text{\eqref{reduced_Burgers_local}}.
\end{aligned}
\ee

\subsection{Analytic derivation of $m$-dimensional $h_\lambda^{(1)}$-based reduced systems for the design of suboptimal controllers}  \label{ss:derivation_local}

{\mk In this subsection, we derive} explicit forms of the reduced suboptimal control problems \eqref{LRBP} and \eqref{LRBP'}. {\hl Details are presented for \eqref{LRBP}, while the analogous derivation for \eqref{LRBP'} is left to the interested reader. For this purpose,} let us first examine the existence of the finite-horizon PM candidate $h^{(1)}_\lambda$. We know from Section~\ref{ss:h1} that the pullback limit $h^{(1)}_\lambda$ associated with the backward-forward system \eqref{LLL} exists {\hl when} the \eqref{NR}-condition holds. For the Burgers equation considered here, {\mkr due} to the nonlinear interaction relations \eqref{nonlinear_interaction}, the \eqref{NR}-condition reads as follows:
\be  \label{NR-Burgers_local}
\Forall \, n > m, \ \Forall  i \in \{1, \cdots, m\},  \Big(n- i \in \{1, \cdots, m\} \Big)\Longrightarrow   \Big(\beta_{i}(\lambda) +  \beta_{n-i}(\lambda) - \beta_n(\lambda) > 0 \Big).
\ee

By using the analytic {\mkr expression} of the eigenvalues {\mkr as} given in \eqref{eq_eigenvalues}, we get
\be
\beta_{i}(\lambda) +  \beta_{n-i}(\lambda) - \beta_n(\lambda) = \lambda + \frac{\nu \pi^2 (n^2 - i^2 - (n-i)^2)}{l^2},
\ee
which is positive for all values of $\lambda$ of interest here ( $\lambda > {\hl \lambda_c := \frac{\nu \pi^2}{l^2}}$). Consequently, the pullback limit $h^{(1)}_\lambda$ always exists for such given $\lambda$, and its analytic form provided in \eqref{h1} reads as follows for the problem considered here:
\be  \label{h1_part1}
h^{(1)}_\lambda(\xi)  = \sum_{n > m} h^{(1),n}_\lambda(\xi) e_n,
\ee
where 
\be  \label{h1_part2}
\boxed{h^{(1),n}_\lambda(\xi) =  \sum_{ \substack{ i_1 + i_2 = n \\ 1\le i_1, i_2 \le m}}  \frac{\xi_{i_1} \xi_{i_2}}{\beta_{i_1}(\lambda) + \beta_{i_2}(\lambda) - \beta_n(\lambda)} \Bigl \langle  B(e_{i_1}, e_{i_2}), e_n \Bigr \rangle.}
\ee
From \eqref{h1_part2}, it is clear that $h^{(1),n}_\lambda = 0$ for all $n > 2m$. Note also that it follows from the nonlinear  interaction laws \eqref{nonlinear_interaction} that
\bes
\Bigl \langle  B(e_{i}, e_{n -i}), e_n \Bigr \rangle + \Bigl \langle  B(e_{n-i}, e_{i}), e_n \Bigr \rangle = - n \alpha,
\ees
where $\alpha = \frac{\gamma \pi}{\sqrt{2}l^{3/2}}$. By using this identity, we can rewrite $h^{(1),n}_\lambda$ for $n = m+1, \cdots, 2m$ as follows:

\bea \label{h2_part2_2}
h^{(1),n}_\lambda(\xi) = \begin{cases}
{\displaystyle - n \alpha \sum_{i = n - m}^{(n-1)/2} \frac{\xi_{i} \xi_{n-i}}{\beta_{i}(\lambda) + \beta_{n-i}(\lambda) - \beta_n(\lambda)}}, & \text{if $n$ is odd}, \vspace{2ex} \\
{\displaystyle - \frac{n \alpha}{2} \biggl( \sum_{i = n - m}^{(n-2)/2} \frac{2 \xi_{i} \xi_{n-i}}{\beta_{i}(\lambda) + \beta_{n-i}(\lambda) - \beta_n(\lambda)}  +  \frac{(\xi_{n/2})^2}{2\beta_{n/2}(\lambda) - \beta_n(\lambda)}\biggr) } , & \text{if $n$ is even}. \\
\end{cases}
\eea
{\mk where the convention that the sum is zero when the lower bound of the summation index is greater than its upper bound, has been adopted}.

Let us denote by $M$ the matrix whose components are given by 
\be \label{M_local}
M(i,j) := \langle {\hl \chi_\Omega} e_i, e_j \rangle,\qquad 1 \le i,j \le m.
\ee
{\mk Let us also introduce}
\be \label{vector_v}
v_{R}(t):=M^{\mathrm{tr}} \ur(t).
\ee

{\mk By rewriting the reduced system \eqref{reduced_Burgers_local} as
\bea \label{reduced_local_v2}
& \frac{\d z_i}{\d t} = \beta_i(\lambda) z_i + \Bigl \langle B\Bigl(z + h^{(1)}_\lambda(z), z + h^{(1)}_\lambda(z)\Bigr), e_i \Bigr \rangle + v_{R,i}(t),  \qquad {\hl t \in (0,T]},  \quad i = 1, \cdots, m, 
\eea
and by using the {\hl expansions} 
\bes
z = \sum_{i=1}^m z_i e_i,  \qquad h^{(1)}_\lambda(z) = \sum_{n = m+1}^{2m} h^{(1),n}_\lambda(z)e_n,
\ees 
along with the nonlinear  interaction relations} \eqref{nonlinear_interaction}, the above system of equations {\mk becomes}:
\begin{equation}  \label{eqn:z_local}
\boxed{
\begin{aligned}
\frac{\d z_i}{\d t} & = \beta_i(\lambda) z_i + 
\overbrace{i \alpha \Bigl( - \sum_{j = 1}^{\lfloor i/2 \rfloor} \omega_{i,j} z_j z_{i-j} +  \sum_{j = i+1}^m z_j z_{j-i} \Bigr)}^{\mathrm{\bf (a)}}
+ \overbrace{ i \alpha \sum_{j = m-i+1}^{m} z_j h_{\lambda}^{(1),j+i}(z)}^{\mathrm{\bf (b)}}  \\
& + \underbrace{i \alpha \sum_{n = m+1}^{2m-i}  h_{\lambda}^{(1),n}(z)   h_{\lambda}^{(1),n+i}(z)}_{\mathrm{\bf (c)}} + v_{R,i}(t), \qquad {\hl t \in (0,T]},  \quad i = 1, \cdots, m,   
\end{aligned}
}
\end{equation}
where $\lfloor x \rfloor$ denotes the largest integer less than $x$;  $h_{\lambda}^{(1),n}$ is provided by \eqref{h2_part2_2}; and the {\mk coefficients $\omega_{i,j}$ are given} by
\bes
\omega_{i,j} := \begin{cases}
1, & \text{if $i$ is odd, or if $i$ is even and $j \neq i/2$,}  \\
1/2, & \text{if $i$ is even and $j  = i/2$}.
\end{cases}
\ees
In the above system, the {\mk terms gathered in $\mathrm{\bf (a)}$ correspond to the self-interactions between the low modes:} $\langle B(z, z), e_i \rangle$,  the {\mk terms gathered in $\mathrm{\bf (b)}$ correspond to the cross-interactions  between the low and (unresolved) high modes such as parameterized by  $h_{\lambda}^{(1)}$:} $\langle B(z, h^{(1)}_\lambda(z)), e_i \rangle + \langle B(h^{(1)}_\lambda(z),z), e_i \rangle$,   and the {\mk terms gathered in $\mathrm{\bf (c)}$ correspond to the self-interactions between the  high modes (still such as parameterized by $h_{\lambda}^{(1)}$) as projected onto $\mathcal{H}^{\c}$: $\langle B(h^{(1)}_\lambda(z), h^{(1)}_\lambda(z)), e_i \rangle$.}

 Note that in the case $m=2$ the system \eqref{eqn:z_local} takes the same functional form as the $h^{(1)}_\lambda$-based reduced system \eqref{eq:Burgers reduced} derived in Section~\ref{ss:reduction-h1} for the globally distributed control case, {\mk only the matrices given in \eqref{M} and \eqref{M_local} differ.} We refer again to Appendix~\ref{Sect_energy_est} for an analysis of the Cauchy problem associated with \eqref{eqn:z_local}, leaving to the interested reader the generalization to the $m$-dimensional case. 

\subsection{Synthesis of $m$-dimensional locally distributed suboptimal controllers}  \label{ss:BVP_local}

We apply {\mk once more} the Pontryagin maximum principle to derive boundary value problems to be satisfied {\mk by an $h^{(1)}_\lambda$-based}  suboptimal controller. We focus again on the case with terminal payoff given by \eqref{LRBP}, and indicate necessary changes for the case of tracking type \eqref{LRBP'} at the end of this subsection.  

Let us denote the RHS of \eqref{eqn:z_local} by $f(z, v_R)$. The Hamiltonian associated with the cost functional \eqref{JJ_reduced_terminal} reads then as follows:
\be \label{H_local}
H(z, p, \ur)  := \frac{1}{2} \|z + h^{(1)}_\lambda(z)\|^2  +  \frac{\mu_1}{2} \|\ur\|^2  + p^{\mathrm{tr}} f(z, v_R)
\ee
where $p:= (p_1, \cdots, p_m)^{\mathrm{tr}}$ is the costate, and $v_{R}=M^{\mathrm{tr}} \ur$; see \eqref{vector_v}.

Recall also that the terminal payoff, denoted by $C_T(z(T), P_{\c}Y)$, reads {\hl in this case}:
\be \label{C_T}
C_T(z(T), P_{\c} Y):= \frac{\mu_2}{2} \sum_{i=1}^m |z_i(T) - Y_i|^2.
\ee
It follows from the Pontryagin maximum principle that for a given pair
\bes
(z_R^\ast, v_R^\ast)  \in L^2(0,T; \mathcal{H}^{\c}) \times L^2(0,T; \mathcal{H}^{\c}) 
\ees
to be optimal for the {\hl reduced} problem \eqref{LRBP}, it must satisfy the following conditions for all $i = 1, \cdots, m$  (see {\it e.g.} \cite[Chap.~5]{Kirk12}):
\begin{subequations}  \label{Pontryagin relation_local}
\begin{align}
& \frac{\d z^\ast_{R}}{\d t}  = \nabla_{p}H(z^\ast_{R}, p^\ast_{R}, {\mk v^\ast_{R}}) = f(z^\ast_{R}, v^\ast_R), \label{Pontryagin_local-a} \\ 
& \frac{\d p^\ast_{R}}{\d t} =  - \nabla_{z}H(z^\ast_{R}, p^\ast_{R}, {\mk v^\ast_{R}})=  g(z^\ast_{R}, p^\ast_{R}), \\
& \nabla_{u_R} H(z^\ast_{R}, p^\ast_{R}, {\mk v^\ast_{R}})   = 0, \label{Pontryagin_local-c}  \\
 &  p_{R}^\ast(T) =  \nabla_z  C_T(z_R^\ast(T), P_{\c}Y),    \label{Pontryagin_local-d}
\end{align}
\end{subequations}
where $v_{R}^\ast=M^{\mathrm{tr}} u_R^\ast$; $p_R^\ast =\sum_{i=1}^m p_{R,i}^\ast e_i$  {\mk denotes} the costate associated with $z_R^\ast$; and the vector field $(g_1, \cdots, g_m)^{\mathrm{tr}}$ is defined by 
\bea  \label{eq:g2}
g_i(z, p) & := - {\mk \frac{\partial H (z,p,v_R)}{\partial z_i}}= - z_i - \sum_{n=m+1}^{2m} h_\lambda^{(1),n}(z)\frac{\partial h_\lambda^{(1),n}(z)}{\partial z_i} - \sum_{j = 1}^m p_j \frac{\partial f_j(z,v_R)}{\partial z_i}, \qquad i = 1, \cdots, m.
\eea
Here the partial derivatives $\frac{\partial h_\lambda^{(1),n}(z)}{\partial z_i}$  can be obtained by using {\hl the expression of $ h_\lambda^{(1),n}$ given in \eqref{h2_part2_2} which leads to}
\be
\frac{\partial h_\lambda^{(1),n}(z)}{\partial z_i} = \begin{cases}
\frac{\displaystyle -j\alpha z_{n-i}}{\displaystyle \beta_i(\lambda) + \beta_{n-i}(\lambda) - \beta_n(\lambda)}, & \text{if $n \in \{m+1, \cdots, 2m\}$ and $i \in \{n-m, \cdots, m\}$,} \vspace{1.3ex} \\
0, & \text{otherwise}.
\end{cases}
\ee 
The formula for $\frac{\partial f_j(z,v_R)}{\partial z_i}$ can be obtained by taking the corresponding partial derivative of the RHS of \eqref{eqn:z_local} {\mk form which we obtain after simplifications}
\be
\frac{\partial f_j(z,v_R)}{\partial z_i} = \beta_{j}(\lambda) \delta_{ij} + j \alpha ( I_{j,i}^a + I_{j,i}^b + I_{j,i}^c ),
\ee
where {\mkr $\delta_{ij}$ denotes the Kronecker delta,} and {\mk $I_{j,i}^a$, $I_{j,i}^b$ and $I_{j,i}^c$ are given by}
\be
I^a_{j,i}  = \frac{\partial}{\partial z_i}\Bigl( - \sum_{k = 1}^{\lfloor j/2 \rfloor} \omega_{j,k} z_k z_{j-k} +  \sum_{k = j+1}^m z_k z_{k-j} \Bigr) = \begin{cases}
z_{i-j}, & \text{if $i > j$}, \\
z_{i+j}, & \text{if $i = j$ and $i+j\le m$}, \\
z_{i+j} - z_{j-i}, & \text{if $i < j$ and $i+j\le m$}, \\
-z_{j-i}, & \text{if $i < j$ and $i+j > m$}, \\
0, & \text{otherwise};
\end{cases}
\ee
\be
 I^b_{j,i}  = \frac{\partial}{\partial z_i}\Bigl( \sum_{k = m-j+1}^{m} z_k h_{\lambda}^{(1),k+j}(z) \Bigr) = \begin{cases}\displaystyle
h_\lambda^{(1),i+j}  +  \sum_{k = m-j+1}^{m} z_k \frac{\partial h_{\lambda}^{(1),k+j}(z)}{\partial z_i}, & \text{if $i + j > m$}, \vspace{1.5ex}\\
\displaystyle \sum_{k = m-j+1}^{m} z_k \frac{\partial h_{\lambda}^{(1),k+j}(z)}{\partial z_i}, & \text{if $i + j \le m$};
\end{cases}
\ee
and
\be
\displaystyle I^c_{j,i}  = \frac{\partial}{\partial z_i}\biggl( \sum_{n = m+1}^{2m-j}  h_{\lambda}^{(1),n}(z) h_{\lambda}^{(1),n+j}(z) \biggr) = 
 \sum_{n = m+1}^{2m-j}  \biggl ( \frac{\partial h_{\lambda}^{(1),n}(z)}{\partial z_i} h_{\lambda}^{(1),n+j}(z)  +   
h_{\lambda}^{(1),n}(z) \frac{\partial h_{\lambda}^{(1),n+j}(z)}{\partial z_i} \biggr).
\ee

We derive next a relation between $\ur^\ast$ and $p_R^\ast$, which when used in \eqref{Pontryagin relation_local} leads to a BVP for $(z_R^\ast, p_R^\ast)$ to be solved in order to find $\ur^\ast$. To this end, note {\mk {\mkr that} from the expression of the Hamiltonian} $H$ given {\mk  in \eqref{H_local},  we obtain  the following expression of $\nabla_{u_R} H(z^\ast_{R}, p^\ast_{R}, u^\ast_{R})$, which written component-wise,  gives:}
\bes
\frac{\partial H}{\partial u_{R,i}}(z^\ast_{R}, p^\ast_{R}, u^\ast_{R}) = \mu_1 \urc{i}^\ast + \sum_{j=1}^m p^\ast_{R,j} \frac{\partial f_j}{\partial u_{R,i}}(z_R^\ast, M^{\mathrm{tr}}\ur^\ast)  = \mu_1 \urc{i}^\ast + \sum_{j=1}^m p^\ast_{R,j} M(i,j), \qquad  i \in \{1, \cdots, m\}.
\ees
{\mk The first-order optimality condition \eqref{Pontryagin_local-c} leads to}
\be \label{opt_u_sect7}
\ur^\ast = - \frac{1}{\mu_1} M p_R^\ast,
\ee
{\hl where $M$ is given by \eqref{M_local}}.

It follows then that the {\mk controller} $v_R^\ast$ in \eqref{Pontryagin_local-a} takes the form:
\be \label{opt_v_sect7}
v^\ast_R = M^{\mathrm{tr}} \ur^\ast = - \frac{1}{\mu_1} M^{\mathrm{tr}} M p_R^\ast.
\ee

To summarize, {\mk corresponding to  the $h^{(1)}_\lambda$-based} reduced optimal control problem \eqref{LRBP}, we have derived the following BVP {\mk to be satisfied by the optimal trajectory $z_R^\ast$ and its costate $ p_R^\ast$}:
\begin{subequations}  \label{BVP_local}
\begin{align} 
& \frac{\d z^\ast_{R,i}}{\d t} = f_i \Bigl( z^\ast_{R}, v^\ast_R \Bigr), \qquad {\hl t \in (0, T]}, \label{BVP_local-a} \\
 & \frac{\d p^\ast_{R,i}}{\d t} = g_i(z^\ast_{R}, p^\ast_{R}),  \qquad {\hl t \in (0, T]},  \label{BVP_local-b}  \\
& z^\ast_{R,i}(0) = y_{0,i}, \quad p^\ast_{R,i}(T) = \mu_2 (z^\ast_{R_i}(T) - Y_{i}),  \qquad i = 1, \cdots, m,
\end{align}
\end{subequations}
where  $v^\ast_R$ is given by \eqref{opt_v_sect7}, {\hl $y_{0,i}$ is the projection of the initial data $y_0$ for the underlying PDE \eqref{eq:Burgers} against $e_i$,} and the boundary condition for $p^\ast_{R}$ is derived from {\hl the terminal condition} \eqref{Pontryagin_local-d} by using the expression of the terminal payoff $C_T$ given in \eqref{C_T}. Once \eqref{BVP_local} is solved, {\mk the $m$-dimensional controller $u_R^\ast$ given by \eqref{opt_u_sect7} constitutes our $h^{(1)}_\lambda$-based suboptimal controller for the optimal control problem \eqref{LBP}}. {\mk Note that $u_R^\ast$ synthesized this way turns out to be the unique optimal controller for the reduced problem \eqref{LRBP} for the same reasons pointed out in Section~\ref{ss:PMP}.}

The corresponding BVP associated with the reduced optimal control problem \eqref{LRBP'} can be derived in the same fashion; and we indicate below the necessary changes. In this case, the Hamiltonian associated with the cost functional \eqref{JJ_reduced_tracking} reads:
\be \label{H_tracking}
\widetilde{H}(z, p, \ur)  := \frac{1}{2} \|z + h^{(1)}_\lambda(z) - Y\|^2  +  \frac{\mu_1}{2} \|\ur\|^2  + p^{\mathrm{tr}} f(z, v_R).
\ee
The resulting BVP {\mk reads}:
\begin{subequations}  \label{BVP_tracking}
\begin{align} 
& \frac{\d z^\ast_{R,i}}{\d t} = f_i \Bigl( z^\ast_{R}, v^\ast_R  \Bigr), \qquad {\hl t \in (0, T]},  \label{BVP_tracking-a} \\
 & \frac{\d p^\ast_{R,i}}{\d t} = \widetilde{g}_i(z^\ast_{R}, p^\ast_{R}), \qquad {\hl t \in (0, T]}, \label{BVP_tracking-b}  \\
& z^\ast_{R,i}(0) = y_{0,i}, \quad p^\ast_{R,i}(T) = 0,  \qquad i = 1, \cdots, m,
\end{align}
\end{subequations}
where $f(z, v_R)$ denotes  the RHS of \eqref{eqn:z_local}, $v^\ast_R$ is {\hl still} given by \eqref{opt_v_sect7}, {\hl but} in contrast to $g_i$ given by \eqref{eq:g}, the {\mk components $\widetilde{g}_i$ of the} vector field {\mk involved in the RHS of the $p$-equations of \eqref{BVP_tracking}, are now} given by
\bea  \label{eq:g'}
\widetilde{g}_i(z, p) & := - \frac{\partial \widetilde{H}}{\partial z_{i}} = - (z_i - Y_i)  - \sum_{n=m+1}^{2m} (h_\lambda^{(1),n}(z) - Y_n) \frac{\partial h_\lambda^{(1),n}(z)}{\partial z_i} - \sum_{j = 1}^m p_j \frac{\partial f_j(z,v_R)}{\partial z_i}, \qquad i = 1, \cdots, m.
\eea

{\hl Once the above BVP \eqref{BVP_tracking} is solved}, we take $u_R^\ast$ given by \eqref{opt_u_sect7} with $p^\ast_R$ obtained from \eqref{BVP_tracking} as the {\mk $h^{(1)}_\lambda$-based} suboptimal controller for the optimal control problem \eqref{LBP'}.

\subsection{ Control performances: Numerical results}\label{Sec_local_num}

To assess the {\mk ability of the $h^{(1)}_\lambda$-based} reduced optimal control problems \eqref{LRBP} {\mk and} \eqref{LRBP'} in synthesizing suboptimal controllers {\mk of good performance} for respectively the optimal control problems \eqref{LBP} and \eqref{LBP'}, we consider the case where the characteristic function $\chi_\Omega$ is supported on the subdomain $\Omega = [0.2l, 0.8l]$, and the target is taken to {\mk be the target $Y$ used in \eqref{Y_target} for the experiments of  Section \ref{Sec_numresults}. As pointed out prior to Section \ref{ss:derivation_local}, to achieve performances comparable to those  achieved in Section  \ref{Sec_numresults},  it turned out that four-dimensional $h^{(1)}_\lambda$-based reduced systems  were required for the design of suboptimal controllers,  instead of the two-dimensional reduced systems of Section  \ref{Sec_numresults}. As explained above, this increase of the dimension of the resolved subspace $\mathcal{H}^{\c}$ results from the spatial localization of the controller dealt with here.}

\begin{figure}[!hbtp]
\centering
\includegraphics[height=0.4\textwidth, width=.9\textwidth]{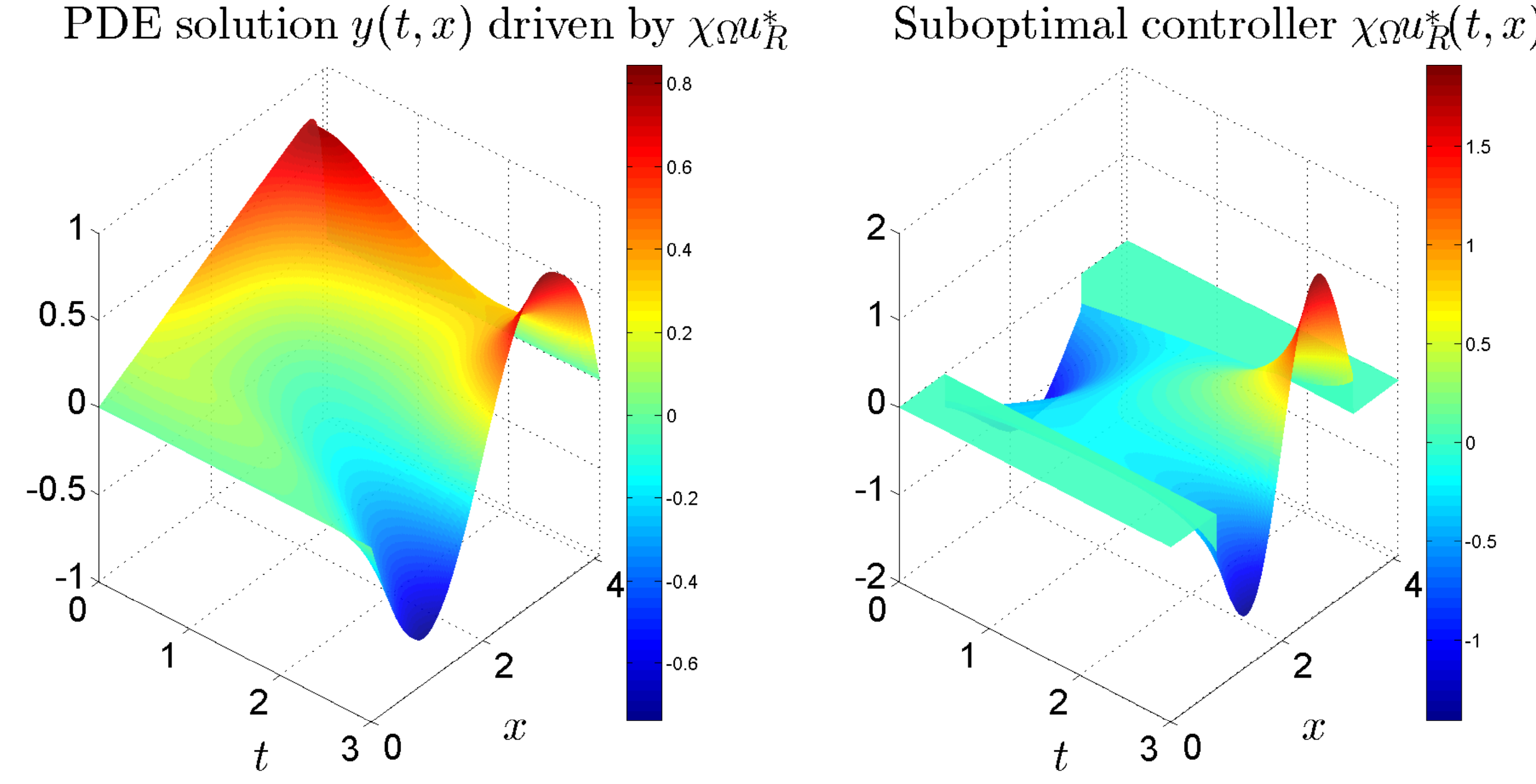}
\vspace{-0.5em}
\caption{{\footnotesize {\mk {\bf Left panel}: The PDE solution field driven by the suboptimal controller $u_R^\ast$ synthesized by solving the $h^{(1)}_\lambda$-based reduced problem \eqref{LRBP}. {\bf Right panel}: The suboptimal  controller $u_R^\ast$ subject to the action of $\chi_\Omega$.
The support of the characteristic function $\chi_\Omega$ is taken to be $\Omega = [0.2l, 0.8l]$.}  {\mk $\mathcal{H}^\c$ is  taken to be spanned by the first four leading eigenmodes ($m=4$)}; the target $Y$ is given by \eqref{Y_target}; and the initial datum is $0.5y^+$. The parameters are $l = 1.3\pi$, $\lambda = 7\lambda_c$, $\nu = 0.25$, $\gamma = 2.5$, and the final time is $T = 3$. The parameters $\mu_1$ and $\mu_2$ in the cost functional \eqref{JJ_terminal} are taken to be $\mu_1 = 1$ and $\mu_2 = 20$.
}} \label{fig:soln_fields_terminal}
\end{figure}

\begin{figure}[!hbtp]
\centering
\includegraphics[height=0.25\textwidth, width=0.75\textwidth]{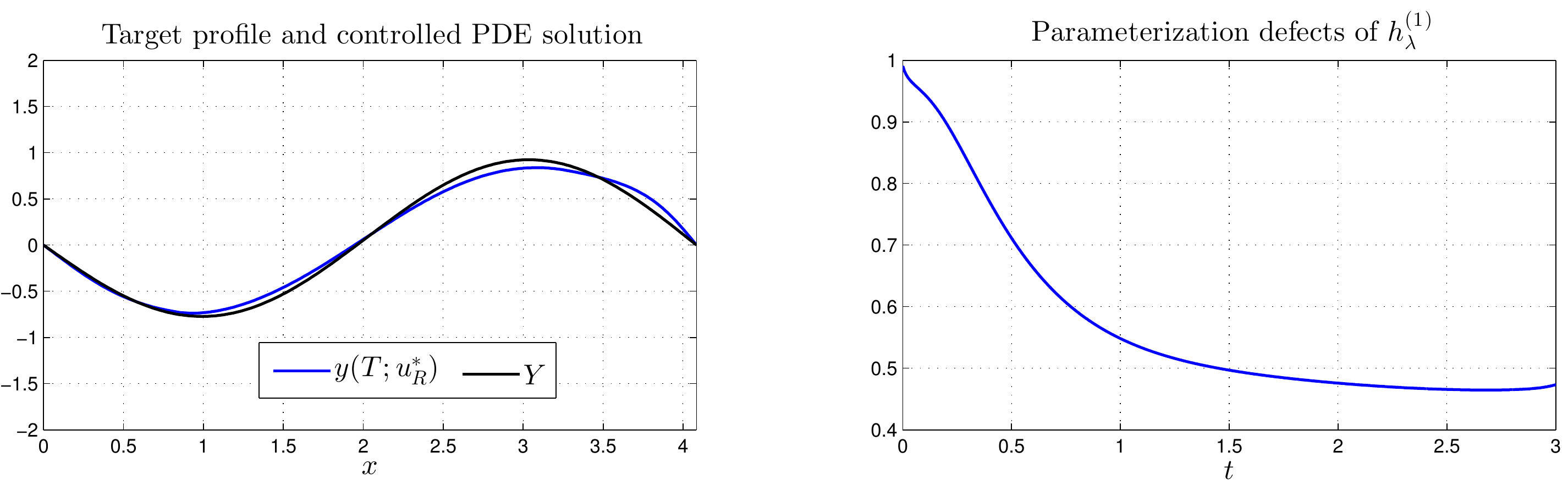}
\vspace{-0.5em}
\caption{{\footnotesize Final time solution profile of the PDE driven by $\chi_{\Omega} u_R^\ast$ compared with  the target $Y$ is given by \eqref{Y_target} ({\bf left panel}); and the parameterization defects associated with the finite-horizon PM, $h^{(1)}_\lambda$, given by \eqref{h1_part2}  for $m=4$ ({\bf right panel}). Parameters are the same as in Fig.~\ref{fig:soln_fields_terminal}.}} \label{fig:yT_terminal}
\end{figure}

{\mk Figures \ref{fig:soln_fields_terminal} and \ref{fig:yT_terminal}  show the performances achieved by the resulting four-dimensional $h^{(1)}_\lambda$-based  suboptimal controllers, corresponding to the cost functional of terminal-payoff type \eqref{JJ_terminal}.}
 {\mk The left panel of Fig.~\ref{fig:soln_fields_terminal} shows the PDE solution field driven by the corresponding suboptimal controller field shown on the right  panel of the same figure.  The left panel of Fig.~\ref{fig:yT_terminal} shows the final-time solution profile, while the right panel shows the corresponding parameterization defect associated with $h^{(1)}_\lambda$}. The {\mk corresponding cost value and relative $L^2$-error of the final time solution profile compared with the target} are given by
\bes
J^{\mathrm{TP}}(y(\cdot; y_0, u_R^\ast), u_R^\ast) = 1.49, \qquad \frac{\|y(T; y_0, u_R^\ast) - Y\|}{\|Y\|} = 9.52\%.
\ees

{\mk As a comparison, by using an $m$-dimensional  Galerkin-based reduced system with $m=16$ to design suboptimal solutions to \eqref{LBP}}, the corresponding cost value and  {\mk relative $L^2$-error} are given by
\bes
J^{\mathrm{TP}}(y(\cdot; y_0, \widetilde{u}_G^\ast), \widetilde{u}_G^\ast) = 1.37, \qquad \frac{\|y(T; y_0, \widetilde{u}_G^\ast) - Y\|}{\|Y\|} = 6.68\%.
\ees

The above numerical results indicate {\mk thus} that the $4$-dimensional $h^{(1)}_\lambda$-based reduced problem \eqref{LRBP} can be used to {\mk design a very good suboptimal controller (for the prescribed target $Y$ given by \eqref{Y_target})} for the optimal control problem \eqref{LBP} {\mk with performance comparable to the (more standard) higher-dimensional Galerkin-based reduced systems}. {\mk This success goes with the relatively small parameterization defect  as well as with the relatively small energy kept in the high-modes (not shown); see right panel of Fig.~\ref{fig:yT_terminal}.} {\mk Note that for these experiments}, the system parameters are chosen to be $l = 1.3\pi$, $\lambda = 7\lambda_c$, $\nu = 0.25$, $\gamma = 2.5$, {\mk while} the final time {\mk is taken to be} $T = 3$. The parameters $\mu_1$ and $\mu_2$ in the cost functional \eqref{JJ_terminal} are taken to be $\mu_1 = 1$ and $\mu_2 = 20$. The initial datum is a scaled version of the corresponding positive steady state $y^+$ of the uncontrolled PDE, namely $y_0 = 0.5 y^+$. 

\begin{figure}[!hbtp]
\centering
\includegraphics[height=0.4\textwidth, width=0.9\textwidth]{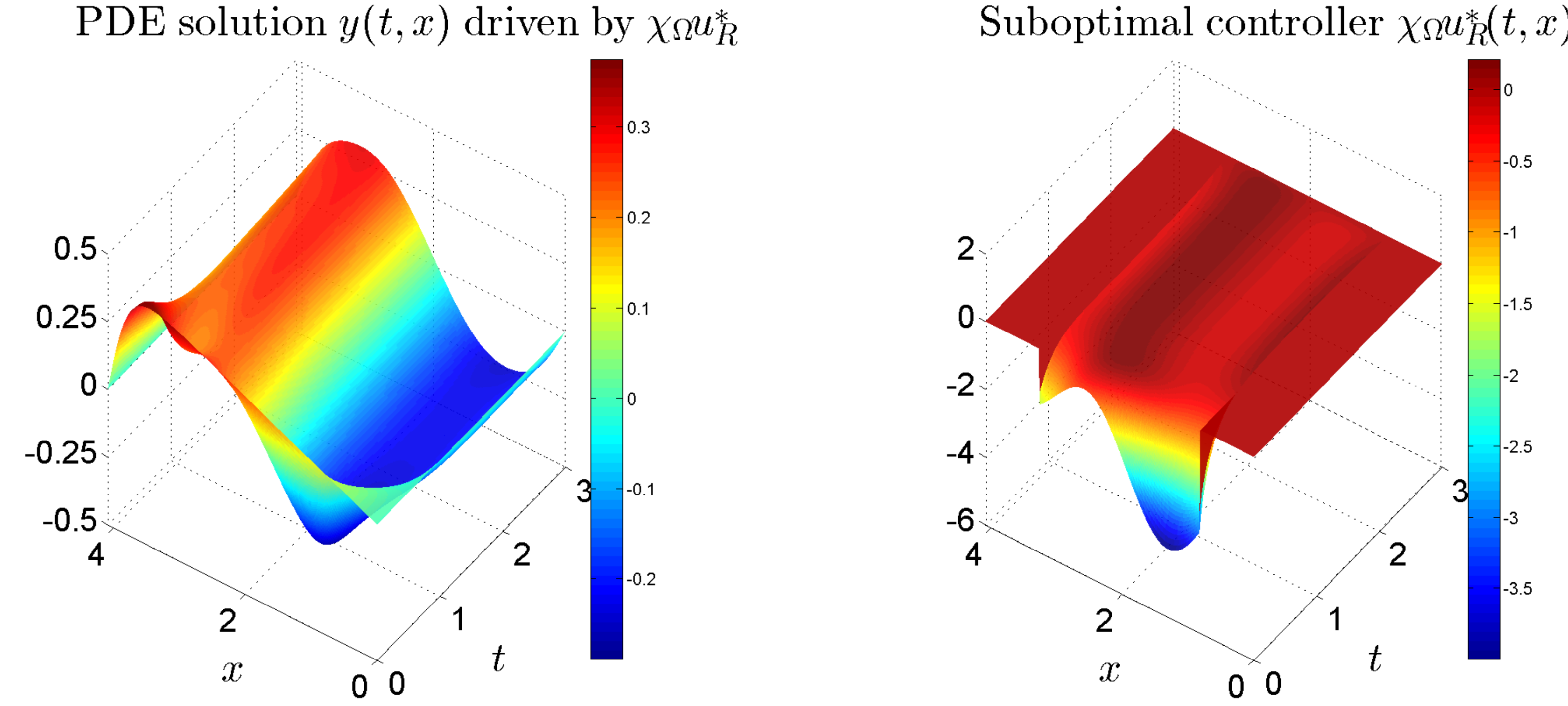}
\vspace{-0.5em}
\caption{{\footnotesize {\mk {\bf Left panel}: The PDE solution field driven by the suboptimal controller $u_R^\ast$ synthesized by solving the $h^{(1)}_\lambda$-based reduced problem \eqref{LRBP'}. {\bf Right panel}: The suboptimal  controller subject to the action of $\chi_\Omega$.
The support of the characteristic function $\chi_\Omega$ is taken to be $\Omega = [0.2l, 0.8l]$.The resolved modes are taken to be the first four leading eigenmodes ($m=4$), the target is $Y = -0.1\langle y^-, e_1\rangle e_1 + 1.6 \langle y^-, e_2 \rangle e_2$ and the initial datum is $0.8 y^+$. The parameters are $l = 1.3\pi$, $\lambda = 3\lambda_c$, $\nu = 0.2$, $\gamma = 2.5$, and the final time is $T = 3$. The parameter $\mu_1$ in the cost functional \eqref{JJ_tracking} is taken to be $\mu_1 = 0.02$. } 
}} \label{fig:soln_fields_tracking}
\end{figure}

\begin{figure}[!hbtp]
\centering
\includegraphics[height=0.25\textwidth, width=0.75\textwidth]{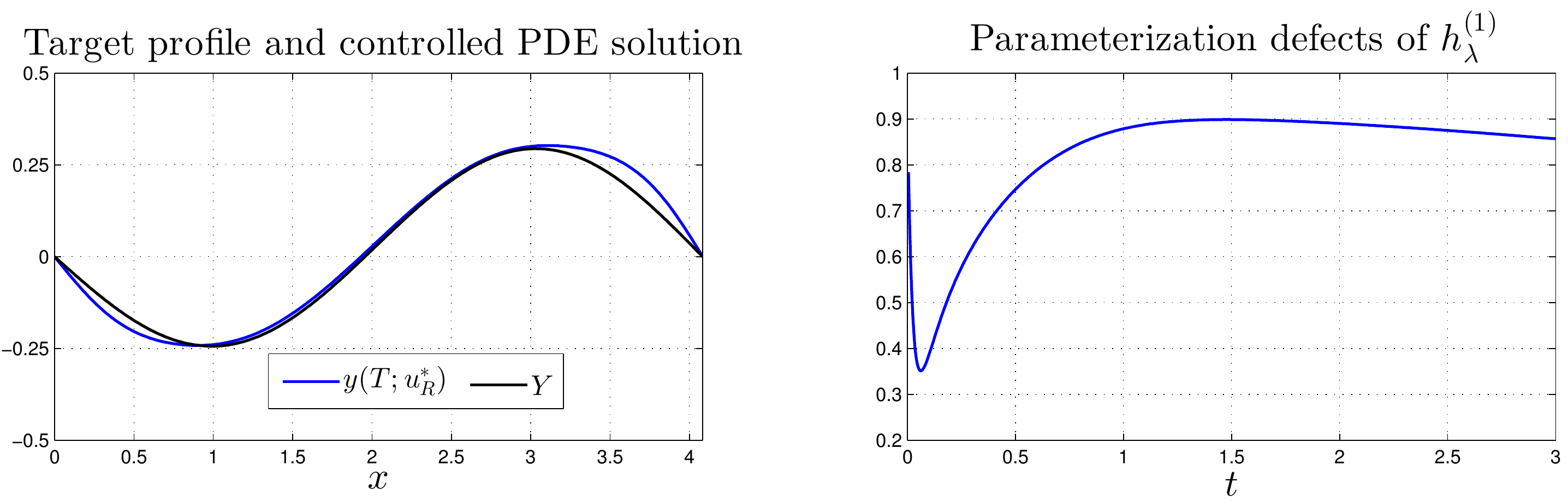}
\vspace{-0.5em}
\caption{{\footnotesize Final time solution profile of the PDE driven by $\chi_{\Omega}u_R^\ast$ compared with  the target $Y$ is given by \eqref{Y_target} ({\bf left panel}); and the parameterization defects associated with the finite-horizon PM, $h^{(1)}_\lambda$,  given by \eqref{h1_part2}  for $m=4$ ({\bf right panel}). Parameters are the same as in Fig.~\ref{fig:soln_fields_tracking}.}} \label{fig:yT_tracking}
\end{figure}

{\mk The performances of the $4$-dimensional $h^{(1)}_\lambda$-based suboptimal controller for \eqref{LRBP'} associated with the cost} functional of tracking type \eqref{JJ_tracking}  {\mk are illustrated} in Figs.~\ref{fig:soln_fields_tracking} and \ref{fig:yT_tracking}. {\mk  The experimental conditions are here chosen to be: $l = 1.3\pi$, $\lambda = 3\lambda_c$, $\nu = 0.2$, $\gamma = 2.5$,  while  the final time is still taken to be $T = 3$. The parameter $\mu_1$ in the cost functional \eqref{JJ_tracking} is taken to be $\mu_1 = 0.02$ and the initial datum is  $y_0 = 0.8 y^+$. }

 {\mk For these experiments, the corresponding cost value and  {\mk relative $L^2$-error} are given by}
\bes
J^{\mathrm{track}}(y(\cdot; y_0, u_R^\ast), u_R^\ast) = 0.032, \qquad \frac{\|y(T; y_0, u_R^\ast) - Y\|}{\|Y\|} = 12.32\%.
\ees

For a high-dimensional Galerkin-based reduced problem with $m=16$, the {\mk corresponding cost value} and {\mk relative $L^2$-error}  are given by
\bes
J^{\mathrm{track}}(y(\cdot; y_0, \widetilde{u}_G^\ast), \widetilde{u}_G^\ast) = 0.025, \qquad \frac{\|y(T; y_0, \widetilde{u}_G^\ast) - Y\|}{\|Y\|} = 10.86\%.
\ees

{\mk Here again, a fairly good performance of the suboptimal controller\footnote{for the optimal control \eqref{LBP'}.} as synthesized by solving the $4$-dimensional $h^{(1)}_\lambda$-based reduced problem  \eqref{LRBP'}, is achieved. } {\mk Due to the deterioration of the parameterization defect of $h^{(1)}_\lambda$ that can be observed  by comparing the right panel of Fig.~\ref{fig:yT_tracking} with the right panel of Fig.~\ref{fig:yT_terminal},  the error estimate  \eqref{cor2:goal} suggests that such a success has to come with a noticeable reduction of the energy contained in the high modes of the PDE solution driven by the suboptimal controller synthesized for \eqref{LRBP'}  compared to the PDE solution driven by the suboptimal controller synthesized for \eqref{LRBP}. Such  theoretical prediction based on Corollary \ref{Cor_2} can actually be empirically confirmed  by looking at the numerical values of these high-mode energies (not shown). }

{\mk Finally, it is worth mentioning that similar to the globally distributed case, the performances of the $h^{(1)}_\lambda$-based reduced systems and the associated parameterization defects of $h^{(1)}_\lambda$ depend on the target and the length of the time horizon; cf.~Figs.~\ref{fig:PM_Sect5}, \ref{fig:contour} and \ref{fig:J}.} {\mk The dependence on the PDE initial datum turned out also to be an important factor}. In particular, it has been observed that for both problems \eqref{LBP} and \eqref{LBP'} the parameterization defects deteriorate when the scaling factors $\delta$ used in the construction of the initial datum $y_0 =\delta y^+$ increases. Based on the results of Section~\ref{Sect_Burgers_h2} for the globally distributed case, it can be reasonably expected that PM functions such as $h^{(2)}_\lambda$ {\mk that bring higher-order terms  compared to $h^{(1)}_\lambda$ (cf.~Theorem \ref{THM_h2}) can allow to reach} better performance for a broader range of initial data and target profiles; {\mk the parameterization defects being reasonably expected to get smaller.}

\section*{Acknowledgments}

We are grateful to Monique Chyba and to Bernard Bonnard for their interest in our works on parameterizing manifolds, which led the authors to propose this article. MDC is also grateful to Denis Rousseau and Michael Ghil for the
unique environment they provided to complete  this work, at the CERES-ERTI, \'Ecole Normale Sup\'erieure, Paris. 
This work has been partly supported by the National Science Foundation
 grant DMS-1049253  and Office of Naval Research grant N00014-12-1-0911.

\appendix

\section{Suboptimal Controller Synthesis Based on Galerkin Projections and Pontryagin Maximum Principle}

To assess the performance of the PM-based reduced systems considered in Sections~\ref{Sect_Burgers} and \ref{Sect_Burgers_h2} in synthesizing suboptimal controllers in the context of a Burgers-type equation, we derive in this appendix suboptimal control  problems associated with the globally distributed optimal control problem \eqref{BP} based on Galerkin approximations. Section~\ref{ss:2D-Galerkin} concerns a two-mode Galerkin approximation; and Section~\ref{ss:ND-Galerkin} deals with the more general {\hl $m$-dimensional} case. The former serves as a basis of comparison to {\hl analyze the performance achieved by the PM-based approach}, while the latter can in principle provide a good indication of the true optimal controller of the {\hl underlying} optimal control problems by taking the dimension sufficiently large. {\hl Results for the general {\hl $m$-dimensional} case will also be used in Section~\ref{Sect_Burgers_local} to derive Galerkin-based  reduced systems for the locally distributed problems \eqref{LBP} and \eqref{LBP'}.}

\subsection{{\mk Suboptimal controller based on a 2D Galerkin reduced optimal problem}}  \label{ss:2D-Galerkin}

{\hl We first present the} reduced optimal control problem based on a two-mode Galerkin approximation of the underlying PDE \eqref{eq:Burgers}, which can be derived by simply setting $h^{(1)}_\lambda$ in \eqref{reduced_Burgers}--\eqref{J2_Burgers} to zero. The corresponding operational forms for the cost functional and reduced system for the low modes can be obtained from \eqref{J2_Burgers-c}--\eqref{eq:Burgers reduced} by setting $\alpha_1(\lambda)$ and $\alpha_2(\lambda)$ to be zero. The resulting cost functional reads:
\be \label{J2_Galerkin}
J_G(v, u_G) = \int_0^T \bigl[ \mathcal{G}^G(v(t))  + \mathcal{E}(u_G(t)) \bigr] \d t +  C_T(v(T), P_{\c}Y),
\ee
where $v = v_{1} e_1 + v_{2} e_2 \in L^2(0,T; \mathcal{H}^{\c})$ is the state variable, $u_G = u_{G,1} e_1 + u_{G,2} e_2 \in L^2(0,T; \mathcal{H}^{\c})$ is the control, $C_T$ is the terminal payoff term defined by \eqref{C_T_sect5}, and
\bea
\mathcal{G}^G(v ) &:= \frac{1}{2}\|v \|^2 = \frac{1}{2} [(v_1 )^2 + (v_2 )^2], \quad \mathcal{E}(u_G) := \frac{\mu_1}{2}\|u_G\|^2 = \frac{\mu_1}{2}[(u_{G,1} )^2 + (u_{G,2})^2].
\eea
The equations for $v_1$ and $v_2$ are given by:
\bea \label{eq:Burgers Galerkin}
& \frac{\d v_1}{\d t} = \beta_1(\lambda) v_1 + \alpha v_1 v_2  + a_{11} u_{G,1}(t) + a_{21} u_{G,2}(t), \\
&  \frac{\d v_2}{\d t} = \beta_2(\lambda) v_2 - \alpha (v_1)^2   + a_{12} u_{G,1}(t) + a_{22} u_{G,2}(t),
\eea
which is subjected to the initial conditions:
\bea \label{bdry:Galerkin}
v_1(0) = \langle y_0, e_1 \rangle,   \qquad v_2(0) = \langle y_0, e_2 \rangle,  
\eea
where  $\alpha = \frac{\gamma \pi}{\sqrt{2}l^{3/2}}$.

The corresponding Galerkin-based reduced optimal control problem for \eqref{BP} reads:
\be  \label{GBP}  
\begin{aligned}
 \hspace{-1em} \min \, J_G(v, u_G)  \; \text{ s.t. } \; (v, u_G) \in L^2(0,T; \mathcal{H}^{\c}) \times L^2(0,T; \mathcal{H}^{\c})  \text{ \ solves   \ \eqref{eq:Burgers Galerkin}--\eqref{bdry:Galerkin}}.
\end{aligned}
\ee

It follows again from the Pontryagin maximum principle that for a given pair
\bes
(v_G^\ast, u_G^\ast)  \in L^2(0,T; \mathcal{H}^{\c}) \times L^2(0,T; \mathcal{H}^{\c})  
\ees
to be optimal for the problem \eqref{GBP}, it must satisfy the following conditions:
\begin{subequations}  \label{Pontryagin relation Galerkin}
\begin{align} 
& \frac{\d v_{G,1}^\ast}{\d t} = \beta_1(\lambda) v_{G,1}^\ast + \alpha v_{G,1}^\ast v_{G,2}^\ast  + a_{11} u_{G,1}^\ast(t) + a_{21}u_{G,2}^\ast(t), \label{Pontryagin-a1} \\
&  \frac{\d v_{G,2}^\ast}{\d t}  = \beta_2(\lambda) v_{G,2}^\ast - \alpha (v_{G,1}^\ast)^2   + + a_{12} u_{G,1}^\ast(t) + a_{22}u_{G,2}^\ast(t),  \label{Pontryagin-a2} \\
& \frac{\d p_{G,1}^\ast}{\d t} =  - v_{G,1}^\ast    - \beta_1(\lambda) p_{G,1}^\ast  -  \alpha  p_{G,1}^\ast v_{G,2}^\ast + 2 \alpha p_{G,2}^\ast v_{G,1}^\ast,  \label{Pontryagin-b1}\\
& \frac{\d p_{G,2}^\ast}{\d t} =  - v_{G,2}^\ast   - \beta_2(\lambda) p_{G,2}^\ast   -  \alpha  p_{G,1}^\ast  v_{G,1}^\ast,  \label{Pontryagin-b2}  \\
& (u_{G,1}^\ast, u_{G,2}^\ast)^{\mathrm{tr}} = - \Bigl( \frac{a_{11} p_{G,1}^\ast(t) + a_{12} p_{G,2}^\ast(t)}{\mu_1},  \frac{a_{21} p_{G,1}^\ast(t) +  a_{22} p_{G,2}^\ast(t)}{\mu_1} \Bigr )^{\mathrm{tr}} = - \frac{1}{\mu_1} M^{\mathrm{tr}}p_G^\ast,  \label{Pontryagin-c2}
\end{align}
\end{subequations}
where $v_{G,1}^\ast = \langle v_G^\ast, e_i \rangle$, $u_{G,i}^\ast = \langle u_G^\ast,  e_i \rangle$, $i = 1, 2$, and $p_G^\ast = p_{G,1}^\ast e_1 + p_{G,2}^\ast e_2$ denotes the costate associated with $v_G^\ast$.

Thanks to \eqref{Pontryagin-c2}, we can express the controller $u_{G,i}^\ast$ in \eqref{Pontryagin-a1}--\eqref{Pontryagin-a2} in terms of the costate $p_{G,i}^\ast$, leading thus to the following BVP for $v_G^\ast$ and $p_G^\ast$: 
\bea \label{BVP-Galerkin}
\frac{\d v_1}{\d t} &= \beta_1(\lambda) v_1 + \alpha v_1v_2  + f_3(p_1,p_2), \\
 \frac{\d v_2}{\d t} & = \beta_2(\lambda) v_2 - \alpha (v_1)^2 +  f_4(p_1,p_2), \\
 \frac{\d p_1}{\d t} &=  -2 v_1    - \beta_1(\lambda) p_1   -  \alpha  p_1  v_2  + 2 \alpha p_2  v_1,  \\
\frac{\d p_2}{\d t} & =  -2 v_2   - \beta_2(\lambda) p_2   -  \alpha  p_1  v_1, 
\eea
subject to the boundary condition
\be \label{BVP-Galerkin bc}
v_1(0) = \langle y_0, e_1 \rangle, \quad v_2(0) = \langle y_0, e_2 \rangle, \quad p_1(T) =\mu_2 (v_{1}(T) -  Y_1), \quad p_2(T) = \mu_2 (v_2(T) -  Y_2),
\ee
where $f_3$ and $f_4$ are defined by \eqref{f3-f4}, and the boundary condition for the costate is derived in the same way as in \eqref{p-condition} thanks to the Pontryagin maximum principle.  Once this BVP is solved,  the corresponding controller $u_{G}^\ast$ is determined by \eqref{Pontryagin-c2} which provides the unique optimal controller for the Galerkin-based reduced optimal control problem \eqref{GBP}, due again to the fact that the cost functional \eqref{J2_Galerkin} is quadratic in $u_G$ and the dependence on the controller is affine for the system of equations \eqref{eq:Burgers Galerkin}; see {\it e.g.}~\cite[Sect.~5.3]{Kirk12} and \cite{Trelat2012}. Note also that analogous results to those presented in Lemma~\ref{Lem:M} hold for the reduced optimal control problem \eqref{GBP} as well.

\subsection{{\mk Suboptimal controller based on an $m$-dimensional Galerkin reduced optimal problem}}  \label{ss:ND-Galerkin}

We derive now a more general reduced optimal control problem based on higher-dimensional Galerkin approximation, where the subspace $\mathcal{H}^{\c}$ is taken to be spanned by the first $m$ eigenmodes:
\be  \label{Hc-m}
\mathcal{H}^{\c} := \mathrm{span}\{e_1, \cdots, e_m\}.
\ee 
The main interest is that by choosing $m$ sufficiently large, such a reduced problem can serve in principle to provide a good estimate of the true optimal controllers of the {\hl globally distributed} optimal control problem \eqref{BP}, which can be taken then as a benchmark for the numerical experiments reported in Sections~\ref{Sect_Burgers} and \ref{Sect_Burgers_h2}. {\hl Analogous reduced problems associated with the locally distributed cases \eqref{LBP} and \eqref{LBP'} considered in Section~\ref{Sect_Burgers_local} can be derived in the same way (and actually the corresponding results are the same as those presented in Section~\ref{ss:BVP_local} by setting $h^{(1)}_\lambda$ therein to be zero).} 

The Galerkin-based reduced optimal control problem \eqref{GBP} when generalized to the case with $m$ {\mkr controlled} modes reads:
\be  \label{GBP'} 
\begin{aligned}
 \hspace{-1em} \min \, \widetilde{J}_G(v, \widetilde{u}_G)  \; \text{ s.t. } \; (v, \widetilde{u}_G) \in L^2(0,T; \mathcal{H}^{\c}) \times L^2(0,T; \mathcal{H}^{\c})  \text{ \ solves   \ \eqref{eq:Burgers Galerkin-m}--\eqref{bdry:Galerkin-m}  below},
\end{aligned}
\ee
where $\mathcal{H}^\c$ is the $m$-dimensional {\mkr reduced phase space} defined in \eqref{Hc-m}, and  
\bes
\widetilde{J}_G(v, \widetilde{u}_G) = \int_0^T \bigl[  \frac{1}{2} \sum_{i=1}^m(v_i)^2  + \frac{\mu_1}{2} \sum_{i=1}^m (\widetilde{u}_{G,i} )^2  \bigr] \d t + \frac{\mu_2}{2} \sum_{i=1}^m |v_i(T) - Y_i|^2.
\ees

The system of equations that $v(\cdot; \widetilde{u}_G)$ satisfies is given by:
\bea \label{eq:Burgers Galerkin-m}
 \frac{\d v_i}{\d t} = \beta_i(\lambda) v_i + \Bigl \langle B \Bigl( \sum_{i=1}^m v_i e_i, \sum_{i=1}^m v_i e_i \Bigr), e_i \Bigr \rangle  +   [M^{\mathrm{tr}}\widetilde{u}_{G}(t)]_i, \qquad i = 1,\cdots, m,
\eea
which is subjected to the initial conditions:
\bea \label{bdry:Galerkin-m}
v_i(0) = \langle y_0, e_i \rangle,   \qquad i = 1,\cdots,  m,
\eea
where the matrix $M_{m\times m}$ is the representation of the linear operator $P_{\c}\mathfrak{C}$ under the basis $e_1, \cdots, e_m$, {\it i.e.} the elements of $M$ are given by $a_{ij} = \langle \mathfrak{C}e_i, e_j \rangle$ (see \eqref{M} for the case $m=2$) and $[M^{\mathrm{tr}}\widetilde{u}_{G}(t)]_i$ denotes the $i^{\mathrm{th}}$-component of the vector $M^{\mathrm{tr}}\widetilde{u}_{G}(t)$. 

As before, by using the Pontryagin maximum principle, we can derive the following BVP to be satisfied by any optimal pair $(v_G^\ast, \widetilde{u}^\ast_{G})$ of \eqref{GBP'}: 
\begin{subequations} \label{bvp-Galerkin-m}
\begin{align}
 & \frac{\d v_i}{\d t}  = \beta_i(\lambda) v_i + i \alpha \Bigl( - \sum_{j = 1}^{\lfloor i/2 \rfloor} \omega_{i,j} v_j v_{i-j} +  \sum_{j = i+1}^m v_j v_{j-i} \Bigr) - \frac{1}{\mu_1} [M^{\mathrm{tr}} M p]_i, &&  i = 1, \cdots, m, \label{bvp-Galerkin-m-1} \\
& \frac{\d p_i}{\d t}  = - v_i - \sum_{j = 1}^m p_j \frac{\partial f_j(v, p)}{\partial v_i}, &&  i = 1, \cdots, m, \label{bvp-Galerkin-m-2} \\
&  v_{i}(0)  = y_{0,i}, \qquad  \qquad p_{i}(T) = \mu_2 (v_{i}(T) - Y_{i}),  &&  i = 1, \cdots, m, \label{bvp-Galerkin-m-3}
\end{align}
\end{subequations}
where the optimal controller $\widetilde{u}^\ast_{G}$ is related to the corresponding costate $p_G^\ast$ by
\be
\widetilde{u}^\ast_{G} = - \frac{1}{\mu_1} M p_G^\ast,
\ee
see \eqref{Pontryagin-c2} for the case $m=2$. Here, $f_i$, $i=1,\cdots, m$, denotes the RHS of \eqref{bvp-Galerkin-m-1} and we have used the nonlinear interactions \eqref{nonlinear_interaction} to derive the quadratic parts of $f_i$. The formula for $\frac{\partial f_j(v, p)}{\partial v_i}$ is given by:
\be
\frac{\partial f_j(v,p)}{\partial v_i} = \beta_{j}(\lambda) \delta_{ij} + j \alpha I_{j,i},
\ee
where $\delta_{ij}$ {\mkr denotes the Kronecker delta}, and
\be
I_{j,i}  = \frac{\partial}{\partial v_i}\Bigl( - \sum_{k = 1}^{\lfloor j/2 \rfloor} \omega_{j,k} v_k v_{j-k} +  \sum_{k = j+1}^m v_k v_{k-j} \Bigr) = \begin{cases}
v_{i-j}, & \text{if $i > j$}, \\
v_{i+j}, & \text{if $i = j$ and $i+j\le m$}, \\
v_{i+j} - v_{j-i}, & \text{if $i < j$ and $i+j\le m$}, \\
-v_{j-i}, & \text{if $i < j$ and $i+j > m$}, \\
0, & \text{otherwise};
\end{cases}
\ee
with $\lfloor x \rfloor$ being the largest integer less than $x$  and the coefficients $\omega_{i,j}$ given by
\bes
\omega_{i,j} := \begin{cases}
1, & \text{if $i$ is odd, or if $i$ is even and $j \neq i/2$,}  \\
1/2, & \text{if $i$ is even and $j  = i/2$}.
\end{cases}
\ees

\section{Global Well-posedness for the Two-dimensional $h^{(1)}_\lambda$-based Reduced System \eqref{eq:Burgers reduced}} \label{Sect_energy_est}

In this appendix, we show that for any given initial datum and any fixed $T>0$, the $h^{(1)}_\lambda$-based reduced system \eqref{eq:Burgers reduced} admits a unique {\hl mild} solution in the space $C([0,T]; \mathbb{R}^2)$.\footnote{For any $
T>0$, a given  continuous function $\mathbf{z}: [0, T] \rightarrow \mathbb{R}^2$ is called a mild solution to the reduced system \eqref{eq:Burgers reduced} if it satisfies the corresponding integral form of the system: $\mathbf{z}(t) = \mathbf{z}(0) + \int_0^t  \mathbf{F}(s,\mathbf{z}(s))\, \d s$, for all $t\in [0,T]$, where $\mathbf{z}:=(z_1, z_2)^{\mathrm{tr}}$ and $\mathbf{F}$ denotes the RHS of \eqref{eq:Burgers reduced}.}  The result follows from {\hl classical ODE theory \cite{amann1990}} once we can establish {\it a priori} bounds for the solution $(z_1(t), z_2(t))$.  Similar (but more tedious) estimates can be used to deal with the Cauchy problem associated with the $h^{(2)}_\lambda$-based reduced system \eqref{eq:Burgers reduced-h2} derived in Section~\ref{Sect_Burgers_h2} and the more general $m$-dimensional $h^{(1)}_\lambda$-based reduced system \eqref{eqn:z_local} encountered in Section~\ref{Sect_Burgers_local}.

Let us first recall that the two-dimensional $h^{(1)}_\lambda$-based reduced system is given by:
\begin{subequations} \label{eq:Burgers reduced-recall}
\begin{align}
& \hspace{-0.77em}\frac{\d z_1}{\d t} = \beta_1(\lambda) z_1 + \alpha [ z_1z_2 + \alpha_1(\lambda) z_1z_2^2 + \alpha_1(\lambda) \alpha_2(\lambda) z_1 z_2^3] + a_{11}\urc{1}(t) + a_{21}\urc{2}(t), \label{eq:Burgers reduced_1}\\
& \hspace{-0.77em} \frac{\d z_2}{\d t} = \beta_2(\lambda) z_2 + \alpha[- z_1^2 + 2 \alpha_1(\lambda) z_1^2z_2 + 2 \alpha_2(\lambda) z_2^3] + a_{12}\urc{1}(t) + a_{22}\urc{2}(t), \label{eq:Burgers reduced_2}
\end{align}
\end{subequations}
where $u_{R}(\cdot):=u_{R,1}(\cdot)e_1 + u_{R,2}(\cdot)e_2 \in L^2(0,T; \mathcal{H}^{\c})$ with $T>0$ being the fixed finite horizon, $\alpha_1(\lambda)$ and $\alpha_2(\lambda)$ are defined in \eqref{alpha1-2}, $\alpha = \frac{\gamma \pi}{\sqrt{2}l^{3/2}}$, and $a_{ij}$, $1\le i,j \le 2$, are elements of the coefficients matrix $M$ associated with the operator $\mathfrak{C}$; see \eqref{C-def}--\eqref{M}.

We check below by energy estimates that no finite time blow-up can occur for solutions to the system \eqref{eq:Burgers reduced-recall} emanating from any initial datum $(z_{1,0}, z_{2,0}) \in \mathbb{R}^2$. For this purpose, let us define
\bes
R := \max \Biggl \{|z_{2,0}|, \; \frac{\alpha}{|2\alpha \alpha_2(\lambda)|}, \; \sqrt{\frac{|\beta_2(\lambda)|}{|2\alpha \alpha_2(\lambda)|}} \Biggr\} \qquad \mbox{and} \qquad
C := \int_0^T |a_{12}\urc{1}(t) + a_{22}\urc{2}(t)| \d t.
\ees
We claim that
\be \label{eq:z2_bound}
|z_2(t)| \le e^{C/R}R \qquad \Forall t\in [0, T].
\ee
It is clear that we only need to deal with those values of $t$ such that $|z_2(t)| > R$. Assume that there exists such time instances, otherwise we are done. Let us fix an arbitrary interval $[t_\ast, t^\ast] \subset [0, T]$ such that
\be \label{eq:condition_z2}
|z_2(t)| \ge R \qquad \Forall t \in [t_\ast, t^\ast].
\ee
Since $R \ge |z_{2,0}|$ and $z_2$ depends continuously on $t$, we can reduce $t_\ast$ such that $z_2(t_\ast) = R$ while the condition \eqref{eq:condition_z2} remains true.

Now by multiplying $z_2(t)$ on both sides of \eqref{eq:Burgers reduced_2}, we obtain
\be
\frac{1}{2} \frac{\d [(z_2)^2]}{\d t} = c(t) (z_2)^2, \qquad \Forall t \in [t_\ast, t^\ast],
\ee
where
\bes
c(t) := \Bigl( \beta_2(\lambda) - \frac{\alpha (z_1)^2}{z_2} + 2\alpha \alpha_1(\lambda) (z_1)^2 + 2 \alpha \alpha_2(\lambda) (z_2)^2 + \frac{a_{12}\urc{1}(t) + a_{22}\urc{2}(t)}{z_2}\Bigr).
\ees
It follows then that
\be \label{eq:contra_basis}
[z_2(t^\ast)]^2 = e^{2\int_{t_\ast}^{t^\ast} c(t) \d t}[z_2(t_\ast)]^2.
\ee

Since $|z_2(t)| \ge R$ for all $t \in [t_\ast, t^\ast]$ by the choices of $t_\ast$ and $t^\ast$, we get
\bes
\int_{t_\ast}^{t^\ast} c(t) \, \d t \le \beta_2(\lambda)( t^\ast - t_\ast) + \int_{t_\ast}^{t^\ast} [ \frac{\alpha}{R} + 2\alpha \alpha_1(\lambda)] (z_1)^2 \d t + 2 \alpha \alpha_2(\lambda) R^2 ( t^\ast - t_\ast) + \frac{\int_{t_\ast}^{t^\ast} |a_{12}\urc{1}(t) + a_{22}\urc{2}(t)| \d t}{R},
\ees
where we have used $|- \frac{\alpha}{z_2}| \le \frac{\alpha }{R}$ and $2 \alpha \alpha_2(\lambda) (z_2)^2 \le 2 \alpha \alpha_2(\lambda) R^2$, which follow from the definition of $R$ and the fact that $\alpha >0$ and $\alpha_2(\lambda)<0$.

According again to the definition of $R$ and the facts that $\alpha > 0$, $\alpha_1(\lambda) < 0$ and $\alpha_2(\lambda)<0$, we get
\bes
\frac{\alpha}{R} + 2\alpha \alpha_1(\lambda) \le 0 \qquad \mbox{ and } \qquad \beta_2(\lambda)( t^\ast - t_\ast) + 2 \alpha \alpha_2(\lambda) R^2 ( t^\ast - t_\ast) \le 0.
\ees
We obtain then
\bes
\int_{t_\ast}^{t^\ast} c(t) \d t \le \frac{\int_{t_\ast}^{t^\ast} |a_{12}\urc{1}(t) + a_{22}\urc{2}(t)| \d t}{R} \le \frac{C}{R}.
\ees
By reporting the above estimate in \eqref{eq:contra_basis} and using $|z_2(t_\ast)| = R$, we obtain
\bes
|z_2(t^\ast)| \le e^{C/R}|z_2(t_\ast)| = e^{C/R}R,
\ees
and \eqref{eq:z2_bound} is thus proven.

Note also that by multiplying $z_1(t)$ on both sides of \eqref{eq:Burgers reduced_1}, we obtain for any $t\in [0, T]$ at which $z_1(t)\neq 0$ that
\be \label{eq:z1_energy}
\frac{1}{2} \frac{\d [(z_1)^2]}{\d t} = (z_1)^2 \Bigl( \beta_1(\lambda) + \alpha z_2 + \alpha \alpha_1(\lambda) (z_2)^2 + \alpha \alpha_1(\lambda) \alpha_2(\lambda) (z_2)^3 + \frac{a_{11}\urc{1}(t) + a_{21}\urc{2}(t)}{z_1}\Bigr).
\ee
It follows then from the boundedness of $z_2$ and \eqref{eq:z1_energy} that $z_1$ can grow at most exponentially. Consequently, no finite time blow-up can occur for the $h^{(1)}_\lambda$-based reduced system \eqref{eq:Burgers reduced-recall}.

\end{document}